\documentclass[11pt]{article}
\usepackage[centertags]{amsmath}
\usepackage{a4,amssymb,amsfonts,theorem,stmaryrd}
\usepackage[all]{xy}

\allowdisplaybreaks

\theorembodyfont{\slshape}
\newtheorem{thm}{Theorem}[section]
\newtheorem{prop}[thm]{Proposition}

\newtheorem{lemma}[thm]{Lemma}
\newtheorem{cor}[thm]{Corollary}
\theorembodyfont{\rmfamily} 
\newtheorem{Def}[thm]{Definition}
\newtheorem{rem}[thm]{Remark}

\newenvironment{pf}{\textit{Proof.}}{\hfill q.e.d.}

\newcommand{\R}{\mathbb{R}}
\newcommand{\N}{\mathbb{N}}
\newcommand{\Z}{\mathbb{Z}}
\newcommand{\C}{\mathbb{C}}
\renewcommand{\H}{\mathbb{H}}

\newcommand{\cp}{cp}
\newcommand{\Ccinf}{C_{\cp}^{\infty}}
\newcommand{\Cinf}{C^{\infty}}

\newcommand{\E}{\mathcal{E}}
\newcommand{\F}{\mathcal{F}}

\newcommand{\cone}{\operatorname{cone}}
\newcommand{\sn}{\operatorname{sn}_{\kappa}}
\newcommand{\cs}{\operatorname{cs}_{\kappa}}

\newcommand{\ct}{\operatorname{ct}_{\kappa}}
\newcommand{\M}{\mathbf{M}^3_{\kappa}}
\newcommand{\Hyp}{\mathbf{H}}
\newcommand{\Sph}{\mathbf{S}}
\newcommand{\Euc}{\mathbf{E}}

\newcommand{\Ric}{\operatorname{Ric}}

\newcommand{\dom}{\operatorname{dom}}
\newcommand{\spec}{\operatorname{spec}}

\newcommand{\im}{\operatorname{im}}

\renewcommand{\min}{min}
\renewcommand{\max}{max}
\newcommand{\Fr}{F}

\newcommand{\Isom}{\operatorname{Isom^+}}
\newcommand{\isom}{\mathfrak{isom^+}}
\newcommand{\g}{\mathfrak{g}}
\renewcommand{\k}{\mathfrak{k}}
\newcommand{\p}{\mathfrak{p}}
\renewcommand{\sl}{\mathfrak{sl}}
\newcommand{\su}{\mathfrak{su}}
\newcommand{\SL}{\operatorname{SL}}
\newcommand{\GL}{\operatorname{GL}}
\newcommand{\SU}{\operatorname{SU}}
\newcommand{\PSL}{\operatorname{PSL}}
\newcommand{\PSU}{\operatorname{PSU}}
\newcommand{\SO}{\operatorname{SO}}

\newcommand{\Stab}{\operatorname{Stab}}
\newcommand{\End}{\operatorname{End}}
\newcommand{\dev}{\operatorname{dev}}
\newcommand{\hol}{\operatorname{hol}}
\newcommand{\id}{\operatorname{id}}
\newcommand{\rot}{\operatorname{rot}}

\newcommand{\tr}{\operatorname{tr}}
\newcommand{\Tr}{\operatorname{Tr}}
\newcommand{\interior}{\operatorname{int}}
\newcommand{\diag}{\operatorname{diag}}
\newcommand{\rank}{\operatorname{rank}}
\renewcommand{\Im}{\operatorname{Im}}
\newcommand{\res}{\operatorname{res}}
\newcommand{\cen}{\operatorname{center}}
\newcommand{\supp}{\operatorname{supp}}

\newfont{\hugemath}{cmsy10 scaled 3000}

\begin{document}

\title{Local rigidity of 3-dimensional cone-manifolds}
\author{Hartmut Weiss}
\maketitle

\begin{abstract}
\noindent We investigate the local deformation space of 3-dimensional cone-manifold structures of constant curvature $\kappa \in \{-1,0,1\}$ and cone-angles $\leq \pi$. Under this assumption on the cone-angles the singular locus will be a trivalent graph. In the hyperbolic and the spherical case our main result is a vanishing theorem for the first $L^2$-cohomology group of the smooth part of the cone-manifold with coefficients in the flat bundle of infinitesimal isometries. We conclude local rigidity from this. In the Euclidean case we prove that the first $L^2$-cohomology group of the smooth part with coefficients in the flat tangent bundle is represented by parallel forms. 
\end{abstract}


\section{Introduction}
A $3$-dimensional cone-manifold is a $3$-manifold $C$ equipped with a singular
geometric structure. More precisely, $C$ carries a length metric, which is in
the complement of an embedded geodesic graph $\Sigma$ induced by a smooth
Riemannian metric of constant sectional curvature $\kappa\in\R$. $\Sigma$ is called the
singular locus and $M=C \setminus
\Sigma$ the smooth part of $C$. Neighbourhoods of singular points are modelled on
cones of curvature $\kappa$ over $2$-dimensional cone-manifolds diffeomorphic
to $S^2$. One associates with each edge
contained in $\Sigma$ the so-called cone-angle, which is a positive real
number. If all cone-angles are $\leq\pi$, then a connected component of
$\Sigma$ is either a (connected) trivalent graph or a circle.

$3$-dimensional cone-manifolds arise naturally in the geometrization of
$3$-dimensional orbifolds, cf.~\cite{Thu}. The concept of cone-manifold can be viewed as a
generalization of the concept of geometric orbifold, where the cone-angles are
no longer restricted to the set of orbifold-angles, which are rational
multiples of $\pi$.

The deformation space of cone-manifold structures on a given cone-3-manifold $C$
with fixed topological type $(C,\Sigma)$ plays a significant role in the proof
of the Orbifold Theorem, which has recently been completed by M.\ Boileau,
B.\ Leeb and J.\ Porti, cf.~\cite{BLP1} and \cite{BLP2}. The proof of the
Orbifold Theorem in the general case requires the analysis of cone-manifold
structures with cone-angles $\leq\pi$, where the singular locus is allowed to have trivalent vertices. The case, where the singular locus is a union of circle components, i.e.~a link in $C$, has earlier been settled by M.\ Boileau and J.\ Porti, cf.~\cite{BP}.

In this article we investigate local properties of the deformation space of cone-manifold structures with cone-angles $\leq\pi$. We consider the general case under this cone-angle restriction, where trivalent vertices are allowed. In particular we prove local rigidity in the spherical and in the hyperbolic case.  

In the hyperbolic case there are some important results known. There is on the
one hand Garland-Weil local rigidity (cf.~\cite{Gar}), which applies in any
dimension $\geq 3$ to the space of complete, finite-volume hyperbolic
structures on a given hyperbolic manifold. On the other hand, C.\ Hodgson and
S.\ Kerckhoff proved a local rigidity result for $3$-dimensional hyperbolic
cone-manifolds, cf.~\cite{HK}. Their proof applies to the case, where the
singular locus $\Sigma$ is a link in $C$, but where the cone-angles are allowed to be $\leq 2\pi$.

Our main technical result is a vanishing theorem for $L^2$-cohomology on the smooth part $M$ of the cone-manifold $C$ with coefficients in the flat vector-bundle of infinitesimal isometries. $L^2$-cohomology is by definition the cohomology of the subcomplex of the de-Rham complex, which consists of those forms $\omega$ such that $\omega$ and $d\omega$ are $L^2$-bounded. 
\begin{thm}\label{L2vanishing}
Let $C$ be a $3$-dimensional cone-manifold of curvature $\kappa \in
\{-1,0,1\}$ with cone-angles $\leq \pi$. Let $(\E,\nabla^{\E})$ be the vector-bundle of infinitesimal isometries of
$M=C\setminus\Sigma$ with its natural flat connection. In the Euclidean case let
$\E_{trans}\subset\E$ be the parallel subbundle of infinitesimal translations.
Then in the hyperbolic and the spherical case
$$
H^1_{L^2}(M,\E) = 0\,,
$$
while in the Euclidean case 
$$
H^1_{L^2}(M,\E_{trans})\cong\{\omega\in \Omega^1(M,\E_{trans})\,\vert\, \nabla\omega=0\}\,.
$$
\end{thm}
The proof of this theorem is analytic in nature. The main difficulty is caused
by the non-completeness of the metric on $M$. On a complete Riemannian manifold the Hodge-Laplace
operator on differential forms is known to be essentially selfadjoint,
cf.~\cite{BL1} and the references therein. This is something we cannot expect to
hold here.
 
On the other hand, the fact that the
singularities of the metric are of iterated cone type allows us to apply separation of variables techniques. This has already been explored by
J.\ Cheeger, cf.~\cite{Ch1}. 

One main ingredient is a Hodge-theorem for cone-manifolds, which allows us to identify $L^2$-cohomology spaces with the
kernel of a certain selfadjoint extension of the Laplacian on forms.
The second one is a Bochner-Weitzenb\"ock formula for the Laplacian on $1$-forms with
values in the flat vector-bundle $\E$, resp.\ the parallel
subbundle $\E_{trans}\subset \E$ in the Euclidean case.

The essence of the Bochner technique is that the Weitzenb\"ock formula may be
used to bound the Laplacian on compactly supported $1$-forms from below: 
$
\left\langle \Delta\omega,\omega \right\rangle_{L^2} \geq C \left\langle \omega,\omega\right\rangle_{L^2}
$
for all $\omega \in \Omega^1_{cp}(M,\E)$ and some $C>0$. If we can show that this lower bound
extends to hold for the selfadjoint extension given to us by the
Hodge-theorem, we can conclude $H^1(M,\E)=0$.
In the Euclidean case, where one does not get a positive lower bound, one has to
vary this argument a little.

In the complete, finite-volume case this settles everything in view of the essential
selfadjointness of the Hodge-Laplacian (cf.~\cite{Gar}). In our case it requires a more detailed study of
the selfadjoint extensions of the Hodge-Laplacian. 
Here we use techniques introduced by
J.\ Br\"uning and R.\ Seeley, cf.~\cite {BS3}, along with some basic functional analytic
properties of the de-Rham complex presented in a very convenient form in $\cite{BL1}$.

In the hyperbolic and in the spherical case we may conclude local rigidity
from this; let us now briefly discuss the results:

If $\Sigma\subset C$ is the singular locus, for $\varepsilon > 0$ let
$U_{\varepsilon}(\Sigma)$ be the smooth part of the $\varepsilon$-tube of $\Sigma$ in
$C$, i.e.~$U_{\varepsilon}(\Sigma)=B_{\varepsilon}(\Sigma)\cap M$. Let $M_{\varepsilon}=M\setminus U_{\varepsilon}(\Sigma)$,
which is topologically a manifold with boundary. Let $\mu_i$ be the meridian
curve around the $i$-th edge of $\Sigma$.
 
In the hyperbolic case, the holonomy representation of the smooth, but incomplete hyperbolic structure
on $M$ lifts to a representation
$$
\hol : \pi_1 M \longrightarrow \widetilde{\Isom}\Hyp^3=\SL_2(\C)\,.
$$
Let $R(\pi_1M,\SL_2(\C))$ denote the set of representations of
$\pi_1M$ in $\SL_2(\C)$ equipped with the compact-open topology.
The set-theoretic quotient of $R(\pi_1M,\SL_2(\C))$ by the conjugation action of $\SL_2(\C)$
equipped with the quotient topology is denoted by
$X(\pi_1M,\SL_2(\C))$. For a representation $\rho \in R(\pi_1M,\SL_2(\C))$ let
$t_{\mu_i}(\rho)=\tr \rho(\mu_i)$. Clearly the functions $t_{\mu_i}$ are
invariant under conjugation and descend to $X(\pi_1M,\SL_2(\C))$.

The above defined spaces may be badly behaved in general, but near the holonomy
representation of a hyperbolic cone-manifold structure we can establish
smoothness and the following parametrization:
\begin{thm}\label{mainthmhyp}
Let $C$ be a hyperbolic cone-3-manifold with cone-angles $\leq\pi$. Let
$\{\mu_1, \ldots, \mu_N\}$ be the family of meridians, where
$N$ is the number of edges contained in $\Sigma$. Then the map $$X(\pi_1M,\SL_2(\C))\rightarrow \C^N, \chi \mapsto
(t_{\mu_1}(\chi),\ldots,t_{\mu_N}(\chi))$$ is locally biholomorphic near $\chi=[\hol]$.
\end{thm}
The quotient space
$X(\pi_1M,\SL_2(\C))$ may be considered, at least locally, as the deformation
space of hyperbolic structures on $M$. Hyperbolic cone-manifold structures correspond 
to representations, where the meridians $\mu_i$ map to elliptic elements in $\SL_2(\C)$. 
Therefore the previous theorem implies
local rigidity in the following strong sense:
\begin{cor}[local rigidity]\label{localrigidityhyp} 
Let $C$ be a hyperbolic cone-3-manifold with cone-angles $\leq\pi$. 
Then the set of cone-angles $\{\alpha_1,\ldots,\alpha_N\}$,  where $N$ is the number of edges contained in $\Sigma$, provides a local parametrization of the space of hyperbolic cone-manifold structures near the given structure on $M$. In particular, there are no deformations leaving the cone-angles fixed. 
\end{cor}
In the spherical case, the holonomy representation of the smooth, but incomplete
spherical structure on $M$ lifts to a product representation 
$$
\hol=(\hol_1,\hol_2):\pi_1M \longrightarrow \widetilde{\Isom}\Sph^3 = \SU(2)\times\SU(2)\,.
$$
Again for a representation $\rho \in R(\pi_1M,\SU(2))$ we set
$t_{\mu_i}(\rho)=\tr \rho(\mu_i)$. As above, the functions $t_{\mu_i}$ are
invariant under conjugation and descend to $X(\pi_1M,\SU(2))$. 

Following \cite{Por} we will say that a cone-3-manifold $C$ is {\em Seifert
  fibered} if $C$ carries a Seifert fibration such that the components of
  $\Sigma$ are leaves of the fibration. In particular $\Sigma$ is a link and $M=C\setminus \Sigma$ is a Seifert fibered 3-manifold.
In the statement of the following result we have to include the additional
hypothesis ``$C$ not Seifert fibered'' to ensure that the representations $\hol_i:\pi_1M\rightarrow
\SU(2)$ are non-abelian.
\begin{thm}\label{mainthmsph}
Let $C$ be a spherical cone-3-manifold with cone-angles $\leq\pi$, which is not Seifert fibered. Let $\{\mu_i,\ldots,\mu_N\}$ be the family of meridians, where
$N$ is the number of edges contained in $\Sigma$.
Then the map $$X(\pi_1M,\SU(2))\rightarrow \R^N, \chi_i \mapsto
(t_{\mu_1}(\chi_i),\ldots,t_{\mu_N}(\chi_i))$$ is a local diffeomorphism near
$\chi_i=[\hol_i]$ for $i \in \{1,2\}$.
\end{thm}
As in the hyperbolic case we conclude local rigidity from this:
\begin{cor}[local rigidity]\label{localrigiditysph}
Let $C$ be a spherical cone-3-manifold with cone-angles $\leq\pi$,
which is not Seifert fibered. 
Then the set of cone-angles $\{\alpha_1,\ldots,\alpha_N\}$,  where $N$ is the number of edges contained in $\Sigma$, provides a local parametrization of the space of spherical cone-manifold structures near the given structure on $M$. In particular, there are no deformations leaving the cone-angles fixed. 
\end{cor}
The geometric significance of the cohomological result in the Euclidean case
is subject to further investigation.\\
\\
The results of this article are contained in my doctoral thesis.
I would like to thank Bernhard Leeb, my thesis advisor, for his support and
encouragement. I am indebted to Joan Porti for answering many of
my questions concerning representation varieties and related
things. Furthermore, I would like to thank Daniel Grieser for explaining
various aspects of analysis on singular manifolds to me. 

\section{Cone-manifolds}\label{conemanifolds}
For $\kappa \in \R$ let $\sn$ and $\cs$ be the unique solutions of the ODE
$$
f''(r)+\kappa f(r)=0
$$
subject to the initital conditions
\begin{align*}
\sn(0)=0 &\text{\quad and\quad}\sn'(0)=1 \\
\cs(0)=1&\text{\quad and\quad}\cs'(0)=0\,.
\end{align*}
If $(N,g^N)$ is a Riemannian manifold we define for $\kappa\in\R$ and
$\varepsilon>0$ (and $\varepsilon < \pi/\sqrt{\kappa}$ if $\kappa >0$) the
$\varepsilon$-truncated $\kappa$-{\em cone} over $N$ to be the space
$$
\cone_{\kappa,(0,\varepsilon)}N = (0,\varepsilon) \times N 
$$
equipped with the Riemannian metric
$$
g=dr^2 + \sn^2(r)g^N\,.
$$

A {\em cone-surface} $S$ of curvature $\kappa\in\R$ is a compact, oriented
surface which carries a length metric with the property that there are a
finite number of points $\{x_1, \ldots ,x_k\}\subset S$ (the {\it
  cone-points}) and numbers  $\{\alpha_1, \ldots ,\alpha_k\}\subset \R^k_{+}$
(the {\it cone-angles}), such that $N=S\setminus\{x_1, \ldots ,x_k\}$ is
a smooth Riemannian manifold of curvature $\kappa$ and furthermore
the smooth part of the $\varepsilon$-ball around each cone-point
$U_{\varepsilon}(x_i)=B_{\varepsilon}(x_i)\cap N$ is isometric to the $\kappa$-cone over the circle of length $\alpha_i$.  

We will also use the notation $\interior S=S\setminus\{x_1, \ldots
,x_k\}$ for the smooth part of a cone-surface $S$. For $\kappa\in\{-1,0,1\}$ we will
call $S$ respectively hyperbolic, Euclidean or spherical. Let us call the
homeomorphism type of
$(S,\{x_1,\ldots,x_k\})$ the topological type of $S$.

Using a version of the Gauss-Bonnet theorem for cone-surfaces, it is easy to classify the spherical cone-surfaces $S$ with cone-angles $\leq\pi$.
The underlying space has to be $S^2$ and there are two types:
$$
S=\left\{ \begin{array}{c@{\quad}c}
\Sph^2(\alpha,\beta,\gamma) & \text{or}\vspace{2mm}\\
\Sph^2(\alpha,\alpha)\,. &  \end{array} \right.
$$ 
$\Sph^2(\alpha,\beta,\gamma)$ is the double of a spherical triangle with angles
$\alpha/2,\beta/2,\gamma/2$.
$\Sph^2(\alpha,\alpha)$ is the double of a spherical bigon with angles $\alpha/2, \alpha/2$.
Spherical cone-surfaces with cone-angles $\leq \pi$ are rigid, i.e.~they are determined up to isometry by the topological
type and the set of cone-angles.

A {\em cone-3-manifold} $C$ of curvature $\kappa\in\R$ is a compact, oriented
3-manifold which carries a length metric with the property that there is a
distinguished subset $\Sigma\subset C$ (the {\it singular locus}) such that
$M=C\setminus\Sigma$ is a smooth Riemannian manifold of curvature $\kappa$ and
furthermore the smooth part of the $\varepsilon$-ball around each singular
point $U_{\varepsilon}(x)=B_{\varepsilon}(x)\cap M$ is isometric to the
$\kappa$-cone over $\interior S_x$ for a
spherical cone-surface $S_x$.  
 
We will also use the notation $\interior C=C\setminus\Sigma$ for the smooth
part of a cone-3-manifold $C$. For $\kappa\in\{-1,0,1\}$ we will
call $C$ respectively hyperbolic, Euclidean or spherical.
Let us call the homeomorphism type of $(C,\Sigma)$ the topological type of $C$. 

If $x\in\Sigma$ is a
singular point
then we call $S_x$ the link of $x$ in $C$.
The hypothesis that the underlying space $C$ is a manifold implies that the
links of singular points are cone-surfaces with underlying space $S^2$. 
If the cone-angles are $\leq \pi$ we in particular obtain that links of singular
points are either $\Sph^2(\alpha,\beta,\gamma)$ or $\Sph^2(\alpha,\alpha)$. This
implies that the singular locus $\Sigma$ is a trivalent graph embedded
geodesically into
$C$.

Cone-manifolds with cone-angles $\leq 2\pi$ satisfy a lower curvature bound in
the triangle comparison sense and may be studied from a synthetic point of view.
This is pursued in \cite{BLP2}.

Foundational material on the geometry of 2- and 3-dimensional cone-manifolds
as well as an outline of the authors' approach to the Orbifold Theorem is contained in \cite{CHK}.

\section{Analysis on cone-manifolds}
By analysis on a cone-manifold $C$ we mean analysis on $M=C\setminus\Sigma$,
the smooth part of $C$. $M$ is a smooth Riemannian manifold, but incomplete if
$\Sigma$ is nonempty. This
causes the main difficulties here.

In this chapter we discuss some functional analytic properties of differential operators on
noncompact manifolds. In contrast to the compact situation one has to
distinguish more carefully between a differential operator acting on smooth,
compactly supported sections of some vector-bundle and its closed realizations
as an unbounded operator on the Hilbert space of $L^2$-sections.

\subsection{Differential operators on noncompact manifolds}
Let $(M,g)$ be a Riemannian manifold (possibly noncompact, possibly incomplete)
and let
$(\E,h^{\E}),(\F,h^{\F})$ be hermitian vector-bundles over $M$. The naturally associated
$L^2$-spaces $L^2(\E)$, resp. $L^2(\F$) depend on the quasi-isometry classes of the metrics $g$ and $h^{\E}$, resp. $h^{\F}$.

We consider a
differential operator $P$ acting on sections of $\E$ as an unbounded,
densely defined operator with domain the compactly supported sections:
$$
P:L^2(\E)\supset\dom P=C_{\cp}^{\infty}(\E)\longrightarrow L^2(\F).
$$
The {\em formal adjoint} of a differential operator $P$
$$
P^t:L^2(\F)\supset\dom P^t=C_{\cp}^{\infty}(\F)\longrightarrow L^2(\E)
$$
is uniquely defined by the relation $\langle Ps,t\rangle=\langle
s,P^tt\rangle$ to hold for all $s\in \Ccinf(\E)$ and $t\in \Ccinf(\F)$
. $P^t$ is again a differential
operator, hence densely defined.
$P$ is said to be {\em symmetric} (or {\em formally selfadjoint}) if $\E=\F$ and $\langle Ps,t\rangle=\langle
s,Pt\rangle$ for all $s,t \in \Ccinf(\E)$.

The formal adjoint is not to be confused with the {\em adjoint} $P^*$ in the sense of unbounded operator theory. The domain of $P^*$ is given as follows:
$$
\dom P^*=\{s\in L^2(\F)\vert u\mapsto \langle Pu,s \rangle \text{ bounded for }
u\in \dom P \}.
$$
Since $P$ is densely defined there is a unique $t\in L^2(\E)$ such that $\langle
Pu,s \rangle=\langle u,t \rangle$ holds for all $u\in \dom P$. Then let $P^*s=t$ by
definition. $P^*$ is a closed operator. Recall that a linear operator $A$ is called (graph-){\em closed} if
$\dom A$ equipped with the graph norm $\| x\|_A=(\| x\|^2+\|
Ax\|^2)^{\frac{1}{2}}$ is complete.

 $P^*$ obviously extends $P^t$
(which we as usual denote by $P^t\subset P^*$), in particular $P^*$ is densely defined. Note that $P$ is symmetric if and only if $P\subset P^*$.
A natural question to ask is if $P$ admits closed extensions, and this is in
fact always the case. Define
$$
P_{\max}=(P^t)^*
$$
and
$$
P_{\min}=P^{**}.
$$
$P^{**}$ is well-defined since $P^*$ is densely defined. $P^{**}$ then equals
$\overline{P}$, the (graph-)closure of $P$, i.e.~the domain of $P_{\min}$ can be
characterized as follows:
\begin{align*}
\dom P_{\min}=& \{s \in L^2(\E) \vert \exists (s_n)_{n\in\N}\subset \dom P \text{ such that }
s_n\rightarrow s \text{ in }L^2(\E)\\
& \text { and }(Ps_n)_{n\in\N}\text{ is a Cauchy sequence
in }L^2(\F)\}\,,
\end{align*}
and $P_{\min}(s)=\lim_{n\rightarrow\infty}Ps_n$.

We say that $Ps=t$ in the distributional sense if $\langle s,P^t u \rangle
=\langle t,u \rangle$ holds for all $u \in \Ccinf(\F)$. The domain of $P_{\max}$ may
then be written as:
$$
\dom P_{\max}=\{s\in L^2(\E)\vert Ps \in L^2(\F) \}\,,
$$
and $P_{\max}(s)=Ps$ in the distributional sense. Clearly $P_{\min}\subset P_{\max}$
and both are closed extensions of $P$. $P_{\max}$ is maximal with respect to
having $\Ccinf(\F)$ in the domain of its adjoint, i.e.~$P_{\max}^*$ still extends $P^t$.

If $P$ is symmetric we ask for selfadjoint extensions. Recall that a closed
symmetric operator $A$ is called {\em selfadjoint} if $A=A^*$. $P$ is called
{\em essentially selfadjoint} if $P_{\min}$ is selfadjoint. Since for a
symmetric operator one has $P_{\max}=P^*$, this is the case if and only if
$P_{\min}=P_{\max}$. Selfadjoint extensions need not exist in general.

On the other hand, if we assume that our operator $P$ is semibounded, there is
alway a distinguished selfadjoint extension which preserves the lower
bound. This feature will turn out to be particularly useful.

$P$ {\em semibounded} means by definition that there exists $c \in \R$ such that
$\langle s,Ps\rangle\geq c\langle s,s \rangle$ for all $s\in \dom P$.
Recall that a semibounded quadratic form $q:\dom q \times \dom q \rightarrow
L^2$ with lower bound $c$ is {\em closed} if and only if $\dom q$ equipped with the norm $\|
x\|_{q}=(q(x)+(1-c)\| x\|^2)^{1/2}$ is complete.
\begin{thm}[the Friedrichs extension]{\rm \cite[Thm.~X.23]{RS}}\label{Friedrichs}
Let $P$ be a semibounded symmetric operator and let $q(s,t)=\langle s,Pt\rangle$
for $s,t \in \dom P$. Then $q$ is a closable quadratic form and the closure $\overline{q}$ is
the quadratic form of a unique selfadjoint operator $P_{\Fr}$, the so-called
Friedrichs extension of $P$. $\dom P_{\Fr}$ is contained in $\dom\,\overline{q}$ and $P_{\Fr}$ is the
only selfadjoint extension of $P$ with this property. Furthermore, $P_{\Fr}$ satisfies the same lower bound as $P$.
\end{thm}
In the formulation of the following theorem as for the rest of the article we
adopt the usual convention $\dom AB = \{x\in\dom B\vert Bx \in \dom A\}$.
\begin{thm}[von Neumann]{\rm \cite[Thm.~X.25]{RS}}\label{Neumann}
Let $A$ be a closed densely
defined operator. Then $A^*\!A$ 
is selfadjoint.
\end{thm}  
For a differential operator of the form $P=D^t\!D$ we obtain for its quadratic form $q(s)=\langle Ds,Ds
\rangle \geq 0$ and therefore $\dom \overline{q}=\dom D_{\min}$.
A consequence of von Neumann's theorem (Theorem \ref{Neumann}) is (with $A=D_{\min}$) that $D^t_{\max}D_{\min}$ is a selfadjoint
extensions of $P$.
On the other hand, $\dom D^t_{\max}D_{\min}$ is certainly contained in $\dom
D_{\min}=\dom\overline{q}$. Therefore we get as an important corollary:
\begin{cor}\label{CorFriedrichs}
$D^t_{\max}D_{\min}$ is the Friedrichs extension of $D^t\!D$.
\end{cor}

\subsection{The de-Rham complex}
Let $(\E,\nabla^{\E})$ be a flat vector-bundle equipped with a hermitian metric
$h^{\E}$. The metric $h^{\E}$ will not necessarily be assumed to be parallel with
respect to $\nabla^{\E}$. 
We denote the
exterior derivative coupled with the flat connection again by $d$. As an
operator
$$
d:\Omega^{\bullet}_{\cp}(M,\E)\rightarrow\Omega^{\bullet+1}_{\cp}(M,\E)\,,
$$
$d$ is uniquely determined by the relation $d(\alpha\varotimes s)=d\alpha\varotimes
s+(-1)^{\vert\alpha\vert}\alpha\varotimes\nabla s$, where $\alpha$ is an ordinary form
and $s$ a section of $\E$.

Since $d^i_{\max}(\dom d^i_{\max}) \subset \dom d^{i+1}_{\max}$ and $d^{i+1}_{\max} \circ d^i_{\max} = 0$, we can consider the $d_{\max}$-complex
$$
\ldots \longrightarrow \dom d^i_{\max} \overset{d^i_{\max}}{\longrightarrow} \dom d^{i+1}_{\max} \longrightarrow \ldots
$$
In fact, $d_{max}$ is a particular choice of {\it ideal boundary condition}, cf.~\cite{Ch1}, and the $d_{\max}$-complex is a particular instance
of a so-called {\it Hilbert complex}, see \cite{BL1} for the definition and a general discussion.

Recall that the {\em Hodge-Dirac} operator $D=d+d^t$ decomposes as a direct sum
$D=D^{ev}\varoplus D^{odd}$, where
$$
D^{ev}:\Omega_{\cp}^{ev}({M,\E})\longrightarrow\Omega_{\cp}^{odd}(M,{\E})
$$
and
$$
D^{odd}=(D^{ev})^t:\Omega_{\cp}^{odd}({M,\E})\longrightarrow\Omega_{\cp}^{ev}(M,{\E})\,.
$$
We obtain closed extensions of
$D,D^{ev}$ and $D^{odd}$ by setting
$$
D(d_{\max})=d_{\max} + d^t_{\min}
$$
and
$$
D(d_{\max})^{ev/odd}=(d_{\max}+d^t_{\min})^{ev/odd}\,.
$$
Here we adopt the usual convention $\dom A+B = \dom A \cap
\dom B$.
Note in particular that $d^t_{\min}=d_{\max}^*$.
Since $d_{\max}$ and $d_{\max}^*$ are closed operators and $(\ker d_{\max})^{\perp}$ and $(\ker d_{\max}^*)^{\perp}$ are
orthogonal, it follows that
$D(d_{\max})^{odd}=(D(d_{\max})^{ev})^*$ and in particular that $D(d_{\max})$ is a
selfadjoint extension of $D$.

Note that we do {\it not} claim that in
general the extension $D(d_{\max})$ equals the maximal extension of $D$ itself.

Recall that the {\em Hodge-Laplace} operator is the square of the Hodge-Dirac operator:
$$
\Delta=D^2=dd^t+d^td\,.
$$
Von Neumann's Theorem (Theorem \ref{Neumann}) implies that
$$
\Delta(d_{max})=D(d_{\max})^2=d_{\max}d_{\min}^t+d_{\min}^td_{\max}
$$
is a selfadjoint extension of $\Delta$.
Note again that this extension need not be equal to the maximal
extension of $\Delta$.
\begin{lemma}
$\Delta_{\Fr}=D_{\max}D_{\min}$
\end{lemma}
\begin{pf}
The assertion follows from Corollary \ref{CorFriedrichs}.
\end{pf}\\
\\
We single out the following consequence since it is the basis for our main line of
argument towards the adaptation of the classical Bochner technique in our singular context.
\begin{cor}\label{CorDeltaFriedrichs}
If $D$ is essentially selfadjoint, then $\Delta_{\Fr}=\Delta(d_{\max})$.
\end{cor}
\begin{pf}
If $D$ is essentially selfadjoint, then since $D(d_{\max})$ is a
selfadjoint extension of $D$, we obtain
$D_{\min}=D(d_{\max})=D_{\max}$. Now the assertion follows from the
previous lemma.
\end{pf}\\
\\
Once essential selfadjointness of $D$ is established, this result allows one
to extend lower bounds obtained for $\Delta$ on compactly supported forms to
$\Delta(d_{\max})$ on its respective domain. Our concern for this particular
extension will become clear from the next section.

\subsection{Hodge theory}

To define $L^2$-{\em cohomology} we consider the following subcomplex of the
de-Rham complex:
\begin{align*}
\Omega^i_{L^2}(M,\E) &= \{\omega \in \Omega^i(M,\E) | w \in L^2 \,\text{and}\, dw \in L^2\}\\
&= \dom d^i_{\max} \cap \Omega^i(M,\E)\,, 
\end{align*}
which we will refer to as the smooth $L^2$-complex.
$L^2$-cohomology is by definition the cohomology of the smooth $L^2$-complex, i.e.
$$
H^i_{L^2}(M,\E)=\ker d^i\cap\Omega^i_{L^2}(M,\E) / d^{i-1}\Omega^{i-1}_{L^2}(M,\E)\,.
$$
Let us denote the cohomology of the $d_{\max}$-complex by
$$
H_{\max}^i=\ker d_{\max}^i / \im d_{\max}^{i-1}.
$$
We define the $d_{\max}$-{\em harmonic} $i$-forms to be
$$
\mathcal{H}_{\max}^i = \ker d_{\max}^i\cap\ker(d^{i-1})^t_{\min}\,.
$$ 
The following theorem is due to Cheeger, cf.~\cite{Ch1},
the corresponding statement in a slightly more general setting may be found in \cite{BL1}.
\begin{thm}
The inclusion $\Omega_{L^2}^i(M,\E) \hookrightarrow \dom d^i_{\max}$ induces
an isomorphism on the level of cohomology: $H^i_{L^2}(M,\E) \cong  H^i_{\max}$.
\end{thm}
There is a basic Hodge
theorem for the $d_{\max}$-complex, which goes back to Kodaira, cf.~\cite{Kod}, while \cite{BL1} prove a similar statement in the context of Hilbert
complexes.
\begin{thm}[weak Hodge-decomposition]\label{weakHodge}
For each $i$
there is an orthogonal decomposition
$$
L^2(\Lambda^iT^*M\otimes\E)=\mathcal{H}_{\max}^i\oplus\overline{\im d^{i-1}_{\max}}\oplus\overline{\im (d^i)^t_{\min}}.
$$
and furthermore
$$
\mathcal{H}_{\max}^i=\ker \Delta^i(d_{\max})=\ker
D(d_{\max})\cap L^2(\Lambda^iT^*M\otimes\E).
$$
\end{thm}
We define a map
\begin{align*}
\iota: \mathcal{H}_{\max}^i  &\longrightarrow  H_{\max}^i\\
\alpha                                  &\longmapsto  \alpha + \im d_{\max}^{i-1}
\end{align*}
Injectivity of $\iota$ is equivalent to $\im d_{\max}^{i-1} \cap \ker(d^{i-1})^t_{\min}=0$, which is always the case,
since
$$\overline{\im d_{\max}^{i-1}}=(\ker(d^{i-1}_{\max})^*)^{\perp} =(\ker(d^{i-1})^t_{\min})^{\perp}\,.$$
Surjectivity of $\iota$ is equivalent to
$$ \im d_{\max}^{i-1}=\overline{\im d_{\max}^{i-1}}\,,$$ 
therefore we obtain the following enhancement of the
Hodge decomposition, which is due to Cheeger (cf.~\cite{Ch1}) in the case of the
$d_{\max}$-complex. Again a more general statement may be found in \cite{BL1}.
\begin{thm}[strong Hodge-decomposition]\label{strongHodge}
If \,$\im d_{\max}^{i-1}$ is closed for all $i$, then for each $i$
there is an orthogonal decomposition
$$
L^2(\Lambda^iT^*M\otimes\E)=\mathcal{H}_{\max}^i\oplus\im d_{\max}^{i-1}\oplus\im(d^i)^t_{\min}
$$
and furthermore $\iota: \mathcal{H}_{\max}^i  \rightarrow  H_{\max}^i$ is an isomorphism.
\end{thm}
A sufficient condition for $d_{\max}^{i-1}$ to have closed range is
finite dimensionality of $H_{\max}^i$ on the one hand, since $\ker d_{\max}^i / \im
d_{\max}^{i-1}$ finite dimensional implies that $\im d_{\max}^{i-1}$
is closed in $\ker d_{\max}^i$, hence in
$L^2(\Lambda^iT^*M\otimes\E)$. Note that by the closed-range theorem
$(d_{\max}^i)^*$ has closed range if and only if $d_{\max}^i$ has
closed range.

On the other hand, if $D(d_{\max})^{ev}$ has closed range, then
$d_{\max}^{i}$ and $(d_{\max}^{i+1})^*$ will have closed range for all $i$
even. Similarly, if
$D(d_{\max})^{odd}$ has closed range,
then $d_{\max}^{i}$ and $(d_{\max}^{i+1})^*$ will have closed range for
all $i$ odd. Since
$D(d_{\max})^{odd}=(D(d_{\max})^{ev})^*$, the closed-range
theorem implies that $D(d_{\max})^{ev}$ has closed range if and only $D(d_{\max})^{odd}$ has closed range.

It is easy to show that $D(d_{\max})^{ev}$ has closed range if $\dom
D(d_{\max})^{ev}$ equipped with the graph norm embeds into
$L^2(\Lambda^{ev}T^*M\otimes\E)$ compactly. 
This latter condition is related
to the question of discreteness of the spectra of the operators $D(d_{\max})$ and
$\Delta(d_{\max})$. Recall that an operator is said to have discrete
spectrum if its spectrum consists of a discrete set of
eigenvalues with finite multiplicities.


\section{Spectral properties of cone-manifolds}

In this chapter we apply the techniques of Br\"uning and Seeley
to analyze the closed extensions of the Hodge-Dirac operator on a $3$-dimensional
cone-manifold. The main reference for the first order case will be \cite{BS3}. 
The analysis relies heavily on the fact that the spaces we consider are {\it locally
conical}, i.e.~neighbourhoods of points are isometric to ($\kappa$-)cones over
spaces of lower dimension. 
This allows us to apply separation of variables techniques.

To keep the exposition
self-contained here, we describe these techniques in detail. 
Furthermore we
adopt a more elementary viewpoint than in \cite {BS3}, in particular we give a direct
argument for discreteness of the relevant operators.

Let us further mention that \cite{BS3} deal with isolated conical
singularities only,
i.e.~the links of singular points are compact smooth Riemannian manifolds, where in
our case we have to allow the links of singular points to be again singular,
namely the spherical cone-surfaces $\Sph^2(\alpha,\beta,\gamma)$
and $\Sph^2(\alpha,\alpha)$. This requires some extra arguments which we will
provide as we expose the theory.

There has been a lot of work on Hodge-theory and $L^2$-cohomology of Riemannian manifolds with conical singularities, besides \cite{Ch1} and \cite{Ch2} see for example \cite{BL2}.


\subsection{Separation of variables}
Let $(N,g^N)$ be a Riemannian
manifold of dimension $n$ and let us consider $U_{\varepsilon}=
\cone_{\kappa,(0,\varepsilon)}N$ with the metric
$
g=dr^2+\sn^2(r)g^N
$.
We may think of $N$ as the (smooth
part of the) link $S_x$ of a singular point $x$ in a cone-manifold, $U_{\varepsilon}$
serves as a model for the (smooth part of the) $\varepsilon$-neighbourhood $U_{\varepsilon}(x)$ of a
singular point $x$ in $M$.

Let $(\E,\nabla^{\E})$ be a flat vector-bundle over $U_{\varepsilon}$. We will
identify the fibers of $\E$ along radial geodesics via parallel
translation using $\nabla^{\E}$, in particular we may canonically identify $\E\vert_{U_{\varepsilon}}=(0,\varepsilon) \times \E\vert_N$. Let
us further assume that $\E$ is equipped with a metric $h^{\E}$, which is not necessarily parallel with respect to $\nabla^{\E}$. We will assume instead:
\begin{enumerate}
\item[{\bf A1}] The limit $h_0^{\E}:=\lim_{r \rightarrow 0}h^{\E}(r)$
exists as a smooth metric on $\E|_N$ and is parallel with respect to
$\nabla^{\E}$. (The limit is defined using the canonical identification
$\E\vert_{U_{\varepsilon}}=(0,\varepsilon) \times \E\vert_N$ as above.)  
\end{enumerate}
Now $h_0^{\E}$ extends to a parallel metric on $\E|_{U_{\varepsilon}}$, which we
continue to denote by $h_0^{\E}$. We may write
$$
h^{\E}(\sigma, \tau) = h^{\E}_0(A\sigma,\tau)
$$
for $\sigma,\tau \in \Gamma(U_{\varepsilon},\E)$, where $A \in \Gamma(U_{\varepsilon},\End \E$) is symmetric
with respect to $h_0^{\E}$. Let us continue to denote the flat connection on $\End \E$ by $\nabla^{\E}$. We will further assume:
\begin{enumerate}
\item[{\bf A2}] $A^{-1}(\nabla^{\E} A) \in
  \Omega^1(U_{\varepsilon},\End \E)$ is bounded with respect to
  $g$ and $h^{\E}$.
\end{enumerate}
\begin{rem}\label{A123}
\begin{enumerate}
\item A2 implies that $h^{\E}$ and $h_0^{\E}$ are quasi-isometric on
  $U_{\varepsilon}$, since for $\sigma \in \Gamma(U_{\varepsilon},\E)$ satisfying
$\nabla^{\E}_{\partial/\partial r} \sigma =0$ we have
$$
\Bigl\vert\, \frac{d}{dr} \log
\frac{h^{\E}(\sigma,\sigma)}{h_0^{\E}(\sigma,\sigma)}\,\Bigr\vert = \frac{\bigl\vert
  h^{\E}(A^{-1}(\nabla_{\partial/\partial r}^{\E} A)
  \sigma,\sigma)\bigr\vert}{h^{\E}(\sigma,\sigma)} \leq C
$$
on the complement of the zero-set of $\sigma$, where $C$ is the bound on
$A^{-1}(\nabla^{\E} A)$ given by A2. 
\item If the cross-section $N$ is compact, then A2 is a direct
  consequence of A1, in the general case A2 is an additional assumption.
\item If $h^{\E}$ is already parallel with respect to $\nabla^{\E}$, then
  $h^{\E}_0 = h^{\E}$ and A1 and A2 are trivially satisified.
\end{enumerate}
\end{rem}
Let $d$ denote the exterior covariant derivative
coupled with $\nabla^{\E}$ and let $d^t$ denote the formal adjoint of $d$ with
respect to $h^{\E}$. Similarly let $d^t_0$ denote the formal adjoint of $d$ with
respect to $h^{\E}_0$. If $\iota(\nabla^{\E} A)$ denotes interior
multiplication with the $\End \E$-valued 1-form $\nabla^{\E} A$, then we have:
\begin{lemma} $d^t = d^t_0 - A^{-1}\iota(\nabla^{\E} A)$.
\end{lemma}
\begin{pf}
If $L^2_0$ denotes the $L^2$-space with respect to $g$ and $h^{\E}_0$, we have
\begin{align*}
\langle A \, d^t \eta, \xi \rangle_{L^2_0} &= \langle A\eta,  d\xi
\rangle_{L^2_0}=\langle \eta, d (A \xi) - \nabla^{\E} A \wedge \xi
\rangle_{L^2_0}\\
&= \langle A (d^t_0
\eta) - \iota(\nabla^{\E}A) \eta, \xi \rangle_{L^2_0}
\end{align*}
for $\eta \in \Omega_{cp}^{p+1}(U_{\varepsilon},\E)$ and $\xi \in
\Omega_{cp}^p(U_{\varepsilon},\E)$. In the last line we have used that
$h_0^{\E}$ is parallel with respect to $\nabla^{\E}$, hence $\nabla^{\End
    \E}A$ has values in the symmetric (w.r.t.~$h^{\E}_0$) endomorphisms of $\E$.
\end{pf}\\
\\
With $D=d + d^t$ and $D_0=d + d_0^t$ we therefore have
$$
D=D_0 - A^{-1}\iota(\nabla^{\E}A) \,.
$$
Following $\cite{BS3}$, we identify $p$-forms on the model neighbourhood $U_{\varepsilon}$ with pairs of $r$-dependent forms on $N$ via
$$
(\phi,\psi) \mapsto  \sn(r)^{(p-1)-\frac{n}{2}} \phi \wedge dr + \sn(r)^{p-\frac{n}{2}} \psi\,,
$$ 
where $\phi \in \Gamma(\pi_N^* \Lambda^{p-1}T^*N\otimes\E)$ and $\psi \in
\Gamma(\pi_N^* \Lambda^{p}T^*N\otimes\E)$. 
This correspondence preserves $L^2$-norms, if we use the parallel metric
$h^{\E}_0$:
$$
\int_0^{\varepsilon}\int_0 \vert\phi\vert_0^2 \,dr dvol_N
= \int_{U_{\varepsilon}} \sn(r)^{2(p-1)-n}\vert \phi\wedge dr
\vert^2_0 \,dvol_{U_{\varepsilon}}
$$
and
$$
\int_0^{\varepsilon}\int_N 
\vert\psi\vert_0^2 \,dr dvol_N
= \int_{U_{\varepsilon}} \sn(r)^{2p-n} \vert\psi\vert^2_0 \,dvol_{U_{\varepsilon}}\,.
$$
With respect to these decompositions the exterior differential has the
following matrix form on $U_{\varepsilon}$:
$$
 d^p=\left[ \begin{array}{cc}
                               \sn(r)^{-1}d_N^{p-1} & (-1)^p \left\{\frac{\partial}{\partial r}+(p-\frac{n}{2})\ct(r)\right\} \\
                               0 & \sn(r)^{-1}d_N^p
                               \end{array} \right] \,.
$$
By passing to the formal adjoints using $h^{\E}_0$ we obtain:
$$
  (d_0^t)_p=\left[ \begin{array}{cc}
                               \sn(r)^{-1}(d^t_N)_{p-1} & 0 \\
                               (-1)^p\left\{\frac{\partial}{\partial r}+(\frac{n}{2}-p+1)\ct(r)\right\}& \sn(r)^{-1}(d^t_N)_p
                               \end{array} \right] \,.
$$
We may identify $r$-dependent forms on $N$ of arbitrary degree with either
even forms on $U_{\varepsilon}$ via
$$
(\phi^0,\ldots,\phi^n) \mapsto \sum_{i}
\sn(r)^{2i+1-\frac{n}{2}}\phi^{2i+1}\wedge dr +
\sum_{i}\sn(r)^{2i-\frac{n}{2}}\phi^{2i}\,,
$$
or odd forms on $U_{\varepsilon}$ via
$$
(\phi^0,\ldots,\phi^n) \mapsto \sum_{i}
\sn(r)^{2i-\frac{n}{2}}\phi^{2i}\wedge dr +
\sum_{i}\sn(r)^{2i+1-\frac{n}{2}}\phi^{2i+1}\,.
$$
We obtain that the even part of the Hodge-Dirac operator associated with
$h_0^{\E}$
may be written on $U_{\varepsilon}$ as
$$
D^{ev}_0 = \frac{\partial}{\partial r} + \frac{1}{\sn(r)} B_{\kappa}(r) \,,
$$
where 
$$
B_{\kappa}(r)=D_N + \left[ \begin{array}{ccc}
   \cs(r)c_0 & & \\
   & \ddots & \\
   & & \cs(r)c_n \end{array} \right]
$$
with
$$
c_p=(-1)^p(p-\textstyle\frac{n}{2}) \,.
$$
Note that $\lim_{r\rightarrow 0}B_{\kappa}(r)$ is independent of $\kappa \in \R$, more precisely we have
$$
\lim_{r\rightarrow 0}B_{\kappa}(r)=D_N + \left[ \begin{array}{ccc}
   c_0 & & \\
   & \ddots & \\
   & & c_n \end{array} \right]\,.
$$
\begin{Def}[model operator]
Let $B=\lim_{r\rightarrow 0}B_{\kappa}(r)$ and
$$
P_B^{\kappa} = \frac{\partial}{\partial r} + \frac{1}{\sn(r)} B\,.
$$
\end{Def}
If the assumptions A1 and A2 hold, the operator
$P_B^{\kappa}$ may be used as a model operator for $D^{ev}$ on $U_{\varepsilon}$, since it captures its essential analytic features. This is
made precise by the following lemma:
\begin{lemma}\label{modeloperator}
If A1 and A2 hold,
$\dom(D^{ev})_{\max/\min}=\dom (P_B^{\kappa})_{\max/\min}$ and   
the graph norms $\|\cdot\|_{D^{ev}}$ and
$\|\cdot\|_{P_B^{\kappa}}$ are equivalent.
\end{lemma}
\begin{pf}
Since
$$
\frac{B_{\kappa}(r)-B}{\sn(r)}=\frac{\cs(r)-1}{\sn(r)}\left[ \begin{array}{ccc}
   c_0 & & \\
   & \ddots & \\
   & & c_n \end{array} \right]
$$
and
$$
\lim_{r\rightarrow 0} \frac{\cs(r)-1}{\sn(r)} = 0\,,
$$
we see that $D_0^{ev}$ differs from $P_B^{\kappa}$ just by a bounded $0$-th
order term. If the assumptions A1 and A2 hold, then the
$L^2$-norms defined by using $h^{\E}$, respectively $h^{\E}_0$, are equivalent and
$D_0^{ev}$ differs from $D^{ev}$ again by a bounded $0$-th order term. This
implies the assertion. 
\end{pf}

\subsection{The radial equation}
The operator $B$ is obviously symmetric on
$\Omega^{\bullet}_{cp}(N,\E)$. Note also that $B$ does not depend on
the radial variable $r\in(0,\varepsilon)$ any more. If $B$ is essentially
selfadjoint and has discrete spectrum, we use
the spectral decomposition of $L^2(\Lambda^{\bullet}T^*N,\E)$ with respect to $B$ 
to transform the model operator $P_B^{\kappa}$ into a family of operators
$P_b^{\kappa}$ on the interval $(0,\varepsilon)$, where $b$ ranges over the
spectrum of $B$.\\
\\
For $b\in\R$ let 
$$P_b^{\kappa}=\frac{\partial}{\partial r} + \frac{b}{\sn(r)}\,.$$ 
We will consider $P^{\kappa}_b$ acting on $\Ccinf(0,1)$. Furthermore let
$P_b=P^0_b$, i.e.
$$
P_b=\frac{\partial}{\partial r} + \frac{b}{r}\,.
$$
It is enough to study the operator $P_b$ in view of the following lemma:
\begin{lemma}\label{radialequiv}
$\dom(P_b^{\kappa})_{\max/\min}=\dom (P_b)_{\max/\min}$ and
the graph norms $\|\cdot\|_{P_b^{\kappa}}$ and $\|\cdot\|_{P_b}$ are
equivalent.
\end{lemma}
\begin{pf}
Since $P_b^{\kappa} - P_b = \varphi(r)b$ with
$$
\varphi(r)=\frac{1}{\sn(r)} - \frac{1}{r}
$$
and
$$
\lim_{r\rightarrow 0} \varphi(r) = 0\,,
$$
we see that $P_b^{\kappa}$ differs from $P_b$ just by a bounded $0$-th order
term. In the same way as before this implies the assertion.
\end{pf}\\
\\
It is useful to observe that
$$\left(P_bf\right)(r)=r^{-b}\textstyle\frac{\partial}{\partial r}(r^b f)\,,$$ 
therefore
$
P_bf=0
$
if and only if
$$
f(r)=f(1)r^{-b}\,,
$$
and 
$
P_bf=g
$
if and only if 
$$
f(r)=f(1)r^{-b}+r^{-b}\int_1^r\varrho^bg(\varrho)d\varrho\,.
$$
For any subinterval $(\delta,1)\subset(0,1)$ the graph norm of
$P_b$ is equivalent to the ordinary $H^1$-norm, since $\frac{1}{r}\in
L^{\infty}(\delta,1)$. $H^1$-functions - more generally: $W^{1,1}$-functions - on
$(\delta,1)$ are absolutely continuous on $[\delta,1]$, hence differentiable almost everywhere.
For absolutely continuous functions the fundamental theorem of calculus holds, i.e.~$\varphi\in AC([\delta,1])$ if and only if
$\varphi(r)=\varphi(1)+\int_1^r \varphi'(\varrho) d\varrho$ for
$r\in[\delta,1]$. Therefore the above integral representation remains valid for
$f\in \dom (P_b)_{\max}$ (take $\varphi(r)=r^{b}f(r)$).
It follows in particular that $f\in\dom (P_b)_{\max}$ is
continuous on $(0,1)$ and has a continuous boundary value at $r=1$, i.e.~$f \in C^0((0,1])$.\\
\\
Following \cite{BS3} we define two integral operators acting on $L^2(0,1)$:
$$
(T_{b_,1}g)(r)=r^{-b}\int_1^r\varrho^bg(\varrho)d\varrho\,,
$$
where $b$ is arbitrary, and
$$
(T_{b,0}g)(r)=r^{-b}\int_0^r\varrho^bg(\varrho)d\varrho\,,
$$
for $b>-\frac{1}{2}$. Note that $b>-\frac{1}{2}$ implies that $r^b\in L^2(0,1)$
and therefore with the Cauchy-Schwarz inequality $\int_0^r\varrho^bg(\varrho)d\varrho<\infty$. 

We start from the following estimates in \cite{BS3}, which easily follow from the Cauchy-Schwarz inequality:
\begin{lemma}{\rm \cite[Lm. 2.1]{BS3}} For $g \in L^2(0,1)$ and $r\in(0,1)$ we have
$$
\vert(T_{b,0}g)(r)\vert\leq r^{\frac{1}{2}}(2b+1)^{-\frac{1}{2}}\left(\int_0^r
\vert g(\varrho)\vert^2 d\varrho\right)^{\frac{1}{2}}
$$
for $b>-\frac{1}{2}$, and
$$
\vert(T_{b,1}g)(r)\vert \leq 
 \left\{ \begin{array}{c@{\quad,\quad}c}
r^{\frac{1}{2}}\vert 2b+1 \vert^{-\frac{1}{2}}\|g\|_{L^2(0,1)} & b
 <-\frac{1}{2}\vspace{.2cm}\\
r^{\frac{1}{2}}\vert \log r \vert^{\frac{1}{2}}\|g\|_{L^2(0,1)} & b
 =-\frac{1}{2}\vspace{.2cm}\\
r^{-b}(2b+1)^{-\frac{1}{2}}\|g\|_{L^2(0,1)} & b >-\frac{1}{2} \end{array} \right.,
$$
in particular $T_{b,1}g \in L^2(0,1)$ if $b<\frac{1}{2}$.
\end{lemma}
From this we may derive decay estimates for $f\in\dom (P_b)_{\max}$:
\begin{lemma}[decay estimates]
Let $f \in \dom (P_b)_{\max}$. Then for $r\in(0,1)$ and with $g=P_bf$ we have
$$
\vert f(r)\vert \leq 
 \left\{ \begin{array}{c@{\quad,\quad}c}
 r^{\frac{1}{2}}(2b+1)^{-\frac{1}{2}}\left(\int_0^r|g(\varrho)|^2\right)^{\frac{1}{2}} & b \geq \frac{1}{2}\vspace{.2cm}\\
 r^{-b}|f(1)|+r^{-b}(2b+1)^{-\frac{1}{2}}\|g\|_{L^2(0,1)} & b \in(-\frac{1}{2},\frac{1}{2})\vspace{.2cm}\\
 r^{\frac{1}{2}}|f(1)|+r^{\frac{1}{2}}|\log r|^{\frac{1}{2}}\|g\|_{L^2(0,1)} & b
 =-\frac{1}{2}\vspace{.2cm}\\
 r^{-b}|f(1)|+r^{\frac{1}{2}}|2b+1|^{-\frac{1}{2}}\|g\|_{L^2(0,1)}
 & b <-\frac{1}{2} \end{array}\right._.
$$
\end{lemma}
\begin{pf}
The estimates for $b<\frac{1}{2}$ follow directly from the integral
representation
$$
f(r)=r^{-b}f(1)+\left(T_{b,1}g\right)(r)
$$
and the corresponding estimates for $T_{b,1}g$ from the previous lemma. For the
case $b\geq\frac{1}{2}$ we observe that for $b\geq\frac{1}{2}$ (in fact already
for $b>-\frac{1}{2}$) $r^b\in L^2(0,1)$, hence $r^bg \in L^1(0,1)$ by the
Cauchy-Schwarz inequality. This implies that $r^bf$ has its distributional derivative in
$L^1(0,1)$ and is therefore absolutely continuous on $[0,1]$.
We obtain
$$
f(r)=r^{-b}C+\left(T_{b,0}g\right)(r)
$$
with $C=\lim_{r\rightarrow 0} r^bf(r)$. Now $r^{-b}\not\in L^2(0,1)$ for $b\geq\frac{1}{2}$, therefore $C=0$, so the
estimate for $T_{b,0}g$ gives the result. 
\end{pf}
\begin{cor}\label{CorDecay}
Let $f\in\dom(P_b)_{\max}$ and $r\in(0,1)$. If $b\not\in(-\frac{1}{2},\frac{1}{2})$, then
$$
|f(r)| \leq C(b) r^{\frac{1}{2}}(1+|\log r|^{\frac{1}{2}})\|f\|_{P_b}\,,
$$
in particular $f\in
C^0([0,1])$ with $f(0)=0$,
while if $b \in(-\frac{1}{2},\frac{1}{2})$, then 
$$
|f(r)| \leq C(b) r^{-b} \|f\|_{P_b}\,.
$$
\end{cor}
\begin{pf}
The case $b\geq\frac{1}{2}$ follows directly from the above estimates. For the
other cases we again refer to the integral representation
$$
f(r)=r^{-b}f(1)+\left(T_{b,1}g\right)(r)
$$
and observe that $r^{-b}f(1)\in L^2(0,1)$ for $b<\frac{1}{2}$. Therefore
the bound on $T_{b,1}g$ translates into a bound on $|f(1)|$ in terms of $\|f\|_{L^2(0,1)}$
and $\|g\|_{L^2(0,1)}$. This plugged into the decay estimates gives the
result, which clearly implies that $f(r)=o(1)$ as
$r\rightarrow 0$ in the first case.
\end{pf}\\
\\
The following statement is implicitly contained in Br\"uning and Seeley's parametrix
construction, cf.~\cite{BS3}:
\begin{prop}[integration by parts]
Let $\varphi\in\Cinf(0,1)$ be a
cut-off function with $\varphi\equiv 1$ near $0$ and $\varphi\equiv 0$ near
$1$. For $u\in\dom(P_b)_{\max}$ let $f=\varphi u\in\dom(P_b)_{\max}$, and let $g\in \dom (P_b^t)_{\max}$. Then for
$b\not\in(-\frac{1}{2},\frac{1}{2})$ the following holds:
$$
\left\langle (P_b)_{\max}f,g\right\rangle_{L^2(0,1)}=\left\langle f,(P_b^t)_{\max}g\right\rangle_{L^2(0,1)}
$$
\end{prop}
\begin{pf}
With $(P_b)^t=-P_{-b}$ we calculate
\begin{align*}
\left\langle (P_b)_{\max}f,g\right\rangle_{L^2(0,1)} 
&= \int_0^1 \left(\frac{\partial f}{\partial r}+\frac{rf}{b}\right)g\\
&=\lim_{\delta\rightarrow 0}\left\{\int_{\delta}^1\left(\frac{\partial f}{\partial
r}\right)g+\int_{\delta}^1 \left(\frac{rf}{b}\right)g \right\}\\
&=\lim_{\delta\rightarrow
0}\left\{\left[fg\right]_{\delta}^1-\int_{\delta}^1f\left(\frac{\partial g}{\partial
r}\right)+\int_{\delta}^1 f\left(\frac{rg}{b}\right)\right\}\\
&=\lim_{\delta\rightarrow
0}\left\{f(1)g(1)-f(\delta)g(\delta)\right\}+\left\langle f,(P_b^t)_{\max}g\right\rangle_{L^2(0,1)}.
\end{align*}
Now $f(1)=0$ and $\lim_{\delta\rightarrow 0} f(\delta)g(\delta)=0$ according to
the decay estimates. Therefore
$$
\lim_{\delta\rightarrow
0}\left\{f(1)g(1)-f(\delta)g(\delta)\right\}=0
$$
and we obtain the result.
\end{pf}\\
\\
This statement becomes wrong, if we allow $b\in(-\frac{1}{2},\frac{1}{2})$. To
see this, let $f(r)=\varphi(r)r^{-b}$ with $\varphi$ as above and
$g(r)=r^b$. Note that $P_b(r\mapsto r^{-b})=P^t_b(r\mapsto r^b)=0$, so clearly $f\in\dom(P_b)_{\max}$
and $g\in\dom(P_b^t)_{\max}$. But on the other hand
$$
\lim_{\delta\rightarrow
0}\left\{f(1)g(1)-f(\delta)g(\delta)\right\}=0-\lim_{\delta\rightarrow 0}f(\delta)g(\delta)=-1\,,
$$
so we have a boundary contribution.\\
\\ 
The preceding result allows us to conclude that we do not have to impose boundary conditions for
$P_b$ at $0$, if (and only if) $b\not\in(-\frac{1}{2},\frac{1}{2})$.  
\begin{cor}\label{radial_max/min}
Let $\varphi\in\Cinf(0,1)$ be a
cut-off function with $\varphi\equiv 1$ near $0$ and $\varphi\equiv 0$ near
$1$. For $u\in\dom(P_b)_{\max}$ let $f=\varphi u\in\dom(P_b)_{\max}$. Then
$f\in\dom(P_b)_{\min}$ for $b\not\in(-\frac{1}{2},\frac{1}{2})$. 
\end{cor}
\begin{pf}
For all $g\in\dom(P_b^t)_{\max}$ we have
$$
\left\langle (P_b)_{\max}f,g\right\rangle_{L^2(0,1)}=\left\langle
f,(P_b^t)_{\max}g\right\rangle_{L^2(0,1)}
$$
This means that $f\in\dom(P_b^t)_{\max}^*=\dom(P_b)_{\min}$.
\end{pf}\\
\\
Let $P^{\kappa}_B=\frac{\partial}{\partial r}+\sn(r)^{-1}B$ acting on
$\Ccinf\left((0,1)\times N\right)$. We will assume that $B$ is
essentially selfadjoint on $\Ccinf(N)$, i.e.~in equivalent terms
$B_{\max}=B_{\min}$, since $B$ is symmetric. We will furthermore assume that $B$ has discrete
spectrum. Let $\{\Psi_b\}_{b\in\spec B}$ be an orthonormal basis of $L^2(N)$
consisting of eigensections of $B$, where as usual each eigenvalue is repeated according to its
multiplicity.
By interior elliptic regularity, the
$\Psi_b$ are smooth. There are orthogonal
decompositions
$$L^2(N)=\overline{\bigoplus_{b\in\spec B}\R\otimes\langle\Psi_b\rangle}$$
and 
$$L^2\left(
(0,1)\times N\right)=\overline{\bigoplus_{b\in\spec B}L^2(0,1)\otimes\langle\Psi_b\rangle}\,,$$
where the closure is taken with respect to the corresponding $L^2$-norm. For $f\in L^2\left(
(0,1)\times N\right)$ we have an $L^2$-convergent expansion
$$f=\sum_{b\in\spec B}f_b\varotimes \Psi_b\,,$$ 
where 
$$f_b(r)=\int_N\left(f(r,x),\Psi_b(x)\right)dx\,.$$
Obviously we have 
$$\|f\|_{L^2\left((0,1)\times N\right)}^2=\sum_{b\in\spec B}\|f_b\|^2_{L^2(0,1)}\,.
$$
\begin{lemma}
Let $f,g \in L^2\left((0,1)\times N\right)$. Then $P^{\kappa}_Bf=g$ 
if and only if $P^{\kappa}_bf_b=g_b$ for all $b\in\spec B$.
In particular $f \in \dom (P^{\kappa}_B)_{\max}$ if and only if $f_b\in\dom(P^{\kappa}_b)_{\max}$ for all $b \in \spec B$.
\end{lemma}
\begin{pf}
Let us assume first that $P^{\kappa}_Bf=g$ with
$f,g\in L^2\left((0,1)\times N\right)$, i.e.~by definition
$\langle f,P^{\kappa,t}_B\phi\rangle_{L^2}=\langle g,\phi\rangle_{L^2}$ for all $\phi \in
\Ccinf\left((0,1)\times N\right)$. 

If 
$\varphi\in\Ccinf(0,1)$ is an arbitrary cut-off function, we claim that this
relation remains valid for $\phi=\varphi \Psi_b$ and $b\in\spec B$. 
Since by assumption $B_{\max}=B_{\min}$ we may choose sequences $\Psi_{b,n}
\in \Ccinf(N)$, which approximate $\Psi_b$ with respect to $\|\cdot\|_B$. 
Then it follows immediately that
$\varphi\Psi_{b,n}$ approximate $\varphi\Psi_b$ with respect to $\|\cdot\|_{P^{\kappa,t}_B}$.
Since
$\varphi\Psi_{b,n}\in\Ccinf\left((0,1)\times N\right)$ we have 
$$
\langle f,P^{\kappa,t}_B(\varphi\Psi_{b,n})\rangle_{L^2} = \langle g, \varphi\Psi_{b,n}\rangle_{L^2}
$$
for all $n$. By continuity we obtain
$$
\langle f,P^{\kappa,t}_B(\varphi\Psi_{b})\rangle_{L^2} = \langle g, \varphi\Psi_{b}\rangle_{L^2}\,,
$$
which proves the subclaim.
Now the left-hand side of this equation equals
\begin{align*}
\int_0^{1}\int_N \left( f, P^{\kappa,t}_B(\varphi\Psi_b)\right)
& =\int_0^{1}P_b^{\kappa,t}\varphi \int_N(f,\Psi_b)
=\int_0^{1}f_b P_b^{\kappa,t}\varphi\,,
\end{align*} 
whereas the right-hand side is given by
$$
\int_0^{1}\int_N \left (g,\varphi\Psi_b\right) =
\int_0^{1}\varphi\int_N \left(g,\Psi_b\right) = \int_0^1 g_b\varphi\,.
$$
Since $\varphi$ was arbitrary, this means that $P^{\kappa}_bf_b=g_b$ for all
$b\in\spec B$.

Conversely, if $P^{\kappa}_bf_b=g_b$ holds for all $b\in\spec B$, we have to show that
$$
\langle f,P_B^{\kappa,t}\phi\rangle_{L^2} = \langle g, \phi \rangle_{L^2}
$$ 
is true for all $\phi \in \Ccinf\left((0,1)\times N\right)$. Now
$$
\langle f,P_B^{\kappa,t}\phi\rangle_{L^2} = \sum_{b\in\spec B}\langle f_b, (P_B^{\kappa,t}\phi)_b\rangle_{L^2(0,1)}
$$
and 
$$
\langle g,\phi\rangle_{L^2} =  \sum_{b\in\spec B}\langle g_b,
\phi_b \rangle_{L^2(0,1)}\,,
$$
so we obtain the result, since $(P_B^{\kappa,t}\phi)_b = P_b^{\kappa,t}
\phi_b$.
\end{pf}
\begin{lemma}
Let $f\in\dom(P^{\kappa}_B)_{\max}$. Then 
$f \in \dom (P^{\kappa}_B)_{\min}$ if and only if $f_b \in \dom (P^{\kappa}_b)_{\min}$ for all $b
\in\spec B$.
\end{lemma}
\begin{pf}
The proof essentially uses the observation that $f\in\dom(P^{\kappa}_B)_{\min}$
if and only if 
 $\langle P^{\kappa}_Bf,g\rangle_{L^2} = \langle f,
P^{\kappa,t}_Bg\rangle_{L^2}$ for all $g\in\dom(P_B^{\kappa,t})_{\max}$.
Now the left-hand side of the
equation in question equals
$$
 \sum_{b\in\spec B} \langle (P_B^{\kappa}f)_b,g_b \rangle_{L^2(0,1)} = \sum_{b\in\spec B} \langle P^{\kappa}_bf_b,g_b \rangle_{L^2(0,1)}\,, 
$$
since $f_b\in\dom(P^{\kappa}_b)_{\max}$ and $g_b\in\dom(P_b^{\kappa,t})_{\max}$,
while the right-hand side is given by
$$
 \sum_{b\in\spec B} \langle f_b,(P_B^{\kappa,t}g)_b \rangle_{L^2(0,1)} =  \sum_{b\in\spec B} \langle f_b,P_b^{\kappa,t}g_b \rangle_{L^2(0,1)}\,.
$$
We obtain that $f\in\dom(P^{\kappa}_B)_{\min}$ if and only if $\langle P^{\kappa}_bf_b,g_g\rangle_{L^2(0,1)} = \langle f_b,
P^{\kappa,t}_Bg_b\rangle_{L^2(0,1)}$ for all $g_b\in\dom(P_b^{\kappa,t})_{\max}$, i.e.~that $f_b \in \dom (P^{\kappa}_b)_{\min}$ for all $b \in \spec B$. 
\end{pf}\\
\\
The following lemma will turn out to be decisive in the question of essential selfadjointness of
$D$ on cone-manifolds.
\begin{lemma}\label{max/min}
Let $\varphi\in\Cinf(0,1)$ be a cut-off function with $\varphi\equiv 1$ near $0$ and $\varphi\equiv 0$ near
$1$. For $u\in\dom(P^{\kappa}_B)_{\max}$ let $f=\varphi u\in\dom(P^{\kappa}_B)_{\max}$. Then
$f\in\dom(P^{\kappa}_B)_{\min}$ if $\spec B\cap(-\frac{1}{2},\frac{1}{2})=\varnothing$. 
\end{lemma}
\begin{pf}
This follows from the above discussion together with Corollary
\ref{radial_max/min} and Lemma \ref{radialequiv}.
\end{pf}\\
\\
In the following we derive certain compactness properties which will be relevant
for the question of discreteness of $D(d_{\max})$ and $\Delta(d_{\max})$
on cone-manifolds.
\begin{lemma}\label{radial_discr}
The embedding $\dom (P_b)_{\max}\hookrightarrow L^2(0,1)$ is compact for all $b\in\R$.
\end{lemma}
\begin{pf}
Given a sequence $f_n \in \dom (P_b)_{\max}$ with bound $\|f_n\|_{P_b}\leq
C$ independent of $n$, we have to extract a subsequence convergent in $L^2(0,1)$. On any 
subinterval $(\delta,1)\subset(0,1)$ the graph norm of $P_b$ is equivalent to the
ordinary $H^1$-norm, since $\frac{1}{r}\in L^{\infty}(\delta,1)$. Recall that the
embedding $H^1(\delta,1)\hookrightarrow C^0([\delta,1])$ is compact by 
Rellich's theorem. Therefore we obtain a locally uniformly convergent
subsequence, which we again denote by $f_n$.
As a consequence of the decay estimates (cf.~Corollary \ref{CorDecay}) we have
$$
|f_n(r)| \leq C(b)r^{\frac{1}{2}}(1+|\log r|^{\frac{1}{2}})\|f_n\|_{P_b} \leq C'(b)r^{\frac{1}{2}}(1+|\log r|^{\frac{1}{2}})
$$
if $b\not\in(-\frac{1}{2},\frac{1}{2})$, and
$$
|f_n(r)| \leq C(b)r^{-b}\|f_n\|_{P_b} \leq C'(b)r^{-b}
$$
if $b\in(-\frac{1}{2},\frac{1}{2})$. The functions $r^{\frac{1}{2}}(1+|\log
r|^{\frac{1}{2}})$ and $r^{-b}$ with $b<\frac{1}{2}$ are certainly in $L^2(0,1)$. In any case we
conclude with Lebesgue's dominated convergence theorem, that $f_n$ is convergent in $L^2(0,1)$. 
\end{pf}
\begin{cor}
The embedding $\dom (P^{\kappa}_b)_{\max}\hookrightarrow L^2(0,1)$ is compact for all $b\in\R$.
\end{cor}
\begin{pf}
This is a direct consequence of the previous lemma in view of Lemma \ref{radialequiv}.
\end{pf}\\
\\
For $b \in \R$ we define
$$
\widetilde{P}^{\kappa}_b = \left\{ \begin{array}{c@{\quad,\quad}c}
(P^{\kappa}_b)_{\max} & b \in (-\frac{1}{2},\frac{1}{2})\vspace{.2cm}\\
(P^{\kappa}_b)_{\min} & b \not\in (-\frac{1}{2},\frac{1}{2}) \end{array}\right.
$$
This determines a closed extension $\widetilde{P}^{\kappa}_B$ of $P^{\kappa}_B$ such that
$$
\dom \widetilde{P}^{\kappa}_B = \overline{\bigoplus_{b\in\spec B} \dom \widetilde{P}^{\kappa}_b \otimes \Psi_b}\,,
$$
where the closure is taken with respect to the graph norm $\|\cdot\|_{P^{\kappa}_B}$.
Note in particular that $\widetilde{P}^{\kappa}_B=(P_B^{\kappa})_{\min}$ if $\spec B \cap
(-\frac{1}{2},\frac{1}{2}) = \varnothing$. 
\begin{lemma}\label{ConeCompactness}
The embedding $\dom \widetilde{P}^{\kappa}_B\hookrightarrow L^2\left((0,1)\times N\right)$
is compact.
\end{lemma}
\begin{pf}
The previous lemma implies that for all $b
\in \spec B$ the embedding $(L_b^{\kappa})_{\max}: \dom(P_b^{\kappa})_{\max}\hookrightarrow L^2(0,1)$
is compact. 
We derive an upper bound for the operator norm of $(L_b^{\kappa})_{\min}:\dom
(P_b^{\kappa})_{\min} \hookrightarrow L^2(0,1)$, where
$\dom(P_b^{\kappa})_{\min}$ is equipped with the graph norm $\|\cdot\|_{P_b^{\kappa}}$. 
For $f
\in \Ccinf(0,1)$ we have
$$
P_b^{\kappa,t}P_b^{\kappa} f =-\frac{\partial^2f}{\partial r^2} + \frac{b(b+\cs(r))f}{\sn^2(r)}\,, 
$$
and therefore integration by parts applied twice yields
\begin{align*}
\left\|P_b^{\kappa}f\right\|^2_{L^2(0,1)} &= \left\langle P_b^{\kappa,t}P_b^{\kappa} f, f \right\rangle_{L^2(0,1)} \\
&= \int_0^1 \left\vert\frac {\partial f}{\partial r}\right\vert^2 +
\int_0^1 \frac{b(b+\cs(r))f^2}{\sn^2(r)}\\
&\geq C_{\kappa}(b) \left\|f\right\|^2_{L^2(0,1)}\,,
\end{align*} 
where $C_{\kappa}(b) \nearrow  \infty$ as $|b| \rightarrow \infty$.
Since $\Ccinf(0,1)$ is dense in $\dom (P_b^{\kappa})_{\min}$ we obtain
\begin{align*}
\left\|(L_b^{\kappa})_{\min}\right\|^2 &= \sup_{ f \in \Ccinf(0,1)\setminus\{0\} } \frac{\|f\|^2}{\|f\|^2 + \|P_b^{\kappa}f\|^2}\\
& \leq \frac {1}{1+C_{\kappa}(b)}\;,
\end{align*}
i.e.~for large eigenvalues of $B$ the operator norm of $(L_b^{\kappa})_{\min}$ is uniformly
small.

Let $L$ denote the embedding $\dom \widetilde{P}^{\kappa}_B\hookrightarrow
L^2\left((0,1)\times N\right)$. Furthermore for $a>0$
let $\pi^{<a}$ denote the projection onto the eigenspaces corresponding to
eigenvalues $b$ with $|b| < a$.  
Since there are only finitely many such eigenvalues,
$$L^{<a}=\pi^{< a}\circ L$$
is a compact operator and by the above estimates
$$\|L-L^{< a}\|^2 =
\sup_{|b|\geq a}\|(L_b^{\kappa})_{\min}\|^2  \leq \frac {1}{1+C_{\kappa}(a)}\,,$$
for $a$ large enough. 
In particular, for $a \rightarrow \infty$ we obtain that $L$ is a limit of compact operators
with respect to the operator norm
and is therefore itself compact. 
\end{pf}

\subsection{Spectral properties of cone-surfaces}

Let now $S$ be a cone-surface and $(\F,\nabla^{\F})$ a flat vector-bundle
over $N = \interior S$ equipped with a metric $h^{\F}$. Particular attention will be paid to
the spherical cone-surfaces 
$\Sph^2(\alpha,\beta,\gamma)$ and $\Sph^2(\alpha,\alpha)$, which appear as
links of singular points in a $3$-dimensional cone-manifold.

We wish to investigate spectral properties of the operators
$D(d_{\max})$ and $\Delta(d_{\max})$ by separation of variables.
In view of Lemma \ref{modeloperator} and Lemma \ref{max/min} the following requirements are natural:

\begin{Def}\label{surfaceconeadmissible}
Let $S$ be a cone-surface and $(\F,\nabla^{\F})$ a flat vector-bundle over
$N=\interior S$ equipped with a metric $h^{\F}$. 
If $\{x_i\}$ are the cone-points, we call $(\F,\nabla^{\F},h^{\F})$  {\em
  cone-admissible} if for all $i$:
\begin{enumerate}
\item Assumptions A1 and A2 hold for
  $(\F,\nabla^{\F},h^{\F})$ restricted to $U_{\varepsilon}(x_i)$, hence the
  model operator $P_{B_i}^{\kappa}$ is defined. 
\item $\spec B_i \cap
  (-\frac{1}{2},\frac{1}{2})=\varnothing$ holds.
\end{enumerate}
\end{Def}
\begin{rem}
Since the cross-section $S^1_{\alpha}$ is compact in this case, it would
  be enough to require A1 here, cf.~Remark \ref{A123}.
\end{rem}
We will see in the following that Definition \ref{surfaceconeadmissible} implicitly contains
  restrictions on the cone-angles of $S$ and the holonomy of the flat bundle
  $(\F,\nabla^{\F})$ around the cone-points:\\
\\
Let $S^1_{\alpha}=\R/\alpha \Z$ be the circle of length $\alpha$ and let $\cone_{\kappa,(0,\varepsilon)}S^1_{\alpha}$ be the
$\varepsilon$-truncated $\kappa$-cone over $S^1_{\alpha}$, i.e.~
$$
\cone_{\kappa,(0,\varepsilon)}S^1_{\alpha}=(0,\varepsilon)\times S^1_{\alpha}
$$
with metric
$$
dr^2+\sn^2(r)d\theta^2
$$
where $r \in (0,\varepsilon)$ and $\theta\in\R/\alpha\Z$.
Recall that if $x$ is a cone-point, the smooth part of the $\varepsilon$-ball around $x$ will be isometric
to
$$U_{\varepsilon}=\cone_{\kappa,(0,\varepsilon)}S^1_{\alpha}\,.$$ 
In this situation the model operator for the even part of the Hodge-Dirac
operator on the cone is given by
$$
P_B^{\kappa}=\frac{\partial}{\partial r} + \frac{1}{\sn(r)}B
$$
with
$$
B=D_{S^1_{\alpha}}+\left[ \begin{array}{cc} -\frac{1}{2} & \vspace{1mm}\\
                                                                       & -\frac{1}{2} \end{array} \right]=\left[ \begin{array}{cc} -\frac{1}{2} & d^t_{S^1_{\alpha}}\vspace{1mm}\\
                                                                      d_{S^1_{\alpha}} & -\frac{1}{2} \end{array} \right]\,.
$$
We determine the spectrum of the operator $B$, let us discuss the case with
trivial coefficient bundle first.
If we identify functions and $1$-forms on $S^1_{\alpha}$ via
\begin{align*}
\Cinf(S^1_\alpha) &\longrightarrow \Omega^1(S^1_{\alpha})\\
g &\longmapsto g\cdot d\theta\,,
\end{align*}
we may write
$$
D_{S^1_{\alpha}} = \left[ \begin{array}{cc} 0 & -\frac{\partial}{\partial \theta}\\
                                                                      \frac{\partial}{\partial \theta} & 0 \end{array} \right]\,.
$$
It is easily verified that
$$
\spec D_{S^1_{\alpha}}=\left\{\frac{2\pi n}{\alpha} , n \in \Z\right\},
$$
and therefore we obtain
$$
\spec B  =\left\{-\frac{1}{2}+\frac{2\pi n}{\alpha} , n \in \Z\right\}.
$$
We see that $\spec B \cap (-\frac{1}{2},\frac{1}{2})= \varnothing$ if
$\alpha\leq 2\pi$ in the case of trivial coefficients.

Let us now add a flat bundle to the situation. Let $\C(a)$ be the flat $U(1)$-bundle over $S^1_{\alpha}$
with holonomy $e^{ia}$, $a \in \R$. Without loss of generality we may assume that $a \in [0,2\pi)$.
Note that the bundles $\C(a)$ are topologically trivial. Any unitarily flat bundle on $S^1_{\alpha}$ decomposes
as a direct sum of these. A flat connection is given by
$$
\nabla^{\C(a)}=d -i\frac{a}{\alpha}d\theta.
$$
The associated Hodge-Dirac operator may be written as
$$
D_{S^1_{\alpha},\C(a)} = \left[ \begin{array}{cc} 0 & -\frac{\partial}{\partial \theta}+i\frac{a}{\alpha}\\
     \frac{\partial}{\partial \theta}-i\frac{a}{\alpha} & 0 \end{array}  \right]
$$
We obtain
$$
\spec  D_{S^1_{\alpha},\C(a)} = \left\{\pm \left|\frac{2\pi n-a}{\alpha}\right| , n \in \Z\right\},
$$
and therefore
$$
\spec B  =\left\{-\frac{1}{2}\pm\left|\frac{2\pi n-a}{\alpha}\right| , n \in \Z\right\}.
$$
We see that $\spec B \cap (-\frac{1}{2},\frac{1}{2})= \varnothing$ if either $a=0$ and $\alpha\leq 2\pi$ or
$\alpha\leq a \leq 2\pi - \alpha$. In the latter case we must in particular have
that $\alpha\leq\pi$.
\begin{rem}\label{remarksurfaceconeadmissible} The previous discussion shows
  that if $S$ has cone-angles $\leq \pi$ and $(\F,\nabla^{\F}, h^{\F})$ is an
  orthogonally flat bundle which decomposes locally around the cone-points as a direct sum of trivial
bundles $\R$ and bundles of type $\C(a)$  with $\alpha \leq a \leq
2\pi-\alpha$, then $(\F,\nabla^{\F}, h^{\F})$ will be cone-admissible in the sense of Definition \ref{surfaceconeadmissible}.
\end{rem}

\subsubsection{Discreteness}
In this section we investigate discreteness of the operators $D(d_{\max})$ and
$\Delta(d_{\max})$ on a cone-surface.
Recall that a selfadjoint operator $A$ is called discrete if
its spectrum is discrete, i.e.~if
$\spec A$ consists of a discrete set of eigenvalues with finite
multiplicities. A necessary and sufficient condition for $A$ to be discrete is the compactness of
the embedding $\dom A \hookrightarrow L^2$, where $\dom A$ is equipped with the
graph norm $\|\cdot\|_A$.

For simplicity we state the results concerning discreteness under the stronger
hypothesis that $(\F,\nabla^{\F},h^{\F})$ is cone-admissible, though we do not need the
assumption  $\spec B_i
\cap(-\frac{1}{2},\frac{1}{2})=\varnothing$ for $i \in \{1, \ldots, k\}$ as far as discreteness is concerned.
\begin{prop}\label{surfacediscreteness}
The embedding $\dom D^{ev}_{\max} \hookrightarrow
L^2(\Lambda^{ev}T^*N\otimes\F)$ is compact if $(\F,\nabla^{\F},h^{\F})$ is cone-admissible.
\end{prop}
\begin{pf}
We construct a partition of unity on $S$ in the
following way: Let $\{x_1,\ldots,x_k\}$ be the set of cone-points, we choose $\varepsilon >0$
such that the $U_{\varepsilon}(x_i)$ are disjoint.
We choose cut-off functions $\varphi_i$ supported inside $U_{\varepsilon}(x_i)$ with
$\varphi_i=\varphi_i(r)$ and $\varphi_i\equiv 1$ near the cone-point $x_i$. Then we define
$\varphi_{int}=1-\sum_{i=1}^k\varphi_i$.
Let $u_n \in \dom D^{ev}_{\max}$ be a sequence with
$\|u_n\|_{D^{ev}}\leq C$.

We claim that $\varphi_{int}u_n$ has a subsequence which is convergent in $L^2$: 
Let $\Omega \subset N$ be a relatively compact domain with smooth boundary, such
that $\supp \varphi_{int} \subset \Omega$. Then by the usual elliptic regularity
results, $\varphi_{int}u_n \in H^1_0(\Omega)$. Furthermore by the standard
elliptic estimate
$$
\|\varphi u_n\|^2_{H^1(\Omega)} \leq C \left(\,\|\varphi u_n\|_{L^2(\Omega)}^2 +
\|D^{ev}\varphi u_n\|^2_{L^2(\Omega)}\,\right) = C\,\|\varphi u_n\|^2_{D^{ev}_{\Omega}}\,.
$$
Now by Rellich's theorem
$H^1_0(\Omega)$ embeds into $L^2(\Omega)$ compactly, which proves the subclaim.

Thus we are reduced to a situation on the cone $U_{\varepsilon}=\cone_{\kappa,(0,\varepsilon)}S^1_{\alpha}$, i.e.~given a sequence
$f_n=\varphi u_n$ with
$\|f_n\|_{P^{\kappa}_B}\leq C$, we have to extract a
subsequence convergent in $L^2((0,1)\times S^1_{\alpha})$. The operator $B$ is essentially selfadjoint and discrete, since the cross-section of the cone is nonsingular in this case. Therefore the discussion from the last section applies.
It is a consequence of Corollary \ref{radial_max/min} that
$\varphi u_n \in \dom \widetilde{P}^{\kappa}_B$, therefore
Lemma \ref{ConeCompactness} yields the result. 
\end{pf}\\
\\
As a consequence we obtain that strong Hodge-decomposition holds for the
$d_{\max}$-complex on a cone-surface if $(\F, \nabla^{\F},
h^{\F})$ is cone-admissible. Here we remind the
reader of Theorem \ref{strongHodge} and the remark thereafter. 

We summarize the results concerning
Hodge-decomposition on a cone-surface relevant to $L^2$-cohomology in the following statement:
\begin{thm}[Hodge-theorem for cone-surfaces]\label{HodgeSurface} If $S$ is a
  cone-surface and
$(\F,\nabla^{\F})$ a flat vector-bundle over $N=\interior S$ together with a
  metric $h^{\F}$ such that $(\F,\nabla^{\F},h^{\F})$ is cone-admissible, then there is an orthogonal decomposition
$$
L^2(\Lambda^iT^*N\otimes\F)=\mathcal{H}^i_{\max}\oplus\im
d^{i-1}_{\max}\oplus\im (d^i)^t_{\min},
$$
and furthermore $\iota: \mathcal{H}^i_{\max} \rightarrow {H}^i_{\max}$ is an
isomorphism. The inclusion of the smooth $L^2$-complex $\Omega_{L^2}^i(N,\F) \rightarrow \dom d^i_{\max}$ induces
an isomorphism $H^i_{L^2}(N,\F) \cong  H^i_{\max}$.     
\end{thm}
Since $D^{odd}=(D^{ev})^t$, the same arguments
yield that the embeddings $\dom D^{odd}_{\max} \hookrightarrow
L^2(\Lambda^{odd}T^*N\otimes\F)$ and $\dom D_{\max} \hookrightarrow
L^2(\Lambda^{\bullet}T^*N\otimes\F)$ are again compact. 
\begin{prop}
The operators $D(d_{\max})$ and $\Delta(d_{\max})$ are
discrete on a cone-surface if $(\F,\nabla^{\F},h^{\F})$ is cone-admissible.
\end{prop}
\begin{pf}
Since $\dom D(d_{\max})$ and
$\dom\Delta(d_{\max})$ are continuously contained in
$\dom D_{\max}$, this follows from compactness of the embedding $\dom D_{\max} \hookrightarrow
L^2(\Lambda^{\bullet}T^*N\otimes\F)$. 
\end{pf}

\subsubsection{Selfadjointness}
In this section we establish essential selfadjointness of the Hodge-Dirac
operator $D$ on a cone-surface if $(\F, \nabla^{\F}, h^{\F})$ is cone-admissible. In contrast to the the previous section, we will now make strong use of the
assumption $\spec B_i\cap(-\frac{1}{2},\frac{1}{2})=\varnothing$ for $i \in
\{1, \ldots, k\}$ in Definition \ref{surfaceconeadmissible} since we
wish to apply Lemma \ref{max/min}.

\begin{prop}\label{turnover/spindle} 
$D^{ev}_{max}=D^{ev}_{min}$ on a cone-surface if $(\F,\nabla^{\F},h^{\F})$ is cone-admissible.
\end{prop}
\begin{pf}
Given $u\in\dom D^{ev}_{\max}$ we have to show that
$u\in\dom D^{ev}_{\min}$.
We choose a partition of unity on $S$ as in the proof of Proposition \ref{surfacediscreteness}.

We claim that $\varphi_{int}u\in\dom D^{ev}_{\min}$: As we have already observed
in the proof of Proposition \ref{surfacediscreteness}, if $\Omega\subset N$
is a relatively compact domain with smooth boundary such that $\supp
\varphi_{int} \subset \Omega$, then $\varphi_{int}u \in
H^1_0(\Omega)$. Now $\Ccinf(\Omega)$ is dense in $H^1_0(\Omega)$, therefore
we find a sequence $f_n \in \Ccinf(\Omega)$ such that $f_n$ approximates
$f=\varphi_{int}u$ with respect to the $H^1$-norm. But since $D^{ev}$
maps $H^1(\Omega)$ continuously to $L^2(\Omega)$, $f_n$ approximates $f$ also
with respect to the graph norm of $D^{ev}$, which proves the subclaim.

It remains to prove that $\varphi_i u \in \dom D^{ev}_{\min}$ for $i
\in\{1,\ldots,k\}$. But here we are again in a situation on the cone $U_{\varepsilon}=\cone_{\kappa,(0,\varepsilon)}S^1_{\alpha}$.
It is therefore
sufficient to show that $f=\varphi
u\in\dom(P^{\kappa}_B)_{\min}$ for
$u\in\dom(P^{\kappa}_B)_{\max}$ and $\varphi$ a cut-off
function of the above type. Now since $(\F,\nabla^{\F},h^{\F})$ is cone-admissible, $\spec
B\cap(-\frac{1}{2},\frac{1}{2})=\varnothing$ will be satisfied. Then 
Lemma \ref{max/min} implies that $f\in\dom(P^{\kappa}_B)_{\min}$, hence
in $\dom D^{ev}_{\min}$.
\end{pf}
\begin{cor}
The operator $D$ is essentially selfadjoint on a cone-surface if $(\F,\nabla^{\F},h^{\F})$ is
cone-admissible.
\end{cor}
\begin{pf}
We have
$$
D= \left[ \begin{array}{cc} 0 & (D^{ev})^t\\
                       D^{ev} & 0 \end{array} \right]
$$
considered as an operator
$$
\Omega_{\cp}^{ev}({N,\F})\oplus\Omega_{\cp}^{odd}({N,\F})\longrightarrow\Omega_{\cp}^{ev}({N,\F})\oplus\Omega_{\cp}^{odd}(N,{\F})
$$
and therefore
$$
D_{\min}= \left[ \begin{array}{cc} 0 & (D^{ev})_{\min}^t\\
                       D^{ev}_{\min} & 0 \end{array} \right]
$$
and
$$
D_{\max}= \left[ \begin{array}{cc} 0 & (D^{ev})_{\max}^t\\
                       D^{ev}_{\max} & 0 \end{array} \right]\,.
$$
This shows that $D_{\max}=D_{\min}$, i.e.~$D$ is
essentially selfadjoint.
\end{pf}
\begin{cor}
$\Delta_{\Fr}=\Delta(d_{max})$ on a cone-surface if $(\F,\nabla^{\F},h^{\F})$ is cone-admissible.
\end{cor}
\begin{pf}
This follows from essential selfadjointness of $D$ together with Corollary
\ref{CorDeltaFriedrichs}.
\end{pf}

\subsubsection{The first eigenvalue}
Let $\lambda_1$ be the smallest positive eigenvalue of
$\Delta^0(d_{\max})$ on the smooth part of $\Sph^2(\alpha,\beta,\gamma)$ (resp. $\Sph^2(\alpha,\alpha)$) with
coefficients in a flat vector-bundle $(\F,\nabla^{\F})$. Here we will derive a lower bound on $\lambda_1$, which will
be sufficient for later purposes.
Comparison with the smooth case suggests
that this bound might not be optimal.

\begin{prop}\label{firsteigenvalue}
Let $S$ be either $\Sph^2(\alpha,\beta,\gamma)$ or $\Sph^2(\alpha,\alpha)$ and $(\F,\nabla^{\F})$
a flat vector-bundle over $N=\interior S$ equipped with a metric $h^{\F}$. If $(\F,\nabla^{\F},h^{\F})$ is orthogonally flat and cone-admissible, then $\mathcal{H}^1_{\max}=0$. Moreover, under the same hypothesis, if $\lambda_1$ denotes the smallest positive eigenvalue of $\Delta^0(d_{\max})$, then $\lambda_1 \geq 1$.
\end{prop} 
\begin{pf}
Since $(\F,\nabla^{\F},h^{\F})$ is orthogonally flat, we may apply the standard Weitzenb\"ock
formula on $\F$-valued 1-forms
$$
\Delta\omega=\nabla^t\nabla\omega+(\Ric\varotimes\id)\omega\,,
$$
where the action of the Ricci tensor on a scalar-valued 1-form $\alpha$ is determined by the relation
$$
g(\Ric(\alpha),\beta)=\Ric(\alpha,\beta)
$$
for all $\beta\in\Omega^1(N,\R)$.
In two dimensions the Ricci tensor of a spherical metric (i.e.~of constant curvature
$\kappa=1$) is given by
$$
\Ric(\cdot\,,\cdot)=g(\cdot\,,\cdot),
$$
so we end up with
$$
\Delta\omega=\nabla^t\nabla\omega+\omega.
$$
For $\omega\in\Omega^1_{\cp}(N,\F)$ integration by parts yields
\begin{align*}
\int_N (\Delta\omega,\omega)
&=\int_N (\nabla^t\nabla\omega,\omega)+ \int_N \vert\omega \vert^2\\
&=\int_N \vert \nabla\omega\vert^2 + \int_N  \vert\omega \vert^2
\geq  \int_N  \vert\omega \vert^2\,.
\end{align*}
This means we have a lower bound for $\Delta$ on $\Omega^1_{\cp}(N,\F)$:
$$
\langle \Delta\omega, \omega\rangle_{L^2}\geq \|\omega\|^2_{L^2}\,.
$$
Since $\Delta(d_{\max})=\Delta_{\Fr}$ if $(\F,\nabla^{\F},h^{\F})$ is cone-admissible and the Friedrichs extension
preserves lower bounds, we obtain
$$
\langle \Delta(d_{\max})\omega,\omega\rangle_{L^2} \geq \|\omega\|^2_{L^2} 
$$
for all $\omega \in \dom \Delta^1(d_{\max})$. This proves the first part of the
assertion. 
Now for $f\in E_{\lambda_1}$, the
$\lambda_1$-eigenspace of $\Delta^0(d_{\max})$, $f\neq 0$, let $\omega=d_{\max}f$. Then
$w\neq 0$ and
$\Delta^1(d_{\max})\omega=d_{\max}d^t_{\min}d_{\max}f=\lambda_1\omega$.
This yields the estimate $\lambda_1\geq 1$.
\end{pf}

\subsection{Spectral properties of cone-3-manifolds}
Let in the following $C$ be a cone-3-manifold and $(\E,\nabla^{\E})$ a flat vector-bundle over
$M=\interior C$ equipped with a metric $h^{\E}$. Again we wish to investigate
spectral properties of the operators 
$D(d_{\max})$ and $\Delta(d_{\max})$ by separation of variables. We require:
\begin{Def}\label{3manifoldconeadmissible}
Let $C$ be a 3-dimensional cone-manifold and $(\E,\nabla^{\E})$ a flat vector-bundle over
$M=\interior C$ equipped with a metric $h^{\E}$. We call
$(\E,\nabla^{\E},h^{\E})$ {\em cone-admissible} if for all $x\in\Sigma$: 
\begin{enumerate}
\item Assumptions A1 and A2 hold for
  $(\E,\nabla^{\E},h^{\E})$ restricted to $U_{\varepsilon}(x)$, hence the
  model operator $P_{B_x}^{\kappa}$ is defined. 
\item $B_x$ is essentially selfadjoint and $\spec B_x \cap (-\frac{1}{2},\frac{1}{2})=\varnothing$ holds.
\end{enumerate}
\end{Def}
\begin{rem}
If we compare this definition with the cone-surface case, we note that a new
issue arises, namely that we have to include essential selfadjointness of the operator $B$ on the
cross-section of the model cone into the definition. This issue was not present
in the cone-surface case, since there the cross-section of the model cone was
compact.
\end{rem}
Let $x \in \Sigma$ be a singular point. For the local analysis around $x$ we consider two cases:
\begin{enumerate}
\item $x$ is a vertex
\item $x$ lies on a singular edge
\end{enumerate}
In the first case, the smooth part of the $\varepsilon$-ball around $x$ will be isometric to
$$
U_{\varepsilon}=\cone_{\kappa,(0,\varepsilon)}\interior\Sph^2(\alpha,\beta,\gamma) \,,
$$
and in the second case to
$$
U_{\varepsilon}=\cone_{\kappa,(0,\varepsilon)}\interior\Sph^2(\alpha,\alpha) \,.
$$
We treat the two cases simultaneously, let $N$ denote either
$\interior\Sph^2(\alpha,\beta,\gamma)$ or $\interior\Sph^2(\alpha,\alpha)$ in the
following.

Suppose that $(\E,\nabla^{\E},h^{\E})$ satisfies assumptions A1 and A2 on
$U_{\varepsilon}$, in particular $h_0^{\E}=\lim_{r \rightarrow 0}h^{\E}(r)$
exists and is parallel with respect to $\nabla^{\E}$.
Recall that the model operator for the even part of the Hodge-Dirac operator
on the $\kappa$-cone with two-dimensional cross-section $N$
is given by
$$
P_B^{\kappa} = \frac{\partial}{\partial r} + \frac{1}{\sn(r)} B
$$
with
\begin{align*}
 B &=D_{N}+\left[ \begin{array}{ccc}
 -1 & & \\
 & 0 &  \\
 & & 1 \end{array} \right]
 = \left[ \begin{array}{ccc} -1 & d^t_{N} &  \\
 d_{N} & 0 & d^t_{N} \\
  & d_{N} &  1 \end{array} \right]\,.
\end{align*}
Let us now assume that $(\E,\nabla^{\E},h_0^{E})$ restricted to the
2-dimensional cross-section $N$ is cone-admissible. Then $D_N$ and in particular the operator $B$ will be essentially selfadjoint. The Hodge-$\star$-operator defines a linear isometry
$$
\star: L^2(\Lambda^p T^*N\otimes\E) \longrightarrow L^2(\Lambda^{n-p} T^*N\otimes\E)\,, 
$$
where in this case $n=2$. Note furthermore that these two conditions together
imply that $\mathcal{H}^1_{\max}=0$ via Proposition \ref{firsteigenvalue}.

We determine the spectrum of $B$ in the following.
For $\lambda \geq 0$ let $E_{\lambda}$ be the $\lambda$-eigenspace of  
$$\Delta(d_{\max})=\Delta^0(d_{\max})\varoplus\Delta^1(d_{\max})\varoplus\Delta^2(d_{\max})\,.$$
Let $\lambda>0$ be an eigenvalue and $f_{\lambda}$ a corresponding eigensection of $\Delta^0(d_{max})$ with $\| f_{\lambda} \|_{L^2}$=1.
Then $$\left\{f_{\lambda}, \frac{1}{\sqrt{\lambda}}df_{\lambda},\frac{1}{\sqrt{\lambda}}\star df_{\lambda}, \star f_{\lambda}\right\}$$ form an orthonormal basis of a $D_N$-invariant subspace
$E_{f_{\lambda}}\subset E_{\lambda}$.
It is a consequence of Theorem \ref{HodgeSurface} that the $E_{f_{\lambda}}$ provide an orthogonal decomposition of $E_{\lambda}$ for $f_{\lambda}$ pairwise orthogonal. 
With respect to the given basis of $E_{f_{\lambda}}$ we have
$$
\left.D_{N}\right|_{E_{f_{\lambda}}} =  \left[ \begin{array}{cccc}
0 & \sqrt{\lambda} &  &  \\
\sqrt{\lambda} & 0 &  &  \\
 &  & 0 & -\sqrt{\lambda} \\
 &  & -\sqrt{\lambda} & 0 \end{array}\right]
$$
and correspondingly
$$
\left.B\right|_{E_{f_{\lambda}}} =  \left[ \begin{array}{cccc}
-1 & \sqrt{\lambda} &  &  \\
\sqrt{\lambda} & 0 &  &  \\
 &  & 0 & -\sqrt{\lambda} \\
 &  & -\sqrt{\lambda} & 1 \end{array}\right].
$$
For $\lambda=0$ we observe that if there is $f_0\in\mathcal{H}^0_{\max}$ with $\|f_0\|_{L^2}=1$, then $\{f_0, f_0\varotimes dvol\}$ form an orthonormal basis of $E_{f_0}\subset E_0=\mathcal{H}^0_{\max}\oplus\mathcal{H}^2_{\max}$
and we obtain
$$
\left.B\right|_{E_{f_0}} =  \left[ \begin{array}{cc}
-1 &  \\
 & 1 \end{array}\right].
$$
Note that possibly $E_0=0$. Therefore we obtain for the spectrum of $B$
$$
\spec B\subset\left\{-1, 1\right\} \cup \left\{\left. \pm\frac{1}{2}\pm \sqrt{\frac{1}{4}+\lambda} \right\vert \lambda \in \spec \Delta^0(d_{\max}), \lambda > 0\right\}.
$$
We see that $\spec B\cap(-\frac{1}{2}, \frac{1}{2})=\varnothing$ if $\lambda_1 \geq \frac{3}{4}$, which we can guarantee under the given conditions by means of Proposition \ref{firsteigenvalue}.

\begin{rem}\label{remark3manifoldconeadmissible}
As a consequence of the previous discussion we observe that a sufficient
condition for $(\E,\nabla^{\E},h^{\E})$ to be cone-admissible in the sense of
Definition \ref{3manifoldconeadmissible} is that assumptions A1 and A2 hold and the restriction of
$(\E,\nabla^{\E},h_0^{\E})$ to the
link $S_x$ of a singular point $x$ is cone-admissible in the
sense of Definition \ref{surfaceconeadmissible} for all $x\in\Sigma$.
\end{rem}
\subsubsection{Discreteness}%
In this section we investigate discreteness of the operators $D(d_{\max})$ and
$\Delta(d_{\max})$ on
a $3$-dimensional cone-manifold. 

For simplicity we state the results concerning discreteness under the stronger
hypothesis that $(\E,\nabla^{\E},h^{\E})$ is cone-admissible, though we do not need the
assumption  $\spec B_x
\cap(-\frac{1}{2},\frac{1}{2})=\varnothing$ for $x \in \Sigma$ as far as discreteness is concerned.

\begin{prop}\label{3manifolddiscreteness}
The embedding $\dom D^{ev}_{\max} \hookrightarrow
L^2(\Lambda^{ev}T^*M\otimes\E)$ is compact if $(\E,\nabla^{\E},h^{\E})$ is cone-admissible.
\end{prop}
\begin{pf}
Since $\Sigma$ is compact we find finitely many $x_i \in \Sigma$ such that the
$B_{\varepsilon}(x_i)$ cover $\Sigma$. Then $\{M, B_{\varepsilon}(x_i) \}$ is a finite
open cover of $C$. We fix a partition of unity $\{\varphi_{int}, \varphi_i \}$
subordinate to this cover. Let $U_{\varepsilon}(x_i)=B_{\varepsilon}(x_i)\cap M$.

Now let $u_n \in \dom D^{ev}_{\max}$ be a sequence with
$\|u_n\|_{D^{ev}} \leq C$. Clearly $\varphi_{int}u_n$ has a 
subsequence convergent in $L^2$: This follows in the same way as in the
cone-surface case (cf.~the proof of Proposition \ref{surfacediscreteness}). 

On the other hand $U_{\varepsilon}(x)$ will be isometric to
$\cone_{\kappa,(0,\varepsilon)}\interior\Sph^2(\alpha,\beta,\gamma)$ if $x$ is a vertex
or $\cone_{\kappa,(0,\varepsilon)}\interior\Sph^2(\alpha,\alpha)$ if $x$ is an edge point. Thus
we are reduced to a situation on the cone $U_{\varepsilon}=\cone_{\kappa,(0,\varepsilon)}N$. 
Without loss of generality we
may assume that $\varphi = \varphi(r)$ if $r$ is the radial variable and
$\varphi(r)= 1$ for $r$ small. 
If this is not the case we just replace $\varphi$
by a second
cut-off function $\widetilde{\varphi}\in\Ccinf(U_{\varepsilon}(x))$ which
satisfies these assumptions and in addition $\widetilde{\varphi} = 1$ near
$\supp \varphi$, and we replace $u_n$ by $\widetilde{u}_n = \varphi u_n$.
Since $(\E,\nabla^{\E},h^{\E})$ is cone-admissible, the operator $B$ will be essentially selfadjoint. $B$ will have discrete spectrum as a consequence of Proposition \ref{surfacediscreteness}. As in the cone-surface case we obtain that $\varphi u_n \in \dom \widetilde{P}_B^{\kappa}$.
We may now use Lemma \ref{ConeCompactness} to conclude the result.
\end{pf}\\
\\
We obtain that strong Hodge-decomposition holds for
the $d_{\max}$-complex on a cone-3-manifold if $(\E,\nabla^{\E},h^{\E})$ is cone-admissible.

We summarize the results concerning
Hodge-decomposition on a 3-dimensional cone-manifold relevant to $L^2$-cohomology in the following statement:
\begin{thm}[Hodge-theorem for cone-3-manifolds]\label{HodgeConeManifold} If $C$ is a cone-3-manifold and
$(\E,\nabla^{\E})$ a flat vector-bundle over $M=\interior C$ together with a
  metric $h^{\E}$ such that $(\E,\nabla^{\E},h^{\E})$ is cone-admissible, then there is an orthogonal decomposition
$$
L^2(\Lambda^iT^*M\otimes\E)=\mathcal{H}^i_{\max}\oplus\im
d^{i-1}_{\max}\oplus\im (d^i)^t_{\min},
$$
and furthermore $\iota: \mathcal{H}^i_{\max} \rightarrow {H}^i_{\max}$ is an
isomorphism. The inclusion of the smooth $L^2$-complex $\Omega_{L^2}^i(M,\E) \rightarrow \dom d^i_{\max}$ induces
an isomorphism $H^i_{L^2}(M,\E) \cong  H^i_{\max}$.     
\end{thm}
Since $D^{odd}=(D^{ev})^t$, the same arguments
yield that the embeddings $\dom D^{odd}_{\max} \hookrightarrow
L^2(\Lambda^{odd}T^*M\otimes\E)$ and $\dom D_{\max} \hookrightarrow
L^2(\Lambda^{\bullet}T^*M\otimes\E)$ are again compact. 
\begin{prop}
The operators $D(d_{\max})$ and $\Delta(d_{\max})$ are
discrete on a cone-3-manifold if $(\E,\nabla^{\E},h^{\E})$ is cone-admissible.
\end{prop}
\begin{pf}
Since $\dom D(d_{\max})$ and
$\dom\Delta(d_{\max})$ are continuously contained in
$\dom D_{\max}$, this follows from compactness of the embedding $\dom D_{\max} \hookrightarrow
L^2(\Lambda^{\bullet}T^*M\otimes\E)$. 
\end{pf}

\subsubsection{Selfadjointness}%

In this section we establish essential selfadjointness of the
Hodge-Dirac operator $D$ on a
cone-3-manifold if
$(\E,\nabla^{\E},h^{\E})$ is cone-admissible. Here the condition $\spec B_x
\cap(-\frac{1}{2},\frac{1}{2})=\varnothing$ for all  $x \in \Sigma$ is essential. 
\begin{prop}\label{C} 
$D^{ev}_{max}=D^{ev}_{min}$ on a cone-3-manifold if $(\E,\nabla^{\E},h^{\E})$ is cone-admissible.
\end{prop}
\begin{pf}
Given $u\in\dom D^{ev}_{\max}$ we have to show that
$u\in\dom D^{ev}_{\min}$.
We choose a partition of unity on $C$ as in the proof of Proposition
\ref{3manifolddiscreteness}.

Clearly
$\varphi_{int}u\in\dom D^{ev}_{\min}$: This follows in the same way
as in the cone-surface case (cf.~the proof of Proposition \ref{turnover/spindle}). 

It remains to prove that $\varphi_i u \in \dom D^{ev}_{\min}$. Again this brings us back to a situation on the cone $U_{\varepsilon}=\cone_{\kappa,(0,\varepsilon)}N$, where $N=\interior\Sph^2(\alpha,\beta,\gamma)$ or $N=\interior\Sph^2(\alpha,\alpha)$.
It is therefore
sufficient to show that $f=\varphi
u\in\dom(P^{\kappa}_B)_{\min}$ for
$u\in\dom(P^{\kappa}_B)_{\max}$ and $\varphi$ a cut-off
function of the above type. Since $(\E,\nabla^{\E},h^{\E})$ is cone-admissible, $B$ is essentially selfadjoint and has discrete spectrum. Moreover, the condition $\spec
B\cap(-\frac{1}{2},\frac{1}{2})=\varnothing$ will be satisfied. Then 
Lemma \ref{max/min} implies that $f\in\dom(P^{\kappa}_B)_{\min}$, hence
in $\dom D^{ev}_{\min}$.
\end{pf}
\begin{cor}\label{essentialselfadjointness}
The operator $D$ is essentially selfadjoint on a cone-3-manifold if $(\E,\nabla^{\E},h^{\E})$ is cone-admissible.
\end{cor}
\begin{pf}
This follows from $D^{ev}_{\max}=D^{ev}_{\min}$ in the same way as in the cone-surface case.
\end{pf}
\begin{cor}\label{Friedrichs=dmax}
$\Delta_{\Fr}=\Delta(d_{max})$ on a cone-3-manifold if $(\E,\nabla^{\E},h^{\E})$ is cone-admissible.
\end{cor}
\begin{pf}
This follows from essential selfadjointness of $D$ together with Corollary
\ref{CorDeltaFriedrichs}.
\end{pf}


\section{The Bochner technique}


\subsection{Infinitesimal isometries}
For simplicity consider $\M$ for $\kappa \in \{-1, 0,1\}$. Let $G=\Isom\M$ and
$\mathfrak{g}$ its Lie-algebra. $\mathfrak{g}$ may be identified with the Lie-algebra of
Killing vectorfields. Note however, that the Lie-bracket in $\mathfrak{g}$ corresponds to
the negative of the vectorfield commutator under this identification:
$$
ad_{\mathfrak{g}}(X)Y=[X,Y]_{\mathfrak{g}}=-[X,Y]=-\mathcal{L}_XY.
$$  
Fix a point $p\in\M$ and let $K=\Stab_G(p)$. Note that $K\cong\SO(T_p\M$), since $G$ acts
simply transitively on frames in constant curvature. Then we have
the usual decomposition
$\mathfrak{g} = \mathfrak{k}\oplus\mathfrak{p}$, where $\mathfrak{k}$ is the Lie-algebra of $K$.
Recall that
$$
\mathfrak{k}=\{X\in\g \,\vert\, X(p)=0\}
$$
and
$$
\mathfrak{p}=\{X\in\g \,\vert\, (\nabla X)(p)=0\}.
$$
There are isomorphisms
$$
\mathfrak{p} \cong T_p\M,\,
X \mapsto X(p)
$$
and (in our constant-curvature situation)
$$
\mathfrak{k} \cong \mathfrak{so}(T_p\M) ,\, X \mapsto A_X(p):=(\nabla X)(p) \,.
$$
We have
$[\mathfrak{k},\mathfrak{k}]\subset\mathfrak{k}$,
$[\mathfrak{k},\mathfrak{p}]\subset\mathfrak{p}$ and
$[\mathfrak{p},\mathfrak{p}]\subset\mathfrak{k}$, since $\mathfrak{k}$
(resp.~$\mathfrak{p}$) is the $+1$ (resp.~$-1$) eigenspace of the Cartan-involution on
$\mathfrak{g}$ induced by the geodesic involution on $\M$ about $p$.

Let $X$ be a Killing vectorfield. Let $\gamma$ be a geodesic with $\gamma(0)=p$
and $\dot{\gamma}(0)=Y(p)$. Then $X$ will be a Jacobi vectorfield along $\gamma$.
We obtain
\begin{align*}
0&=\nabla_{\dot{\gamma}}\nabla_{\dot{\gamma}}X+R(X,\dot{\gamma})\dot{\gamma}\\
&=(\nabla_{\dot{\gamma}}A_X)\dot{\gamma}+R(X,\dot{\gamma})\dot{\gamma}
\end{align*}
Therefore we have
$$
(\nabla_YA_X)Y+R(X,Y)Y=0\,.
$$
The expression $(\nabla_YA_X)Z+R(X,Y)Z$ is symmetric in $Y$ and
$Z$. 
Therefore we obtain by polarization
$$
(\nabla_YA_X)Z+R(X,Y)Z=0\quad\quad\quad\quad(\bigstar)
$$
if $X$ is a Killing vectorfield.
\begin{lemma}
Under the identification $\mathfrak{g}=\mathfrak{so}(T_p\M)\oplus T_p\M$
the Lie-bracket corresponds to 
$$
[(A,X),(B,Y)]=([A,B]-R(X,Y),AY-BX),
$$ 
where $[A,B]$ is the commutator in $\mathfrak{so}(T_p\M)$ and $R$ the Riemannian
curvature tensor.
\end{lemma}\label{Lie}
\begin{pf}
Let $X,Y\in\mathfrak{p}, Z\in\mathfrak{p}$. From equation $(\bigstar)$ we obtain \begin{align*}
A_{[X,Y]_{\g}}Z(p)&=-\nabla_Z([X,Y])(p)=-\nabla_Z\nabla_XY(p)+\nabla_Z\nabla_YX(p)\\
&=-(\nabla_Z A_Y)X(p)+(\nabla_Z A_X)Y(p)\\
&=R(Y,Z)X(p)+R(Z,X)Y(p)=-R(X,Y)Z(p)
\end{align*}
Let $X,Y\in\mathfrak{k}$, $Z\in\mathfrak{p}$.
\begin{align*}
A_{[X,Y]_{\g}}Z(p) &=-\nabla_Z([X,Y])(p)=-\nabla_{[X,Y]}Z(p)-[Z,[X,Y]](p)\\
& = [X,[Y,Z]](p)+[Y,[Z,X]](p)\\
&=[X,\nabla_YZ-\nabla_ZY](p)+[Y,\nabla_ZX-\nabla_XZ](p)\\
&=\nabla_{\nabla_ZY}X(p)-\nabla_{\nabla_ZX}Y(p)=[A_X,A_Y]Z(p)
\end{align*}
Let $X\in\mathfrak{k}, Y\in\mathfrak{p}$.
\begin{align*}
[X,Y]_{\g}(p) &= -\left(\nabla_XY-\nabla_YX\right)(p)\\
&=\nabla_YX(p)=A_XY(p)
\end{align*}
This is sufficient since
$[\mathfrak{k},\mathfrak{k}]\subset\mathfrak{k}$,
$[\mathfrak{k},\mathfrak{p}]\subset\mathfrak{p}$ and $[\mathfrak{p},\mathfrak{p}]\subset\mathfrak{k}$.
\end{pf}\\
\\
Note that the usual formula for the curvature tensor of a symmetric space
$$
R(X,Y)Z(p)=-[[X,Y],Z](p),\;X,Y,Z\in\mathfrak{p}
$$
is contained in the statement.
\begin{cor}\label{Ad}
$Ad_G(g)(A,X)=(Ad_K(g)A,gX)$ for $g\in K$.
\end{cor}
Let $\E=\mathfrak{so}(T\M)\oplus T\M$. $\E$ is a bundle of
Lie-algebras with a flat connection $\nabla^{\E}$,
such that a section $\sigma=(A,X)$ is parallel if and only if $X$ is a
Killing vectorfield and $A=A_X$. 
\begin{lemma}\label{flat}
The flat connection on $\E$ is given by
$$
\nabla_Y^{\E}(A,X)=(\nabla_YA-R(Y,X),\nabla_YX-AY) \,,
$$
where $\nabla$ denotes the Levi-Civita connection on $T\M$ and on $\mathfrak{so}(T\M)$.
\end{lemma}
\begin{pf}
If $\nabla^0$ and $\nabla^1$ are connections on a vector-bundle
$\E$, then the difference $\alpha=\nabla^0-\nabla^1$ is a 1-form with values in
$\End \E$. If $\nabla^0\sigma=0$, then $-\nabla_Y^1\sigma=\alpha(Y)\sigma$ for
all $Y\in T\M$.

Let $\nabla^0=\nabla^{\E}$ and $\nabla^1=\nabla$. A Killing vectorfield $X$ determines a parallel
section $\sigma_X=(A_X,X)$. 
From equation $(\bigstar)$ we have
$$
(\nabla_YA_X)Z=-R(X,Y)Z=R(Y,X)Z \,,
$$
and from the very definition 
$$\nabla_YX=A_XY ,$$
hence $\alpha(Y)(A,X)=(-R(Y,X),-AY)$.
\end{pf}\\
\\
In fact $\E=\M\times\mathfrak{g}$ and $\nabla^{\E}$ is just the trivial connection $d$
written in terms of the subbundles $T\M$ and $\mathfrak{so}(T\M)$.
\begin{cor}$\nabla_Y^{\E}\sigma=\nabla_Y\sigma+ad_{\g}(Y)\sigma$ for $\sigma\in\Gamma(\E),
Y\in T\M$.
\end{cor}
\begin{pf}
Lemma \ref{Lie} implies that $\alpha(Y)\sigma=ad(Y)\sigma$. 
\end{pf}\\
\\
We have a natural metric on $\E$, namely
$$
h^{\E}= (\,\cdot\,,\cdot\,)_{\mathfrak{so}(T\M)}\oplus(\,\cdot\,,\cdot\,)_{T\M},
$$
where
$$
(A,B)_{\mathfrak{so}(T\M)}=-\frac{1}{2}\tr(AB).
$$
Recall the definition of the Killing form
$$
B_{\g}(a,b)=\tr(ad_{\g}(a)ad_{\g}(b))
$$
for $a,b\in\g$. $B_{\g}$ is a symmetric bilinear form, which is
$Ad_G(g)$-invariant for all $g\in G$. This implies in particular that
$ad_{\g}(a)$ is antisymmetric with respect to $B_{\g}$ for all $a \in \g$.

We wish to express $B_{\g}$ in terms of the decomposition
$\g=\mathfrak{k}\oplus\mathfrak{p}$. First of course the relations $[\mathfrak{k},\mathfrak{k}]\subset\mathfrak{k}$,
$[\mathfrak{k},\mathfrak{p}]\subset\mathfrak{p}$ and
$[\mathfrak{p},\mathfrak{p}]\subset\mathfrak{k}$ imply that $\k$ and $\p$ are $B_{\g}$-orthogonal. The following computation is left to the reader:

\begin{lemma}\label{Killing} The restrictions of $B_{\g}$ to $\k=\mathfrak{so}(T_p\M)$ and
$\p= T_p\M$ are given as follows:
\begin{align*}
\left.B_{\g}\right\vert_{\k}(\,\cdot\,,\cdot\,)&=-4(\,\cdot\,,\cdot\,)_{\mathfrak{so}(T_p\M)}\\
\left.B_{\g}\right\vert_{\p}(\,\cdot\,,\cdot\,)&=-4\kappa(\,\cdot\,,\cdot\,)_{T_p\M}\,.
\end{align*}
\end{lemma}
We obtain as an immediate consequence:
\begin{cor}
If $\kappa=1$, then $ad(Y)$ is antisymmetric with respect to $h^{\E}$ for $Y\in T\M$, in
particular $\nabla^{\E}h^{\E}=0$. 
If $\kappa=-1$, then $ad(Y)$ is symmetric with respect to $h^{\E}$ for $Y\in T\M$.
\end{cor}
For $\kappa=-1$ we want
to calculate the precise deviation of $h^{\E}$ from being parallel. With $\nabla^{\E}=\nabla + ad$ we get using
the fact that $h^{\E}$ is parallel with respect to $\nabla$:
\begin{align*}
(\nabla^{\E}_X h^{\E})(\sigma,\tau) &= X(h^{\E}(\sigma,\tau)) -
h^{\E}(\nabla_X^{\E}\sigma,\tau) - h^{\E}(\sigma,\nabla_X^{\E}\tau)\\
&= -h^{\E}(ad(X)\sigma,\tau)-h^{\E}(\sigma,ad(X)\tau)\\
&= -2h^{\E}(ad(X)\sigma,\tau)\,.
\end{align*}
Let $h^{\E}_0$ denote the metric on $\E$ obtained by parallel extension of
$h^{\E}(p)$ with respect to $\nabla^{\E}$ for $p \in \M$. If we write $h^{\E}(\sigma,\tau) = h_0^{\E}(A\sigma,\tau)$ with $A\in
\Gamma(\End \E)$ symmetric, we obtain using the fact that $h^{\E}_0$ is
parallel with respect to $\nabla^{\E}$:
\begin{align*}
h_0^{\E}((\nabla_X^{\E}A)\sigma,\tau) &= (\nabla^{\E}_X h^{\E})(\sigma,\tau) = -2h^{\E}(ad(X)\sigma, \tau)\\
 & =-2h_0^{\E}(A \,ad(X)\sigma,\tau) \,,
\end{align*}
in particular we have proved:
\begin{lemma}\label{bounded}
If $\kappa=-1$, then $A^{-1}(\nabla^{\E}A) = -2 ad$ and is therefore
bounded on $\M$ with respect to $h^{\E}$.
\end{lemma}
Let us now consider $M$, the nonsingular part of a cone-3-manifold $C$. The condition
that $M$ is locally modelled on $\M$ is usually expressed in terms of the
{\it developing map} 
$$
\dev:(\widetilde{M},p_0) \longrightarrow (\M,p)
$$
and the {\it holonomy representation}
$$
\hol:\pi_1(M,x_0) \longrightarrow G=\Isom\M\,,
$$
where $\dev$ is a local isometry and $\pi_1(M)$-equivariant with respect to the
deck-action on $\widetilde{M}$ and the action via $\hol$ on $\M$. For details
we ask the reader to consult Section \ref{XGstructures}.

We again denote by $\E$ the bundle $\mathfrak{so}(TM)\oplus TM$. Since being a
Killing vectorfield is a local condition, we again have a flat connection
$\nabla^{\E}$ on $\E$ with the property that parallel sections correspond to
Killing vectorfields. The formula for $\nabla^{\E}$ given in Lemma \ref{flat}
applies as well. In contrast to the model-space situation, $\E$ will now have holonomy. It is easy
to see that the holonomy of $\E$ along a loop $\gamma\in\pi_1(M,x_0)$ is given by
$Ad\circ\hol(\gamma)$
if we identify $\E_{x_0}$ with $\mathfrak{g}$. Therefore we obtain an alternative
description of $\E$:
$$
\E=\widetilde{M} \times_{Ad\circ\hol} \mathfrak{g}
$$
The Lie-algebra structure on $\E$ induced by this representation coincides with the one given in
Lemma \ref{Lie}.\\
\\
The same considerations apply to the two-dimensional situation as
well if we replace $\M$ and its isometry group with the corresponding two-dimensional
objects. Here we restrict our attention to the spherical case. 
Let 
$$
S=\left\{ \begin{array}{c@{\quad}c}
\Sph^2(\alpha,\beta,\gamma) & \text{or}\vspace{2mm}\\
\Sph^2(\alpha,\alpha) &  \end{array} \right.
$$ 
in the following. Since $\Isom
\Sph^2=\SO(3)$ we have a holonomy representation
$$
\hol: \pi_1(\interior S) \longrightarrow \Isom\Sph^2=\SO(3) 
$$
and developing map
$$
\dev: \widetilde{\interior S} \longrightarrow \Sph^2\,.
$$
Let us denote the vector-bundle of infinitesimal isometries with its natural flat connection
in this situation
by $(\F,\nabla^{\F})$. We have
$$
\F=\widetilde{\interior S} \times_{Ad\circ\hol} \mathfrak{so}(3)\,.
$$
Since the adjoint representation of $\SO(3)$ on
$\mathfrak{so}(3)$ is isomorphic to the standard representation of $\SO(3)$ on $\R^3$, we
have alternatively 
$$
\F=\widetilde{\interior S} \times_{\hol} \R^3\,.
$$
Since $\hol$ preserves the standard scalar product on $\R^3$, we have a natural
metric $h^{\F}$ on $\F$ which is parallel with respect to $\nabla^{\F}$.

Now if $x_i \in S$ is a cone-point with cone-angle $\alpha_i$ and $\gamma_i \in
\pi_1(\interior S)$ a loop around $x_i$, then 
$\hol(\gamma_i)$ is just rotation about the cone-angle $\alpha_i$ around
some fixed axis in $\R^3$. Note that the axis of $\hol(\gamma_i)$ and the axis of
$\hol(\gamma_j)$ need not coincide for $x_i \neq x_j$. 
This gives us a quite explicit
description of $\F$. In particular we see that locally around the cone-points
we have the following splitting
$$
\F\vert_{S^1_{\alpha_i}}=\C(\alpha_i)\oplus \R\,,
$$
where $\C(\alpha_i)$ denotes the flat $U(1)$-bundle over $S^1_{\alpha_i}$
with holonomy $e^{i\alpha_i}$.\\
\\
Next we describe the restriction of $\E$ to the links of singular points. Recall
that if $x\in \Sigma$ is a singular point and $S_x$ is its link, then 
$$
S_x = \Sph^2(\alpha,\beta,\gamma) 
$$ 
if $x$ is a vertex, and
$$
S_x = \Sph^2(\alpha,\alpha) 
$$ 
if $x$ is an edge point.
\begin{lemma}\label{Esplitting}
Let $S_x$ be the link of a singular point $x\in\Sigma$. Then the restriction of
 $\E$ to $\interior S_x$ is given by:
$$
\left.\E\right\vert_{\interior S_x}=\F\oplus\F,
$$
where $\F$ is the flat vector-bundle of infinitesimal isometries on $S_x$.
\end{lemma}
\begin{pf}
The holonomy of $\pi_1(\interior S_x)$ fixes a point $p\in\M$ and is
therefore contained in $K=\Stab_G(p)\cong \SO(T_p\M)$. We have seen in Corollary \ref{Ad} that
$Ad_G(g)=(Ad_K(g),g)$ for $g\in K$ with respect to the splitting
$\mathfrak{g}=\mathfrak{k}\oplus\mathfrak{p}$. Again, since the adjoint
representation and the standard representation of $\SO(3)$ are isomorphic, we obtain two
copies of $\F$.
\end{pf}
\begin{prop}\label{Econeadmissible_spherical}
Let $C$ be a spherical cone-3-manifold with cone-angles
  $\leq \pi$. Then $(\E,\nabla^{\E},h^{\E})$, the vector-bundle of
  infinitesimal isometries of $M=\interior C$ with its natural flat connection
  and metric, is cone-admissible.
\end{prop}
\begin{pf}
Let $x \in \Sigma$ be a singular point, we have
$\E|_{U_{\varepsilon}(x)} = \F \oplus \F$ via Lemma \ref{Esplitting}. 
Since $h^{\E}$ is parallel in the spherical case, we have
  $h_0^{\E}=h^{\E}$ and assumptions A1 and A2 are trivially
  satisfied on $U_{\varepsilon}(x)$, cf.~Remark \ref{A123}. Clearly
  $h^{\E}=h^{\F} \oplus h^{\F}$. $(\F,\nabla^{\F},h^{\F})$ is orthogonally flat and therefore via Remark \ref{remarksurfaceconeadmissible} cone-admissible over $\interior S_x$ if the cone-angles are
$\leq \pi$. Then we may apply Remark \ref{remark3manifoldconeadmissible}
to conclude that $(\E,\nabla^{\E},h^{\E})$ is cone-admissible over $M$.
\end{pf} 
\begin{prop}\label{Econeadmissible_hyperbolic}
Let $C$ be a hyperbolic cone-3-manifold with cone-angles
  $\leq \pi$. Then $(\E,\nabla^{\E},h^{\E})$, the vector-bundle of
  infinitesimal isometries of $M=\interior C$ with its natural flat connection
  and metric, is cone-admissible.
\end{prop}
\begin{pf}
Let $x \in \Sigma$ be a singular point, we have
$\E|_{U_{\varepsilon}(x)} = \F \oplus \F$ via Lemma
\ref{Esplitting}. $(\F,\nabla^{\F},h^{\F})$ is orthogonally flat, clearly $h_0^{\E}=\lim_{r \rightarrow
  0}h^{\E}(r)$ exists and $h_0^{\E}=h^{\F} \oplus h^{\F}$, i.e.~assumption
A1 is satisfied. In view of Lemma \ref{bounded}, assumption A2 is also
satisfied. The assertion follows now as in the spherical case.
\end{pf}\\
\\
In the Euclidean case for fixed $p\in\Euc^3$ we have a group homomorphism
\begin{align*}
\rot: \Isom\Euc^3 &\longrightarrow \Stab_G(p)\cong \SO(T_p\Euc^3)\\ 
g &\longmapsto g+(p-g(p))
\end{align*}
We may form the rotational part of the holonomy 
$$
\rot\circ\hol: \pi_1(M)\longrightarrow  \Stab_G(p)\cong
\SO(T_p\Euc^3)\,.
$$
On the other hand 
$$
\E_{trans}:= TM \subset \E = \mathfrak{so}(TM)\oplus TM
$$
is via the explicit formula for $\nabla^{\E}$ in  Lemma \ref{flat}
easily seen to be a parallel subbundle of $\E$.
Note that in contrast 
$$
\E_{rot}:= \mathfrak{so}(TM) \subset \E = \mathfrak{so}(TM)\oplus TM
$$
is {\it not} parallel. 

Since the rotational part of the
holonomy is nothing but the
holonomy of the flat tangent bundle, we obtain
$$
\E_{trans} = \widetilde{M} \times_{\rot\circ\hol}\R^3\,.
$$
In the same way as before one shows:
\begin{lemma}\label{Etransconeadmissible}
Let $C$ be a Euclidean cone-3-manifold.
The restriction of $\E_{trans}$ to the link $S_x$ of a singular point $x \in \Sigma$ is given as
$$
\left.\E_{trans}\right\vert_{\interior S_x}=\F,
$$
where $\F$ is the flat vector-bundle of infinitesimal isometries on $S_x$.
Furthermore $\E_{trans}$ is cone-admissible if the cone-angles are $\leq \pi$.
\end{lemma}

\subsection{Weitzenb\"ock formulas}
\subsubsection{The spherical and the Euclidean case}
Let $(\E,\nabla^{\E},h^{\E})$ be an orthogonally flat vector-bundle. Recall the standard Weitzenb\"ock
formula on $\E$-valued 1-forms: 
$$
\Delta\omega=\nabla^t\nabla\omega+(\Ric\varotimes\id)\omega
$$
For this formula to hold without extra terms we really need that the metric $h^{\E}$
is parallel with respect to $\nabla^{\E}$.
Recall that the action of the Ricci tensor on a scalar-valued 1-form $\alpha$ is determined by the relation
$$
g(\Ric(\alpha),\beta)=\Ric(\alpha,\beta)
$$
for all $\beta\in\Omega^1(M,\R)$.
In three dimensions the Ricci tensor of a metric with constant curvature
$\kappa$ is given by
$$
\Ric(\cdot\,,\cdot)=2\kappa\cdot g(\cdot\,,\cdot),
$$
so we end up with\\
\\
($\kappa=1$)
$$
\Delta\omega=\nabla^t\nabla\omega+2\cdot\omega
$$
($\kappa=0$)
$$
\Delta\omega=\nabla^t\nabla\omega
$$
in the spherical and the Euclidean case.
\subsubsection{The hyperbolic case}
In the hyperbolic case we use a different type of Weitzenb\"ock formula, due to Y.~Matsushima and S.~Murakami, cf.~\cite{MM}. We use the notation of \cite{HK}. Let $\E=\mathfrak{so}(TM) \oplus TM$ be the vector-bundle of
infinitesimal isometries and $\nabla^{\E}$ its natural flat connection. We continue to denote by
$\nabla^{\E}$ the tensor-product connection on $\Lambda^{\bullet}T^*M\otimes\E$
induced by the Levi-Civita connection on $M$ and the connection $\nabla^{\E}$ on $\E$, whereas we denote by $\nabla$ the Levi-Civita connection on $\Lambda^{\bullet}T^*M\otimes\E$.

Recall the relation $\nabla^{\E}_Y=\nabla_Y+ad(Y)$ for
$Y\in TM$, where the endomorphism $ad(Y)$ is symmetric with respect to $h^{\E}$. Let in
the following
$$
\varepsilon : T^*M \otimes \Lambda^{\bullet}T^*M \rightarrow  \Lambda^{\bullet+1}T^*M
$$
denote exterior multiplication, and
$$
\iota :TM \otimes \Lambda^{\bullet}T^*M \rightarrow  \Lambda^{\bullet-1}T^*M
$$
denote interior multiplication. Let $\{e_1,e_2,e_3\}$ be a local orthonormal frame and $\{e^1,e^2,e^3\}$ the dual coframe.
Then we have 
$$
d=\sum_{i=1}^3\varepsilon(e^i)\nabla_{e_i}^{\E}=\sum_{i=1}^3\varepsilon(e^i)\left(\nabla_{e_i}+ad(e_i)\right).
$$ 
This implies
$$
d^t=-\sum_{i=1}^3\iota(e_i)\left(\nabla_{e_i}-ad(e_i)\right).
$$
Define
$$
\mathcal{D}:=\sum_{i=1}^3 \varepsilon(e^i)\nabla_{e_i}\;\text{and}\;T:=\sum_{i=1}^3\varepsilon(e^i)ad(e_i),
$$
this implies
$$
\mathcal{D}^t=-\sum_{i=1}^3 \iota(e^i)\nabla_{e_i}\;\text{and}\;T^t=\sum_{i=1}^3\iota(e_i)ad(e_i).
$$
We obviously have $d=\mathcal{D}+T$ and $d^t=\mathcal{D}^t+T^t$, let $\Delta_\mathcal{D}=\mathcal{D}\mathcal{D}^t+\mathcal{D}^t\mathcal{D}$ and
$H=TT^t+T^tT$. H is symmetric and non-negative. From the definitions we have
\begin{align*}
\Delta&=dd^t+d^td\\
&=\Delta_\mathcal{D}+H+\mathcal{D} T^t+T\mathcal{D}^t+\mathcal{D}^tT+T^t\mathcal{D}
\end{align*}
A computation in a local orthogonal frame shows that
$$\mathcal{D} T^t+T\mathcal{D}^t+\mathcal{D}^tT+T^t\mathcal{D}=0$$
and
 $$
H=\sum_{i=1}^3 ad(e_i)^2+\sum_{i,j=1}^3\varepsilon(e^i)\iota(e_j)ad\left([e_i,e_j]\right)\,.
$$
This implies the following Weitzenb\"ock formula, where a priori $\Delta_\mathcal{D}$ and $H$ are non-negative.
\begin{lemma}{\rm\cite{MM}} $\Delta=\Delta_\mathcal{D}+H$
\end{lemma}
The following positivity property of $H$ on $1$-forms makes this formula particularly useful for us. The proof may again be obtained by a calculation in a local orthonormal frame.
\begin{prop}{\rm \cite{MM}} There is a constant $C>0$ such that 
$$
(H\omega,\omega)_x\geq C (\omega,\omega)_x
$$
for all $\omega\in\Omega^1(M,\E)$ and $x\in M$.
\end{prop}


\subsection{A vanishing theorem}
In this section we prove our main result about $L^2$-cohomology
spaces of $3$-dimensional cone-manifolds with coefficients in the flat vector-bundle of infinitesimal
isometries.
This completes the analytic part of our argument. For convenience we discuss the proof case by case.

\subsubsection{The spherical case}
\begin{thm} Let $C$ be a spherical cone-3-manifold with cone-angles $\leq
\pi$. Let $M=C\setminus\Sigma$ and $(\E,\nabla^{\E})$ be the vector-bundle of
infinitesimal isometries of $M$ with its natural flat connection. Then
$$
H^1_{L^2}(M,\E)=0\,.
$$
\end{thm}
\begin{pf}
We recall the Weitzenb\"ock formula for the Hodge-Laplace operator on $\E$-valued $1$-forms, which in the spherical case (i.e.~$\kappa=1$) amounts to
$$
\Delta \omega = \nabla^t\nabla \omega + 2 \omega
$$ 
for $\omega \in \Omega^1(M,\E)$.
For $\omega\in\Omega^1_{\cp}(M,\E)$ integration by parts yields
\begin{align*}
\int_M (\Delta\omega,\omega)
&=\int_M (\nabla^t\nabla\omega,\omega)+  2\int_M \vert\omega \vert^2\\
&=\int_M \vert \nabla\omega\vert^2 +  2 \int_M  \vert\omega \vert^2\\
&\geq  2 \int_M  \vert\omega \vert^2
\end{align*}
This means we have a positive lower bound for $\Delta$ on $\Omega^1_{\cp}(M,\E)$:
$$
\left\langle \Delta\omega, \omega\right\rangle_{L^2}\geq C \left\langle \omega, \omega \right\rangle_{L^2}
$$
with $C=2$.
Since $(\E,\nabla^{\E},h^{\E})$ is cone-admissible according to Proposition \ref{Econeadmissible_spherical}, we obtain
$$\Delta_{\Fr}=\Delta(d_{\max})$$ via Corollary \ref{Friedrichs=dmax}. Since the Friedrichs extension
preserves lower bounds, we conclude
$$
\mathcal{H}^1_{\max}=\ker \Delta^1(d_{\max})=0\,.
$$ 
Finally Theorem \ref{HodgeConeManifold} identifies $L^2$-cohomology with the $d_{\max}$-harmonic forms. This implies $H^1_{L^2}(M,\E)=0$ and proves the theorem.
\end{pf}
\subsubsection{The Euclidean case}
\begin{thm} Let $C$ be a Euclidean cone-3-manifold with cone-angles $\leq
\pi$. Let $\E_{trans}\subset\E$ be the parallel subbundle of infinitesimal translations of
$M=C\setminus\Sigma$. Then
$$
H^1_{L^2}(M,\E_{trans})\cong\{\omega\in \Omega^1(M,\E_{trans})\,\vert\, \nabla\omega=0\}\,.
$$ 
\end{thm}
\begin{pf}
The Weitzenb\"ock formula for the Hodge-Laplace operator on $\E_{trans}$-valued $1$-forms in the Euclidean case (i.e.~$\kappa=0$) amounts to
$$
\Delta \omega = \nabla^t\nabla \omega
$$ 
for $\omega \in \Omega^1(M,\E_{trans})$.
This implies with Corollary \ref{CorFriedrichs}
that
$$
\Delta^1_{\Fr}=\nabla^t_{\max}\nabla_{\min}\,.
$$
Since $\E_{trans}\subset\E$ is cone-admissible according to Lemma \ref{Etransconeadmissible}, 
we obtain
$$
\Delta^1_{\Fr} = \Delta^1(d_{max})\,.
$$
via Corollary \ref{Friedrichs=dmax}. This implies that
$$
\Delta^1(d_{max})=\nabla^t_{\max}\nabla_{\min}\,.
$$
For $\omega\in \ker \Delta^1(d_{\max})$ we have
\begin{align*}
0=\langle \Delta(d_{\max})\omega,\omega\rangle_{L^2}=\langle\nabla^t_{\max}\nabla_{\min}\omega,\omega\rangle_{L^2}=\|\nabla_{\min}\omega\|^2_{L^2}\,.
\end{align*}
We conclude that $\omega \in \ker \nabla_{\min}$. On the other hand, if
$\omega \in \ker \nabla_{\max}$, then 
clearly $\omega \in ker D_{\max}$. Since $D$ is essentially selfadjoint according to Corollary \ref{essentialselfadjointness},
 $\ker D_{\max}=\ker D_{\min}=\mathcal{H}_{\max}$.
We obtain
$$
\ker \nabla_{\max} \subset \mathcal{H}^1_{\max} \subset \ker \nabla_{\min}\,,
$$
which proves the theorem 
via Theorem \ref{HodgeConeManifold},
since $\mathcal{H}_{\max}$ consists of smooth forms. Note also that a parallel form $\omega$ will automatically be $L^2$-bounded, since $\nabla$ is compatible with the metric on $\E_{trans}$.
\end{pf}
\subsubsection{The hyperbolic case}
\begin{thm} Let $C$ be a hyperbolic cone-3-manifold with cone-angles $\leq
\pi$.  Let $M=C\setminus\Sigma$ and $(\E,\nabla^{\E})$ be the vector-bundle of
infinitesimal isometries of $M$ with its natural flat connection. Then
$$
H^1_{L^2}(M,\E)=0\,.
$$
\end{thm}
\begin{pf}
The proof follows the same scheme as in the spherical case. For convenience of the reader we also give full details in this case.

We recall that in the hyperbolic case we have a Weitzenb\"ock formula for the Hodge-Laplace operator
for $\E$-valued $1$-forms of the type
$$
\Delta \omega = \mathcal{D}^t\mathcal{D}\omega+\mathcal{D}\mathcal{D}^t\omega+H\omega\,,
$$
where
$$
\left \langle H\omega,\omega\right\rangle_{L^2} \geq C \left \langle \omega,\omega\right\rangle_{L^2}
$$
for $C>0$ independent of $\omega\in\Omega^1(M,\E)$.
For $\omega\in\Omega^1_{\cp}(M,\E)$ integration by parts yields
\begin{align*}
\int_M (\Delta\omega,\omega)
&=\int_M (\mathcal{D}^t\mathcal{D}\omega,\omega)+\int_M (\mathcal{D}\mathcal{D}^t\omega,\omega) + \int_M (H\omega,\omega)\\
&=\int_M \vert \mathcal{D}\omega\vert^2 +\int_M \vert \mathcal{D}^t\omega\vert^2+ \int_M (H\omega,\omega)\\
&\geq C \int_M  \vert\omega \vert^2
\end{align*}
This means we have a positive lower bound for $\Delta$ on $\Omega^1_{\cp}(M,\E)$:
$$
\langle \Delta\omega, \omega\rangle_{L^2}\geq C \left\langle \omega, \omega \right\rangle_{L^2}
$$
for $C>0$. Since $(\E,\nabla^{\E},h^{\E})$ is cone-admissible according to Proposition \ref{Econeadmissible_hyperbolic}, we obtain
$$\Delta_{\Fr}=\Delta(d_{\max})$$ via Corollary \ref{Friedrichs=dmax}. Since the Friedrichs extension
preserves lower bounds, we conclude
$$
\mathcal{H}^1_{\max}=\ker \Delta^1(d_{\max})=0\,.
$$ 
Finally Theorem \ref{HodgeConeManifold} identifies $L^2$-cohomology with the $d_{\max}$-harmonic forms. This implies $H^1_{L^2}(M,\E)=0$ and proves the theorem. 
\end{pf}

\section{Deformation theory}

In this chapter we study the deformation space of cone-manifold structures on
a $3$-dimensional cone-manifold of given topological type $(C,\Sigma)$. It is
convenient to use the more general framework of $(X,G)$-structures and
deformations thereof, in particular since there is a quite general theorem of
\cite{Gol}, which relates the local structure of the deformation space of
$(X,G)$-structures to the local structure of $X(\pi_1M,G)$. By $X(\pi_1M,G)$
we denote the quotient of $R(\pi_1M,G)$, the space of representations of $\pi_1M$ in $G$, by the conjugation action of $G$.

The $(X,G)$-structures relevant for our situation will be $X=\M$ and $G=\Isom\M$, in fact by a theorem of \cite{Cul}, the holonomy representation of a $3$-dimensional cone-manifold structure may always be lifted to the universal covering group of $\Isom\M$, which in the hyperbolic case is $\SL_2(\C)$ and in the spherical case $\SU(2)\times\SU(2)$.

We will use the $L^2$-vanishing theorem to analyze local properties of $\SL_2(\C)$- and $\SU(2)$-representation spaces. From this we will be able to conclude local rigidity in the hyperbolic and in the spherical case.

\subsection{$(X,G)$-structures}\label{XGstructures}
Let $(X,g^X)$ be a Riemannian manifold upon which a Lie group $G$ acts
transitively by isometries. Let $M$ be manifold of the same dimension as
$X$. Then we say that $M$ carries an $(X,G)$-structure if $M$ is locally
modelled on $X$, i.e.~there is a covering of $M$ by charts
$\{\varphi_i:U_i\rightarrow X\}_{i\in I}$ such that for each connected
component of $C$ of $U_i\cap U_j$ there exists  $g_{C,i,j}\in G$ such that
$g_{C,i,j}\circ \varphi_i = \varphi_j$ on $C$. The collection of charts
$\{\varphi_i:U_i\rightarrow X\}_{i\in I}$ is called an $(X,G)$-atlas and an
$(X,G)$-structure on $M$ is a maximal $(X,G)$-atlas. A detailed discussion of
this kind of structure may be found in $\cite{Gol}$, which we will use as the main reference for this section.  

Let us fix basepoints $x_0\in M$ and $p_0\in \pi^{-1}(x_0)$, where $\pi: \widetilde{M}\rightarrow M$ is the universal covering of $M$. Then an $(X,G)$-structure on $M$ 
together with the germ of an $(X,G)$-chart $\varphi:U \rightarrow X$ around $x_0$
determines by analytic continuation of $\varphi$ a local diffeomorphism
$$
\dev: \widetilde{M} \longrightarrow X\,,
$$  
the developing map, and a representation
$$
\hol: \pi_1(M,x_0) \longrightarrow G\,,
$$
the holonomy representation, such that $\dev$ is equivariant with respect to $\hol$, i.e.
$$
\dev \circ \gamma = \hol(\gamma) \circ \dev
$$
for all $\gamma \in \pi_1(M,x_0)$.
Conversely, a local diffeomorphism $\dev:\widetilde{M} \rightarrow X$ equivariant with respect to some representation $\hol: \pi_1(M,x_0) \rightarrow G$ as above, defines an $(X,G)$-structure on $M$ together with the germ of an $(X,G)$-chart at $x_0$. Note that $\hol$ is uniquely determined by $\dev$ and the equivariance condition.

Let $\mathcal{D}'_{(X,G)}(M)$ be the space of developing maps with the topology of $C^{\infty}$-convergence on compact sets. As usual we equip $R(\pi_1(M,x_0),G)$, the set of representations of $\pi_1(M,x_0)$ in $G$, with the compact-open topology. Associating its holonomy representation with a developing map yields a continuous map
\begin{align*}
\mathcal{D}'_{(X,G)}(M) &\longrightarrow R(\pi_1(M,x_0),G)\\
\dev & \longmapsto \hol\,.
\end{align*}
For simplicity we assume that $M$ is diffeomorphic to the interior of a compact manifold with boundary $M\cup\partial M$, which is certainly the case for the object of our main concern, namely the smooth part of a $3$-dimensional cone-manifold. 

Following $\cite{CHK}$ we introduce the equivalence relation $\sim$ on 
the space of developing maps, which is generated by {\it isotopy} and {\it thickening}. Clearly $\operatorname{Diff}_0(M)$, the group of diffeomorphisms of $M$ isotopic to the identity, acts on the space of developing maps, two structures equivalent under this action will be called {\it isotopic}. On the other hand, if an $(X,G)$ structure on $M$ extends to $M\cup \partial M \times [0,\varepsilon)$ for some $\varepsilon >0$, this gives rise to an $(X,G)$-structure on $M$, which we will call a {\it thickening} of the original structure. 
Let $$ \mathcal{D}_{(X,G)}(M)= \mathcal{D}'_{(X,G)}(M)/\!\sim\,.$$ 
We obtain a $G$-equivariant map
\begin{align*}
\mathcal{D}_{(X,G)}(M) & \longrightarrow R(\pi_1(M,x_0),G)\\
[\dev] & \longmapsto \hol\,.
\end{align*}
We define the {\em deformation space} of $(X,G)$-structures to be the quotient
$$
\mathcal{T}_{(X,G)}(M):=\mathcal{D}_{(X,G)}(M)/G\,.
$$
Let $X(\pi_1(M,x_0),G)$ denote the $G$-quotient of $R(\pi_1(M,x_0),G)$ by conjugation, properties of this quotient in our particular context will be discussed in greater detail in subsequent sections.
 
Let us assume that the action of $G$ on $R(\pi_1(M,x_0),G)$ by conjugation is proper, this implies in particular by the $G$-equivariance of the above map, that the action of $G$ on $\mathcal{D}_{(X,G)}(M)$ is also proper. In this situation the arguments of $\cite{Gol}$ (cf.~also the discussion in {\cite{CHK}) yield the following theorem about the local structure of the deformation space of $(X,G)$-structures:
\begin{thm}[deformation theorem]{\rm\cite{Gol}} If the action of $G$ on $R(\pi_1(M,x_0),G)$ by conjugation is proper, then the map
\begin{align*}
\mathcal{T}_{(X,G)}(M) & \longrightarrow X(\pi_1(M,x_0),G)\\
[\dev] & \longmapsto [\hol]
\end{align*}
is a local homeomorphism.
\end{thm}
This theorem explains the meaning of representation varieties in the study of
deformations of $(X,G)$-structures: Local properties of the deformation space
of $(X,G)$-structures on $M$ translate into local properties of
$X(\pi_1(M,x_0),G)$ and vice versa.\\
\\
By a theorem of M. Culler (cf.~\cite{Cul}) the holonomy representation of a cone-3-manifold may be lifted to 
the universal covering group of $\Isom \M$: 
$$
\widetilde{\hol}: \pi_1 M \longrightarrow \widetilde{\Isom} \M
$$
In the hyperbolic case $\widetilde{\Isom} \Hyp^3=\SL_2(\C)$. We obtain that the flat vector-bundle of infinitesimal isometries may be written as
$$
\E= \widetilde{M}\times_{Ad\circ\widetilde{\hol}}{\sl}_2(\C)\,.
$$
As a consequence $\E$ has a parallel complex structure, such that in particular all the cohomology spaces
$H^i(M,\E)$ are complex vector spaces.

In the spherical case $\widetilde{\Isom}\Sph^3=\SU(2)\times\SU(2)$. Therefore the lift of the holonomy splits
as a product representation
$$
\widetilde{\hol}=(\hol_1,\hol_2) : \pi_1M \longrightarrow \SU(2)\times\SU(2)\,,
$$
in particular the flat vector-bundle of infinitesimal isometries splits as a direct sum of parallel subbundles:
$$
\E=\E_1 \oplus \E_2\,,
$$
where
$$
\E_i=\widetilde{M}\times_{Ad\circ\hol_i}\su(2)\,.
$$
Consequently $H^i(M,\E)=H^i(M,\E_1)\oplus H^i(M,\E_2)$ for all $i$.\\
\\
For notational convenience we will drop the distinction between $\hol$ and
$\widetilde{\hol}$ from here.

\subsection{The representation variety}
Let $\Gamma$ be a finitely generated discrete group. Once and for all we fix a presentation
$\langle \gamma_1,\ldots,\gamma_n | (r_i)_{i\in I}\rangle
$ of $\Gamma$. The cardinality of the indexset $I$ may a priori be infinite,
however most of the groups we deal with will turn out to be finitely
presented. Let $G=\SL_2(\C)$ or $\SU(2)$. The {\em representation variety} $R(\Gamma,G)$ is defined to be the
set of group homomorphisms $\rho:\Gamma\rightarrow G$. $R(\Gamma,G)$
endowed with the compact-open topology is a Hausdorff space, compact in the
case of
$\SU(2)$.

The relations $r_i$ define functions $f_i:G^n \rightarrow G$ such that
$R(\Gamma,G)$ may be identified with the set $\{(A_1,\ldots,A_n)\in G^n |\,
f_i(A_1,\ldots,A_n)=1 \}$. Since $\SL_2(\C)$ is a $\C$-algebraic
(resp. $\SU(2)$ a $\R$-algebraic) group and the $f_i$ are polynomial maps,
$R(\Gamma,G)$ acquires the structure of a $\C$-algebraic (resp. $\R$-algebraic)
set. Note that $R(\Gamma,G)$ won't be a smooth space in general.

The action of $G$ on $G^n$ by simultaneous conjugation leaves the set
$R(\Gamma,G)\subset G^n$ invariant. Therefore the quotient 
$X(\Gamma,G)=R(\Gamma,G)/G$ is well defined. We endow $X(\Gamma,G)$ with the
quotient topology. $X(\Gamma,G)$ will in general be neither smooth nor even
Hausdorff. $X(\Gamma,G)$ as we have defined it should not be confused with a quotient constructed in the algebraic category. This usually requires
arguments from geometric invariant theory, which we can avoid using here.

A smooth family of representations $\rho_t:\Gamma\rightarrow G$ with
$\rho_0=\rho$ defines a group
$1$-cocycle $z:\Gamma \rightarrow \g$, where
$$
z(\gamma)=\left.\textstyle\frac{d}{dt}\right\vert_{t=0}\rho_t(\gamma)\rho(\gamma)^{-1}
$$
for $\gamma \in \Gamma$. Recall that $Z^1(\Gamma,\g)$, the space of $1$-cocycles of $\Gamma$ with coefficients
in the representation $Ad\circ\rho:\Gamma\rightarrow\GL(\g)$, is the the space of maps $z : \Gamma \rightarrow \g$
such that
$$
 z(ab)=z(a)+ \left(Ad\circ\rho(a)\right)z(b)
$$
for all $a,b \in \Gamma$. A cocycle $z$ is a coboundary if there exists some $v\in\g$ such that
$$
z(a)=v-\left(Ad\circ\rho(a)\right)v
$$
for all $a\in\Gamma$. Let $B^1(\Gamma,\g)$ be the space of $1$-coboundaries. Now by definition
$$
H^1(\Gamma,\g) =  Z^1(\Gamma,\g)/B^1(\Gamma,\g)
$$
is the first {\em group cohomology} group of $\Gamma$ with coefficients in the representation $Ad\circ\rho: \Gamma\rightarrow \GL(\g)$. $H^1(\Gamma,\g)$ is a real vector space. Recall further that
$$
H^0(\Gamma,\g) = Z^0(\Gamma,\g) = \{ v \in \g \vert (Ad \circ \rho(\gamma))v = v \, \forall \gamma \in \Gamma \}\,.
$$
For more details on group cohomology see \cite{Bro}.

We refer to $Z^1(\Gamma,\g)$ as the space of infinitesimal deformations of the
representation $\rho$. We call a $1$-cocycle $z$ integrable, if there exists a
(local) deformation $\rho_t$, which is tangent to $z$ in the above sense.
 
It is easy to see that $z \in B^1(\Gamma,\g)$ if and only if $z$ is tangent to
the orbit of $G$ through $\rho$, i.e.~there exists a smooth curve $g_t$ in $G$ with
$g_0=1$ such that
$$
z(\gamma)=\left.\textstyle\frac{d}{dt}\right\vert_{t=0}g_t\rho(\gamma)g_t^{-1}\rho(\gamma)^{-1}
$$ 
for $\gamma \in \Gamma$. A deformation $\rho_t(\gamma)=g_t\rho(\gamma)g_t^{-1}$ will be considered trivial.

We use the following observation due to A. Weil (cf.~\cite{Wei}): 
A map $z: \Gamma \rightarrow
\g$ defines a group 1-cocycle if and only if the map 
\begin{align*}
(Ad\circ\rho,z): \Gamma &\longrightarrow \GL(\g)\ltimes\g\\
\gamma &\longmapsto (Ad\circ\rho(\gamma), z(\gamma))
\end{align*}
is a group homomorphism. $\GL(\g)\ltimes\g$ is the affine group of the
vector-space $\g$. Using the fixed presentation of $\Gamma$, this identifies $Z^1(\Gamma,\g)$
with a linear subspace of $\g^n$. More precisely,
the relations $r_i$ determine linear functions $g_i:\g^n \rightarrow
\g$, such that 
$$
Z^1(\Gamma,\g)=\left\{ (a_1,\ldots,a_n) \in \g^n \vert
g_i(a_1,\ldots,a_n)=0  \, \forall i\in I \right\}.
$$
On the other hand, $\ker df_i$ may be identified with a subspace of
$\g^n$ via 
$$
( \dot{A}_1,\ldots,\dot{A}_n) \mapsto (\dot{A}_1A_1^{-1},
\ldots, \dot{A}_nA_n^{-1})\,. 
$$
With these identifications we have the following lemma:
\begin{lemma}
$Z^1(\Gamma,\g)=\cap_{i\in I}\ker df_i$.
\end{lemma}
\begin{pf}
A straightforward calculation shows that
$df_i(\dot{A}_1,\ldots,\dot{A}_n)=0$ for $\dot{A}_i \in T_{A_i}G$ if and
only if $g_i(a_1,\ldots,a_n)=0$, where
$a_i=\dot{A}_iA_i^{-1}$. 
\end{pf}\\
\\
If the equations $(f_i)_{i\in I}$ cut out $R(\Gamma,G)$ transversely near $\rho$, then the previous lemma identifies $Z^1(\Gamma,\g)$ with the tangent space of $R(\Gamma,G)$ at the point $\rho$. In particular $\rho$ will be a smooth point.
If furthermore the $G$-action on $R(\Gamma,G)$ by conjugation is free and proper, then $X(\Gamma,G)$ will be smooth near $\chi=[\rho]$ and the tangent space at $\chi$ may be identified with $H^1(\Gamma,\g)$.

\subsection{Integration and group cohomology}

We wish to represent group cocycles of $\pi_1M$
with coefficients in the representation $Ad\circ\hol : \pi_1M \rightarrow \g=\isom \M$
by differential forms
on $M$ with values in $\E$. This will be achieved by means of integration.

Let $x_0$ be a base point in $M$, then for $\gamma \in \pi_1(M,x_0)$ and a
closed $1$-form $\omega \in
\Omega^1(M,\E)$ we define
$$
\int_\gamma \omega = \int_0^1
\tau_{\gamma(t)}^{-1}\omega(\dot{\gamma}(t)) dt \in \E_{x_0},
$$
where $\tau_{\gamma(t)}$ denotes the parallel transport along $\gamma$ from
$x_0=\gamma(0)$ to $\gamma(t)$. Since $\omega$ is closed, the integral
depends only on the homotopy class of $\gamma$. If we identify $\E_{x_0}$ with $\g$, then
we may set 
$$
z_{\omega}(\gamma) = \int_\gamma \omega \in \g\,.
$$ 
Alternatively, we may proceed as follows:
The flat bundle $\E$ may be described as an associated
bundle $\E = \widetilde{M} \times_{Ad\circ\hol}\g$. 
$1$-forms $\omega \in \Omega^1(M,\E)$ correspond to $1$-forms
$\widetilde{\omega} \in \Omega^1(\widetilde{M},\g)$ satisfying the following
equivariance condition:
$$\gamma^*\widetilde{\omega}=\left(Ad\circ\hol(\gamma)\right)\widetilde{\omega}$$
for all $\gamma \in \pi_1(M,x_0)$.  
For $\omega \in \Omega^1(M,\E)$ closed consider $\widetilde{\omega} \in
\Omega^1(\widetilde{M},\g)$, which will again be closed. Let $p_0\in\pi^{-1}(x_0)$ be a base point in
$\widetilde{M}$. Now since $\widetilde{M}$ is simply connected, there exists a primitive $F\in C^{\infty}(M,\g)$ such that $dF=\widetilde{\omega}$.
For $\gamma \in \pi_1(M,x_0)$ we define
$$
z_{\omega}(\gamma) = \int_{\gamma}\omega = F(\gamma p_0)-F(p_0)\in\g\,.
$$
Since $F$ is determined up to an additive constant, this is well defined. 

Both definitions of
the map $z_{\omega}:\pi_1M \rightarrow \g$ associated with the closed form
$\omega\in\Omega^1(M,\E)$ clearly agree. The proof of the following lemma is
straightforward and left to the reader:
\begin{lemma} If $\omega \in \Omega^1(M,\E)$ is closed, then $z_{\omega}$
defines a group cocycle, i.e.~$z_{\omega} \in
Z^1(\pi_1M,\g$). $\omega$ is exact if and only if $z_{\omega} \in B^1(\pi_1M,\g)$. 
\end{lemma}
As a consequence of the preceding lemma, we obtain that the period map
\begin{align*}
P: H^1(M,\E) & \longrightarrow H^1(\pi_1M,\g)\\
[\omega] & \longmapsto [\gamma\mapsto\textstyle\int_{\gamma}\omega]
\end{align*}
is well defined and injective. Since we know from general considerations
(cf.~\cite[Thm.~5.2]{Bro}) that
$H^i(M,\E) \cong H^i(\pi_1M,\g)$ for $i \in \{0,1\}$, we find that the period
map provides an explicit isomorphism between $H^1(M,\E)$ and $H^1(\pi_1M,\g)$. 


\subsection{Isometries}

\subsubsection{Isometries of $\Hyp^3$}
The action of $\SL_2(\C)$ on $\Hyp^3$ by Poincar\'{e} extension identifies
$\SL_2(\C)$ with the universal cover of $\Isom \Hyp^3=\PSL_2(\C)$. Here we use the upper half space model. 
Let $\phi:\SL_2(\C)\rightarrow \Isom \Hyp^3$ denote the covering projection.

Semisimple elements in $\SL_2(\C)$ project to semisimple isometries. A semisimple
isometry $\phi$ has an invariant axis; this is the unique geodesic, where
$\delta_{\phi}$, the
displacement function of $\phi$, assumes its minimum. If this minimum is positive, we
call $\phi$ hyperbolic, otherwise elliptic. Parabolic elements in $\SL_2(\C)$
project to parabolic isometries. Parabolic isometries have a unique fixed point
at infinity. The following is well-known:  
\begin{lemma}
$A,B\in \SL_2(\C)$ commute if and only if $\phi(A),\phi(B)$ are either semisimple isometries and preserve the same axis $\gamma$ or $\phi(A),\phi(B)$ are parabolic isometries with the same fixed point at infinity.
\end{lemma}
The stabilizer of an oriented geodesic $\gamma$ is isomorphic to $\C^*$, more precisely, if we work in the upper half space model $\Hyp^3=\C \times \R_+$, then for $\gamma=\{0\}\times\R_+$ we obtain
$$
\Stab_{\SL_2(\C)}(\gamma) = \left\{\left( \begin{array}{cc}
\lambda & 0\\
0 & \lambda^{-1} \end{array}\right) : \lambda \in \C^*\right\}\,.
$$  
$S^1 \subset \C^*$ corresponds to pure rotations around $\gamma$, while $\R
\subset \C^*$ corresponds to pure translations along $\gamma$.
Recall that for a Killing vectorfield $X$ on $\Hyp^3$ we denote by $ \sigma_X =(\nabla X,
X)\in {\sl_2(\C)}$ the corresponding parallel section.
In particular, if we choose cylindrical coordinates $(r,\theta,z)$ around
$\gamma$, we see that 
\begin{align*}
\sigma_{\partial/\partial \theta} &= \frac{1}{2}\left( \begin{array}{cc}
i & 0\\
0 & -i \end{array}\right)\in \sl_2(\C)
\end{align*}
and
\begin{align*}
\sigma_{\partial/\partial z} &=  \frac{1}{2}\left( \begin{array}{cc}
1 & 0\\
0 & -1 \end{array}\right)\in \sl_2(\C)\,.
\end{align*}
Note in particular that $\sigma_{\partial/\partial
  \theta}=i\sigma_{\partial/\partial z}$. The factor $1/2$ comes from
the fact that $\SL_2(\C)$ is a twofold cover of $\Isom \Hyp^3$.\\
\\
Let $A \in \SL_2(\C)$ be semisimple and $\phi=\phi(A) \in \Isom \Hyp^3$. Then $A$ is conjugate to
$\diag(\lambda,\lambda^{-1})$ in $\SL_2(\C)$ for $\lambda \in \C^*$. Now let $z \in
\C/2\pi i\Z$ such that
$
\lambda = \exp(z)
$. We define $\mathcal{L}(A)=2z \in
\C/2\pi i\Z$.
Then $\mathcal{L}(A)$ is determined by $A$ up to
sign. $\mathcal{L}(A)$ is called the {\it complex length} of $A$. 

For $ A \neq
\pm \id$ we can orient the axis $\gamma$ of $\phi$ and remove 
the sign ambiguity
of $\mathcal{L}$ consistently in a
neighbourhood of $A$ in $\SL_2(\C)$. The real part of $\mathcal{L}(A)$ equals the (signed) translation length
of $\phi$ along $\gamma$, while the imaginary part equals
the angle of rotation around $\gamma$.
We obtain
$$
\tr A = 2 \cosh (z) = \pm 2\cosh(\mathcal{L}(A)/2)
$$
and by the inverse function theorem:
\begin{lemma}\label{complexlengthhyp}
Let $A\ \neq \pm \id \in \SL_2(\C)$ be semisimple. There exist
neighbourhoods $U$ of $A$ in $\SL_2(\C)$ and $V$ of $\tr A$ in $\C$ and a biholomorphic map $\phi: V \rightarrow \mathcal{L}(V)\subset\C$ such
that $\tr(U) \subset V$ and $\phi \circ \tr = \mathcal{L}$ on $U$. 
\end{lemma}

\subsubsection{Isometries of $\Sph^3$}
We identify $\Sph^3$ with the unit quaternions, i.e.~$\Sph^3=\{x\in\H:\vert x
\vert=1\}$. If we view the quaternions as a subalgebra of $\C^{2\times 2}$ via$$
1 \mapsto  \left(  \begin{array}{cc}
1 & 0\\
0 & 1 \end{array}\right)\,,\,
i \mapsto  \left(  \begin{array}{cc}
i & 0\\
0 & -i \end{array}\right)\,,\,
j \mapsto  \left(  \begin{array}{cc}
0 & 1\\
-1 & 0 \end{array}\right)\,,\,
k \mapsto  \left(  \begin{array}{cc}
0 & i\\
i & 0 \end{array}\right)
\,,$$ 
$\Sph^3$ gets identified with the group $\SU(2)$ via
\begin{align*}
\Sph^3 & \longrightarrow \SU(2)\\
a+bj &\longmapsto \left(  \begin{array}{cc}
a & b\\
-\overline{b} & \overline{a} \end{array}\right),
\end{align*}
where $a,b\in\C$ with $\vert a\vert^2+\vert b\vert^2=1$. The map
\begin{align*}
\phi: \SU(2)\times \SU(2) &\longrightarrow \SO(4)\\
(A,B) &\longmapsto (x\mapsto AxB^{-1})
\end{align*}
exhibits $\SU(2)\times \SU(2)$ as the universal cover of $\Isom
\Sph^3=\SO(4)$. 
Note that the diagonal matrices 
$$
\left\{
\left( \begin{array}{cc}
\lambda & 0\\
0 & \overline{\lambda} \end{array}\right) : \lambda \in S^1\right\}\subset \SU(2)
$$
correspond to the geodesic $\gamma=\C \cap \Sph^3$, where as usual $\C$ is identified with
$\R\oplus \R i \subset\H$. For any geodesic $\gamma \subset \Sph^3$ let us
denote by $\gamma^{\perp}$ the geodesic which lies in the plane orthogonal to
$\gamma$. In the above case $\gamma^{\perp}=\C j\cap \Sph^3=(\R j\oplus \R
k)\cap \Sph^3$, which corresponds to the set of matrices 
$$
\left\{\left( \begin{array}{cc}
0 & \lambda\\
-\overline{\lambda} & 0 \end{array}\right) : \lambda \in S^1 \right\} \subset \SU(2).
$$
A spherical isometry may be put in a standard form, namely if an isometry is
represented as $\phi=\phi(A,B)$
with $A,B \in 
\SU(2)$, then by conjugation we may achieve that
$A=\diag(\lambda,\overline{\lambda})$ and $B=\diag(\mu,\overline{\mu})$ with
$\lambda,\mu \in S^1$. The matrix $A$ corresponds to $\lambda \in \C\cap\Sph^3$ and $B$ to
$\mu \in \C\cap\Sph^3$ if we identify $\SU(2)$ with $\Sph^3$ as above.
Then for $x \in \Sph^3$ we have
$
\phi(x)=\lambda x\overline{\mu}
$, 
such that $\phi$ preserves the Hopf-fibrations,
which are associated with the complex structures $x \mapsto ix$ and $x \mapsto xi$
on $\H$.

In particular, $\phi$ preserves $\gamma=\C\cap\Sph^3$ and
$\gamma^{\perp}=\C j\cap\Sph^3$, more precisely we have
$
\phi(\eta)=\lambda\overline{\mu}\eta
$
for $\eta\in S^1=\C\cap\Sph^3$, and 
$
\phi(\eta j)=(\lambda\mu\eta) j
$
for $\eta j\in \C j\cap\Sph^3$. Note that $\gamma$ and $\gamma^{\perp}$ are
the common fibers of the two fibrations, which are transverse everywhere else.

If $\mu=1$, then $\phi$
translates along the fibers of the Hopf-fibration obtained by
left-multiplication with $S^1$, in particular the
displacement of $\phi$ is constant on $\Sph^3$. Similarly, if $\lambda=1$, then $\phi$
translates along the fibres of the Hopf-fibration obtained by right-multiplication with
$S^1$. Again the 
displacement of $\phi$ will be constant on $\Sph^3$. 

If $\lambda=\mu$, then $\phi$ is a pure rotation around $\gamma$,
or equivalently, a pure translation along $\gamma^{\perp}$. Similarly, if
$\lambda=\overline{\mu}$,  then $\phi$ is a pure rotation around $\gamma^{\perp}$,
or equivalently, a pure translation along $\gamma$.

Recall that for a Killing vectorfield $X$ on $\Sph^3$ we denote by $ \sigma_X =(\nabla X,
X)\in \su(2) \oplus \su(2)$ the corresponding parallel section.
In particular, if we choose cylindrical coordinates $(r,\theta,z)$ around
$\gamma$, we see that 
\begin{align*}
\sigma_{\partial/\partial \theta} 
&=\left(  \frac{1}{2}\left( \begin{array}{cc}
i & 0\\
0 & -i \end{array}\right),  \frac{1}{2}\left( \begin{array}{cc}
i & 0\\
0 & -i \end{array}\right) \right)
 \in \su(2)\oplus \su(2)
\end{align*}
and
\begin{align*}
\sigma_{\partial/\partial z} 
&= \left( \frac{1}{2}\left( \begin{array}{cc}
i & 0\\
0 & -i \end{array}\right), \frac{1}{2}\left( \begin{array}{cc}
-i & 0\\
0 & i \end{array}\right) \right)
\in \su(2)\oplus \su(2)\,.
\end{align*}
The factors $1/2$ arise from the fact that $\SU(2)\times\SU(2)$ is a twofold cover of $\Isom \Sph^3$. The following is immediate from the above discussion:
\begin{lemma}
$\phi_1,\phi_2 \in \Isom \Sph^3$ commute if and only they preserve the same pair of orthogonal axes $\{\gamma,\gamma^{\perp}\}$.
\end{lemma}
Since $\phi(A,B)1=1$ if and only if $A=B\in\SU(2)$ we obtain:
\begin{lemma}
$\phi=\phi(A,B)\in \Isom \Sph^3$ has a fixed point if and only if $A$ is conjugate to $B$ within $\SU(2)$.
\end{lemma}
We want to define an analogue of the complex length in the spherical case. If $\phi=\phi(A,B)$ with $A$ conjugate to $\diag(\lambda,\overline{\lambda})$ and $B$ conjugate to $\diag(\mu,\overline{\mu})$, then let $x\in\R/2\pi\Z$ such that
$
\lambda= \exp(ix)
$
and $y\in\R/2\pi\Z$ such that
$
\mu= \exp(iy)
$.
We define $\mathcal{L}_1(A,B) = x-y$ and $\mathcal{L}_2(A,B) = x+y$.
Then $\mathcal{L}(A,B)= (\mathcal{L}_1(A,B),\mathcal{L}_2(A,B)) \in \R^2/2\pi\Z^2$ is
determined by $A$ and $B$ up to an overall sign and up to switching
components. 

Let in the following $A \neq \pm \id$ and $B \neq \pm \id$. If $\phi$ preserves a pair of orthogonal axes $\{\gamma, \gamma^{\perp}\}$, these ambiguities can be removed in a neighbourhood of $(A,B)$ by orienting $\gamma$. Let us again call $\mathcal{L}(A,B)$ the "complex" length of $(A,B)\in \SU(2)\times \SU(2)$. $\mathcal{L}_1(A,B)$ equals the (signed) translation length along $\gamma$, while $\mathcal{L}_2(A,B)$ equals the (signed) translation length along $\gamma^{\perp}$. We obtain
$$
\tr A = 2 \cos x = \pm 2 \cos\bigl((\mathcal{L}_1(A,B)+\mathcal{L}_2(A,B))/2\bigr)
$$
and
$$
\tr B = 2 \cos y = \pm 2 \cos\bigl((-\mathcal{L}_1(A,B)+\mathcal{L}_2(A,B))/2\bigr)\,.
$$
We set $\Tr_1(A,B)=\tr A$, $\Tr_2(A,B)=\tr B$ and $\Tr=(\Tr_1,\Tr_2)$.
By the inverse function theorem we obtain:
\begin{lemma}\label{complexlengthsph}
Let $(A,B) \in \SU(2) \times \SU(2)$ with $A \neq \pm \id$ and $B \neq \pm \id$. There exist
neighbourhoods $U$ of $(A,B)$ in $\SU(2) \times \SU(2)$  and $V$ of $\Tr
(A,B)$ in $\R^2$ and a diffeomorphism $\phi: V \rightarrow \mathcal{L}(V)\subset\R^2$ such
that $\Tr(U) \subset V$ and $\phi \circ \Tr = \mathcal{L}$ on $U$. 
\end{lemma}

\subsection{Cohomology computations}
Let $C$ be a $3$-dimensional cone-manifold with cone-angles $\leq \pi$. Under
this cone-angle bound, a connected component
of the singular locus $\Sigma$ will either be a circle or a (connected) trivalent
graph, cf.~Chapter \ref{conemanifolds}, see also \cite {CHK} and \cite{BLP2}. Let $M_{\varepsilon}=M\setminus B_{\varepsilon}(\Sigma)$, where
$B_{\varepsilon}(\Sigma)$ is the open $\varepsilon$-tube around $\Sigma$. Let $U_{\varepsilon}(\Sigma)=B_{\varepsilon}(\Sigma)\setminus \Sigma$. Then $M_{\varepsilon}$
is topologically a manifold with boundary, which is a deformation retract of $M$. $\partial
M_{\varepsilon}$ consists of tori and
surfaces of higher genus. $\partial M_{\varepsilon}=\partial
U_{\varepsilon}(\Sigma)$ is a deformation retract of $U_{\varepsilon}(\Sigma)$.

Without loss of generality we may assume in the following that $\Sigma$ is connected.
\subsubsection{The torus case}
Let $\Sigma = S^1$. Then $U_{\varepsilon}(\Sigma)$ is given as $(0,\varepsilon)\times T^2$, where $T^2=\R^2/\Lambda$ and $\Lambda$ is the lattice generated by $(\theta,z)\mapsto(\theta+\alpha,z)$ and $(\theta,z)\mapsto(\theta-t,z+l)$. The metric is given as $g=dr^2+\sn^2(r)d\theta^2+\cs^2(r)dz^2$. Here $\alpha,t$ and $l$ are the parameters, which determine the geometry of $U_{\varepsilon}(\Sigma)$, namely the cone-angle, the twist and the length of the singular tube. Note that a function $f$ in the coordinates $(r,\theta,z)$ descends to a function on $U_{\varepsilon}(\Sigma)$ if and only if $f(r,\theta,z)=f(r,\theta+\alpha,z)$ and $f(r,\theta,z+l)=f(r,\theta+t,z)$. Note also that $H^i(U_{\varepsilon}(\Sigma),\cdot)=H^i(T^2,\cdot)$ for any local coefficient system.

The forms $d\theta$ and $dz$ are invariant under $\Lambda$ and descend to
forms on $T^2$, which generate the de-Rham cohomology of the torus in degree $1$, i.e.~$H^1(T^2,\R)=\R\cdot[d\theta]\oplus\R\cdot[dz]$.

Similarly, $\partial/\partial \theta$ and $\partial/\partial z$ descend to Killing-vectorfields on $U_{\varepsilon}(\Sigma)$. To be more specific, $\partial/\partial \theta$ is an infinitesimal rotation around the singular axis and $\partial/\partial z$ an infinitesimal translation along the same axis. Consequently, $\sigma_{\partial/\partial \theta}$ and $\sigma_{\partial/\partial z}$ make up parallel sections of the bundle $\E$, i.e.~$\sigma_{\partial/\partial \theta}, \sigma_{\partial/\partial z}\in H^0(T^2,\E)$.
\begin{lemma} If the cone-angles are $\leq\pi$, then
in the hyperbolic and the spherical case
$$H^0(T^2,\E)=\R\cdot\sigma_{\partial/\partial\theta}\oplus\R\cdot\sigma_{\partial/\partial z}\,.$$
\end{lemma}
\begin{pf}
Let $\lambda$ be the longitudinal and $\mu$ be the meridian loop. Clearly
$H^0(T^2,\E)\cong Z^0(\pi_1
T^2,\g)=\{v\in\g:(Ad\circ\hol(\gamma))v=v\; \forall \gamma \in \pi_1 T^2\}$,
which we view as the infinitesimal centralizer of the holonomy representation
restricted to the torus. 
We compute the centralizer $Z({\hol}(\pi_1T^2))$ in each case.

In the hyperbolic case, let $A={\hol}(\lambda)\in\SL_2(\C)$ and $B={\hol}(\mu)\in\SL_2(\C)$. Since $\hol$ is the holonomy of a hyperbolic cone-manifold structure with
cone-angles $\leq \pi$, we may assume that $A=\diag(\eta,\eta^{-1})$ and
$B=\diag(\xi,\xi^{-1})$ with $\eta, \xi \neq \pm 1$. Then it is easy to see
that $Z({\hol}(\pi_1T^2))=\{\diag(\zeta,\zeta^{-1}), \zeta \in \C^{\ast}\}$,
hence $Z^0(\pi_1T^2,\mathfrak{sl}_2(\C))\cong\R^2$. Since $\sigma_{\partial/\partial\theta}$ and $\sigma_{\partial/\partial z}$ are closed and linearly independent, the result follows.

In the spherical case, ${\hol}:\pi_1T^2\rightarrow \SU(2)\times \SU(2)$ splits
as a product representation ${\hol}=({\hol}_1,{\hol}_2)$ with
${\hol}_i:\pi_1T^2\rightarrow \SU(2)$ for $i\in\{1,2\}$. We then have
$Z({\hol}(\pi_1T^2))=Z({\hol}_1(\pi_1T^2))\times Z({\hol}_2(\pi_1T^2))$. Let
$A_i={\hol}_i(\lambda)\in\SU(2)$ and $B_i={\hol}_i(\mu)\in\SU(2)$. Without loss
of generality we assume that $A_i=\diag(\eta_i,\overline{\eta}_i)$ and
$B_i=\diag(\xi_i,\overline{\xi}_i)$ with $\eta_i,\xi_i \in S^1$. Since $\hol$ is
the holonomy of a spherical cone-manifold structure with cone-angles $\leq \pi$, $\hol(\mu)$ must be a nontrivial rotation. This
implies that
$\{\xi_1,\overline{\xi}_1\}=\{\xi_2,\overline{\xi}_2\}\neq\{\pm1\}$. Then it
follows that $Z({\hol}_i(\pi_1T^2))=\{\diag(\zeta,\overline{\zeta}), \zeta\in
S^1\}$. This implies that $Z^0(\pi_1T^2,\su(2))\cong\R$. As above, $\sigma_{\partial/\partial\theta}$ and $\sigma_{\partial/\partial z}$ provide a basis for $H^0(T^2,\E)$.
\end{pf}\\
\\
We define forms
\begin{align*}
  \omega_{ang}&=d\theta\varotimes\sigma_{\partial/\partial\theta}\\
  \omega_{shr}&=d\theta\varotimes\sigma_{\partial/\partial z}\\
  \omega_{tws}&=d z\varotimes\sigma_{\partial/\partial\theta}\\
  \omega_{len}&=d z\varotimes\sigma_{\partial/\partial z}\,.
\end{align*}
Since $\sigma_{\partial/\partial\theta}$ and $\sigma_{\partial/\partial z}$ are parallel, these forms are closed. These forms will be tangent to the corresponding geometric deformations of the singular tube, i.e.~$\omega_{ang}$ is supposed to change the cone-angle $\alpha$, similarly for $t$ and $l$. $\omega_{shr}$ will be tangent to a deformation, which leads out of the class of cone-metrics (which may be called a ``shearing''-deformation). This will be made precise.

\begin{lemma}
The forms $\omega_{ang}$ and $\omega_{shr}$ are {\it not} $L^2$ on
$U_{\varepsilon}(\Sigma)$, whereas the forms $\omega_{tws}$ and $\omega_{len}$ are
bounded on $U_{\varepsilon}(\Sigma)$ and hence $L^2$.
\end{lemma}
\begin{pf}
The metric on $U_{\varepsilon}(\Sigma)$ is given by
$g=dr^2+\sn^2(r)d\theta^2+\cs^2(r)dz^2$. Hence $dvol=\sn(r)\cs(r)dr\wedge
d\theta\wedge dz$. For a $1$-form $\omega=\alpha\varotimes\sigma _X$
with $\alpha\in\Omega^1(U_{\varepsilon}(\Sigma))$ and
$X\in\Gamma(TU_{\varepsilon}(\Sigma))$ we have $|\omega|^2 = |\alpha|^2 \left(|\nabla X|^2 + |X|^2 \right)$.
Clearly 
$$|d\theta|^2=\frac{1}{\sn^2(r)}\,,\; |dz|^2=\frac{1}{\cs^2(r)}\,,\;
\left|\textstyle\frac{\partial}{\partial \theta}\right|^2=\sn^2(r)\,,\; \left|\textstyle\frac{\partial}{\partial z}\right|^2=\cs^2(r)\,.$$
Let $$\left\{e_1=\textstyle\frac{\partial}{\partial r},\,
e_2=\sn(r)^{-1}\frac{\partial}{\partial \theta},\,
e_3=\cs(r)^{-1}\frac{\partial}{\partial z}\right\}$$ be an orthonormal frame for
$TU_{\varepsilon}(\Sigma)$. 
A
straightforward calculation shows that with respect to this frame
$$
\nabla \textstyle\frac{\partial}{\partial \theta}=\left( \begin{array}{ccc}
0 & -\cs(r) & 0 \\
\cs(r) & 0 & 0 \\
0 & 0 & 0\\ \end{array}\right) \in \Gamma(\mathfrak{so}(TU_{\varepsilon}(\Sigma))
$$
and
$$
\nabla \textstyle\frac{\partial}{\partial z}=\left( \begin{array}{ccc}
0 & 0 & \kappa\sn(r) \\
0 & 0 & 0 \\
-\kappa\sn(r) & 0 & 0\\ \end{array}\right) \in \Gamma(\mathfrak{so}(TU_{\varepsilon}(\Sigma))\,,
$$
such that
$$
\left|\nabla\textstyle\frac{\partial}{\partial \theta}\right|^2=\cs^2(r)\,, \,\left|\nabla\textstyle\frac{\partial}{\partial z}\right|^2=\kappa^2\sn^2(r)\,.  
$$
We obtain
$$
|\omega_{ang}|^2=\frac{\sn^2(r)+\cs^2(r)}{\sn^2(r)}\,,\; |\omega_{shr}|^2=\frac{\cs^2(r)+\kappa^2\sn^2(r)}{\sn^2(r)}
$$
and
$$
|\omega_{tws}|^2=\frac{\sn^2(r)+\cs^2(r)}{\cs^2(r)}\,,\; |\omega_{len}|^2=\frac{\cs^2(r)+\kappa^2\sn^2(r)}{\cs^2(r)}\,.
$$
In the first case we observe that $|\omega_{ang}|^2dvol\sim
|\omega_{shr}|^2dvol \sim \sn(r)^{-1}$, which is not integrable for
$r\in(0,\varepsilon)$. In the second case we find $\omega_{tws}$ and $\omega_{len}$
bounded and therefore $L^2$-integrable.
\end{pf}
\begin{lemma}\label{torus} If the cone-angles are $\leq \pi$, then in the hyperbolic
and the spherical case
$$ H^1(T^2,\E)=\R\cdot[\omega_{ang}]
  \oplus \R\cdot[\omega_{shr}]
   \oplus \R\cdot[\omega_{tws}]
      \oplus \R\cdot[\omega_{len}]\,.$$
\end{lemma}
\begin{pf}
Since $H^0(T^2,\E)=\R\cdot\sigma_{\partial/\partial \theta} \oplus \R\cdot\sigma_{\partial/\partial z}$, we have a short exact sequence of flat vector-bundles
$$
0 \rightarrow \underline{\R}^2 \rightarrow \E \rightarrow \E/\underline{\R}^2 \rightarrow 0\,.
$$
Here we denote by $\underline{\R}^k$ the trivial vector-bundle of real rank
$k$ together with the trivial flat connection. 

We claim that the natural map
$
H^1(T^2,\R^2) \rightarrow H^1(T^2,\E)
$ is an isomorphism. In the spherical
case we can use the parallel metric on $\E$ to split the
short exact
coefficient sequence. Then clearly $H^0(\E/\underline{\R}^2)=0$ and we may use
Poincar\'{e} duality to conclude that $\E/\underline{\R}^2$ is acyclic. Now the result follows from the long exact cohomology sequence.

In the hyperbolic case we can use the parallel Killing form $B$
to split the coefficient sequence, if $B$ restricted to $\underline{\R}^2$ is nondegenerate. We use the
local formula for the Killing form in Lemma \ref{Killing} with $\kappa=-1$:
$$
B(\sigma_X,\sigma_Y)=-4(\nabla X,\nabla Y)+4(X,Y)\,.
$$
From the calculations in the previous lemma we obtain
\begin{align*}
B(\sigma_{\partial/\partial \theta},\sigma_{\partial/\partial
\theta}) &=-4\cosh^2(r)+4\sinh^2(r)=-4\\
B(\sigma_{\partial/\partial z},\sigma_{\partial/\partial
z}) &=-4\sinh^2(r)+4\cosh^2(r)=4\\
B(\sigma_{\partial/\partial \theta},\sigma_{\partial/\partial
z}) &=0\,,
\end{align*}
which shows that $B|_{\underline{\R}^2}$ is nondegenerate. Then the result follows as above.
\end{pf}\\
\\
We wish to calculate the periods of the differential forms $\omega_{ang}, \omega_{shr},
\omega_{tws}$ and $\omega_{len}$. Let $x_0=(0,0)$ be the basepoint of $T^2$. For $\gamma
\in \pi_1T^2$ and $\omega \in \Omega^1(T^2,\E)$ closed, we have a well-defined integral
$$
\int_\gamma \omega = \int_0^1
\tau_{\gamma(t)}^{-1}\omega(\dot{\gamma}(t)) dt,
$$
where $\tau_{\gamma(t)}$ denotes the parallel transport along $\gamma$ from
$x_0=\gamma(0)$ to $\gamma(t)$. Recall that the map $\gamma\mapsto
z_{\omega}(\gamma)=\int_{\gamma}\omega$ defines a group cocycle, if we
identify $\E_{x_0}$ with $\g$.

Note that if $\omega$ is of the form $\omega = \alpha \varotimes \sigma$ with
$\nabla\sigma=0$, then $\int_{\gamma}\omega$ is very easy to compute:
$$
\int_{\gamma}\omega=\int_\gamma \alpha \cdot \sigma_{x_0}.
$$
This remark applies in particular to $\omega_{ang}, \omega_{shr},
\omega_{tws}$ and $\omega_{len}$. We concentrate on the values of the
corresponding group cocycles $z_{ang}, z_{shr},
z_{tws}$ and $z_{len}$ on the meridian $\mu \in \pi_1T^2$, $\mu(0)=x_0$. We obtain
\begin{align*}
z_{ang}(\mu) &= \int_{\mu}\omega_{ang} = \alpha \cdot (\sigma_{\partial /
  \partial \theta})_{x_0}\\
z_{shr}(\mu) &= \int_{\mu}\omega_{shr} = \alpha\cdot (\sigma_{\partial / \partial z})_{x_0}\
\end{align*}
and
$$
z_{tws}(\mu)=z_{len}(\mu)=0\,.
$$

\subsubsection{The higher genus case}
Let $\Sigma$ be a connected graph with trivalent vertices. Then $F_g=\partial
U_{\varepsilon}(\Sigma)$ is a surface of genus $g=(N+3)/3$, where $N$ is the
number of edges contained in $\Sigma$. $U_{\varepsilon}(v)$, the smooth part of the $\varepsilon$-ball around a vertex $v\in\Sigma$, is homotopy equivalent to a pair of pants $P$.
\begin{lemma}
If the cone-angles are $\leq\pi$, then $H^0(F_g,\E)=0$.
\end{lemma}
\begin{pf}
If we restrict the holonomy of $M$ to $U_{\varepsilon}(v)$, the smooth part of the $\varepsilon$-ball around a vertex $v\in\Sigma$, then 
$\hol(\pi_1(U_{\varepsilon}(v))$ fixes a point $p \in
\M$. $U_{\varepsilon}(v)$ deformation-retracts to $P\subset\partial
U_{\varepsilon}(\Sigma)=F_g$. Using the presentation $\pi_1(P)=\langle
\mu_1,\mu_2,\mu_3 \vert \mu_1\mu_2\mu_3=1 \rangle$, we obtain that the
$\hol(\mu_i)$, $i\in\{1,2,3\}$, project to nontrivial rotations with mutually distinct axes. This implies
that $Z(\hol(\pi_1F_g))=\{\pm 1\}$.
\end{pf}
\begin{cor}
If the cone-angles are $\leq\pi$, then in the hyperbolic case
$H^1(F_g,\E)\cong\C^{6g-6}$ and in the spherical case
$H^1(F_g,\E_i)\cong\R^{6g-6}$.
\end{cor}
\begin{pf}
Using the parallel Killing form $B$ on $\E$ in the hyperbolic case, resp.\ the
parallel metric on $\E_i$ in the spherical case, we conclude that $H^2(F_g,\E)=H^2(F_g,\E_i)=0$
using Poincar\'{e} duality. Now for any flat bundle $\F$ over $F_g$ one has $\chi(F_g,\F)=\dim \F \cdot\chi(F_g) = \dim \F\cdot (2-2g)$, 
this implies in particular that $\dim H^1(F_g,\F)=-\dim \F\cdot (2-2g)$ if $H^0(F_g,\F)=H^2(F_g,\F)=0$.
\end{pf}\\
\\
Away from the vertices, the singular locus $U_{\varepsilon}(\Sigma)$ can be
given coordinates $(r,\theta_i,z_i)$ with $r\in(0,\varepsilon)$, $\theta_i \in
\R/\alpha_i\Z$ and $z_i\in (\delta,l_i-\delta)$ for some $\delta >0$. Here $\alpha_i$ is the cone-angle around the
$i$-th edge and $l_i$ its length. Then the metric is given by
$g=dr^2+\sn^2(r)d\theta_i^2+\cs^2(r)dz_i^2$.

We choose a function $\varphi_i=\varphi_i(z_i)$ such that $\varphi_i(\delta)=0$,
$\varphi_i(l_i-\delta)=l_i-\delta$ and
$d\varphi_i\vert_{(\delta,2\delta)}=d\varphi_i\vert_{(l_i-2\delta,l_i-\delta)}=0$.
Then $d\varphi_i \in \Omega^1(U_{\varepsilon}(\Sigma))$ is
well-defined and so are
\begin{align*}
\omega_{tws}^i &= d\varphi_i \varotimes \sigma_{\partial/\partial\theta_i}\\
\omega_{len}^i &= d\varphi_i \varotimes \sigma_{\partial/\partial z_i}\,.
\end{align*}  
Note that these forms are supported away from the vertices of the singularity.
\begin{lemma}
The differential forms $\omega_{tws}^i$ and $\omega_{len}^i$ are bounded on
$U_{\varepsilon}(\Sigma)$, hence in particular $L^2$. 
\end{lemma}
\begin{pf}
This essentially amounts to the same computation as in the torus case.
\end{pf}
\begin{lemma}
The cohomology classes of the closed differential forms 
$\{ \omega_{tws}^1, \omega_{len}^1, \ldots , \omega_{tws}^N,
\omega_{len}^N \}$
are linearly independent in $H^1(F_g, \E)$.
\end{lemma}
\begin{pf}
Suppose we have a nontrivial linear relation between the above classes in
$H^1(F_g,\E)$, say
$$
t_1\omega_{tws}^1+l_1\omega_{len}^1+\ldots+t_N\omega_{tws}^N+l_N\omega_{len}^N
=d\sigma
$$
for some $\sigma \in \Gamma(F_g,\E)$. Since the forms $\omega_{tws}^i$ and
$\omega_{len}^i$ are supported away from the vertices, we obtain $d\sigma=0$ in
a neighbourhood of each vertex $v_i$. A neighbourhood $U_{\varepsilon}(v_i)$ of
a vertex is
homotopy equivalent to the thrice-punctured sphere $P$. Since $H^0(P,\E)=0$, we have
$\sigma\vert_{U_{\varepsilon}(v_i)}=0$ for each vertex. 
Therefore we obtain a nontrivial linear relation on at least one of the tori
$T^2_i=\R^2/\alpha_i\Z+l_i\Z$, where $\sigma_i$ denotes the restriction of
$\sigma$ to a
neighbourhood of the $i$-th edge:
$$
t_i\omega_{tws}^i+l_i\omega_{len}^i = d\sigma_i\,,
$$
which is a contradiction in view of Lemma \ref{torus}, since $d\varphi_i$ is
cohomologous to $dz_i$ on $T^2_i$.
\end{pf}
\subsection{Local structure of the representation variety}

\subsubsection{The torus case}

Let $\iota:T^2\rightarrow M$ be the inclusion of a torus boundary
component. The map $\iota$ induces a group homomorphism $\iota_*:\pi_1T^2\rightarrow \pi_1M$ and hence a map $\iota^*:R(\pi_1M,G)\rightarrow R(\pi_1T^2,G)$ for $G=\SL_2(\C)$ or $\SU(2)$ respectively.
\begin{lemma}
$\rho=\iota_{T^2}^*\hol: \pi_1T^2\rightarrow \SL_2(\C)$ is a smooth point of $R(\pi_1T^2,\SL_2(\C))$. The local $\C$-dimension of $R(\pi_1T^2,\SL_2(\C))$ around $\rho$ equals $4$. $T_{\rho}R(\pi_1T^2,\SL_2(\C))$ may be identified with $Z^1(\pi_1T^2,\sl_2(\C))$. 
\end{lemma}
\begin{pf}
We identify $R(\pi_1T^2,\SL_2(\C))$ with the (affine algebraic) set $\{(A,B)\in \SL_2(\C)\times \SL_2(\C)\vert\, [A,B]=1\}$. The kernel of the differential of the commutator map $\ker d_{(A,B)}[\,\cdot\,,\cdot\,]$ may be identified with the space of $1$-cocycles $Z^1(\pi_1T^2,\sl_2(\C))$. We have $\dim_{\C} Z^1(\pi_1T^2,\sl_2(\C))=4$ from the cohomology computations. Note that this implies that $d_{(A,B)}[\,\cdot\,,\cdot\,]$ is {\it not} surjective at $(A,B)=\rho$. Without loss of generality we may assume that $\rho=(\diag(\lambda,\lambda^{-1}),\diag(\mu,\mu^{-1}))$ with $\lambda,\mu \in \C^*$. We define a map
\begin{align*}
F: \C^*\times\C^*\times \SL_2(\C) &\longrightarrow \SL_2(\C)\times \SL_2(\C)\\
(\lambda,\mu,A) &\longmapsto \left( A\diag(\lambda,\lambda^{-1})A^{-1},A\diag(\mu,\mu^{-1})A^{-1}\right)\,.
\end{align*}
We claim that $\rank_{\C}F=4$ at $(\lambda,\mu,1)$. The image of $F$ is
certainly contained in $R(\pi_1T^2,\SL_2(\C))$, such that an easy application of the
implicit function theorem (cf.~\cite{Wei}, \cite[Lemma 6.8]{Rag}) yields the
result. Consider the standard $\C$-basis
of $\sl_2(\C)$: 
$$
\left\{x=\left(\begin{array}{cc}0 & 1\\ 0 & 0\end{array}\right),
h=\left(\begin{array}{cc}1 & 0\\ 0 & -1\end{array}\right), y=\left(\begin{array}{cc}0 & 0\\ 1 & 0\end{array}\right)\right\}.
$$
Clearly $\C\cdot h$ exponentiates to
$Z(\rho(\pi_1T^2))=\{\diag(\eta,\eta^{-1})\vert\eta \in\C^*\}$, the
stabilizer of $\rho$ under the conjugation action of $\SL_2(\C)$. Now it is easily
verified that 
$$\left\{dF(1,0,0),dF(0,1,0),dF(0,0,x),dF(0,0,y)\right\}_{(\lambda,\mu,1)}$$
are linearly
independent if $\lambda\neq\pm 1$ or $\mu\neq\pm 1$. This implies that $\rank_{\C}F$ at $(\lambda,\mu,1)$ is at least $4$, but since
$\im d_{(\lambda,\mu,1)}F\subset Z^1(\pi_1T^2,\sl_2(\C))$, it has to equal $4$. 
\end{pf}
\begin{cor} $\chi=[\iota_{T^2}^*\hol]$ is a smooth point of
$X(\pi_1T^2,\SL_2(\C))$. The local $\C$-dimension of $X(\pi_1T^2,\SL_2(\C))$ around $\chi$ equals $2$. Furthermore $T_{\chi}X(\pi_1T^2,\SL_2(\C))$ may be identified with $H^1(\pi_1T^2,\sl_2(\C))$.
\end{cor}
\begin{pf}
The restriction of $F$ to $\C^*\times\C^*\times\{1\}$ provides a local slice to
the action through $\rho$, upon which the stabilizer of $\rho$ acts trivially. The tangent space to the orbit through
$\rho$ may be
identified with $B^1(\pi_1T^2,\sl_2(\C))$. We know that $\dim_{\C}H^1(\pi_1T^2,\sl_2(\C))=2$ from the cohomology
computations.
\end{pf}\\
\\
For $\gamma \in \Gamma$ we define a function $t_{\gamma}:
R(\Gamma,\SL_2(\C))\rightarrow \C$ by $t_{\gamma}(\rho)=\tr \rho(\gamma)$. If
$\rho$ is a smooth point of $R(\Gamma,\SL_2(\C))$, then $t_{\gamma}$ is
smooth near $\rho$.
Since $\tr$ is invariant under conjugation, $t_{\gamma}$
descends to a map on the quotient $X(\Gamma,\SL_2(\C))$, which we again refer to
as $t_{\gamma}$. If
$\chi=[\rho]$ is a smooth point of $X(\Gamma,\SL_2(\C))$, then $t_{\gamma}$ is
smooth near $\chi$.\\
\\
Let $\rho=\iota_{T^2}^*\hol$ and let $z \in Z^1(\pi_1T^2,{\sl}_2(\C))$ be
given. If we have a deformation of $\rho$, i.e.~a family of representations
$\rho_t:\pi_1T^2\rightarrow \SL_2(\C)$ with $\rho_0=\rho$, which is tangent to
$z$, i.e.~$z(\gamma)=\left.\frac{d}{dt}\right\vert_{t=0} \rho_t(\gamma)\rho(\gamma)^{-1}$ for all
$\gamma \in \pi_1T^2$, we have that
the infinitesimal change of the trace of $\rho(\gamma)$ is given as
$$
dt_{\gamma}(z)=\textstyle\left.\frac{d}{dt}\right\vert_{t=0}\tr \rho_t(\gamma) = 
\tr \left( z(\gamma)\rho(\gamma) \right).
$$
We wish to apply this to $z_{ang}, z_{shr}, z_{tws}$ and $z_{len}$. Let
$\mu \in \pi_1T^2$ be the meridian and $\lambda \in \pi_1T^2$ the longitude.
We assume that
$$
\rho(\lambda)=\left( \begin{array}{cc}
\eta & 0\\
0 & \eta^{-1} \end{array}\right)\in \SL_2(\C)
$$
and
$$
\rho(\mu)=\left( \begin{array}{cc}
\xi & 0\\
0 & \xi^{-1} \end{array}\right)\in \SL_2(\C)
$$
with $\eta,\xi \neq \pm 1$. Then $\rho$ preserves the axis $\gamma=\{0\}\times\R_+\subset \Hyp^3$, if
we work in the upper half-space model $\Hyp^3=\C \times \R_+$. If we use cylindrical coordinates
$(r,\theta,z)$ around $\gamma$,
then we have already observed that
\begin{align*}
\sigma_{\partial/\partial \theta} &= \frac{1}{2}\left( \begin{array}{cc}
i & 0\\
0 & -i \end{array}\right)\in \sl_2(\C)
\end{align*}
and
\begin{align*}
\sigma_{\partial/\partial z} &=  \frac{1}{2}\left( \begin{array}{cc}
1 & 0\\
0 & -1 \end{array}\right)\in \sl_2(\C)\,.
\end{align*}
Let us concentrate on the value of the cocycles $z_{ang}, z_{shr}, z_{tws}$ and $z_{len}$
on the meridian $\mu \in \pi_1 T^2$.
We obtain
\begin{align*}
z_{ang}(\mu) &= \frac{\alpha}{2}\left( \begin{array}{cc}
i & 0\\
0 & -i \end{array}\right)\in \sl_2(\C)\\
z_{shr}(\mu) &= \frac{\alpha}{2}\left( \begin{array}{cc}
1 & 0\\
0 & -1 \end{array}\right)\in \sl_2(\C)\,,
\end{align*}
while
$$
z_{tws}(\mu)=z_{len}(\mu)=0\,.
$$
As a consequence we obtain for the infinitesimal change of trace
\begin{align*}
dt_{\mu}(z_{ang}) &= (i\alpha/2)(\xi-\xi^{-1})\in\C\\
dt_{\mu}(z_{shr}) &= (\alpha/2)(\xi-\xi^{-1})\in\C\,,
\end{align*}
while
$$
dt_{\mu}(z_{tws})=dt_{\mu}(z_{len})=0\,.
$$
Note that $\xi-\xi^{-1} \neq 0$ since $\xi \neq \pm 1$. Since the cohomology classes
of the cocycles $\{z_{ang}, z_{shr}, z_{tws},z_{len}\}$ provide a $\R$-basis of
$H^1(\pi_1T^2,\sl_2(\C))$, we obtain as a
consequence of the above
 calculations:
\begin{lemma}
The function $t_{\mu}$ has $\C$-rank $1$ at $\chi=[\iota_{T^2}^*
\hol]$. In particular, the level-set $V=\{t_{\mu}\equiv t_{\mu}(\chi)\}$ is locally around
$\chi$ a smooth, half-dimensional submanifold of $X(\pi_1T^2, \SL_2(\C))$. Furthermore,
the cohomology class of the cocycle $z_{len}$ provides a $\C$-basis for $T_{\chi}V$. The cohomology classes of the
cocycles $\{z_{tws},z_{len}\}$ provide a $\R$-basis of $T_{\chi}V$. 
\end{lemma}
We now turn to the spherical case.
\begin{lemma}
Let $\rho_i=\iota_{T^2}^*\hol_i: \pi_1T^2\rightarrow
\SU(2)$. Then $\rho_i$ is a smooth point of $R(\pi_1T^2,\SU(2))$. The local $\R$-dimension of $R(\pi_1T^2,\SU(2))$ around $\rho_i$ equals $4$. Furthermore $T_{\rho_i}R(\pi_1T^2,\SU(2))$ may be identified with $Z^1(\pi_1T^2,\su(2))$. 
\end{lemma}
\begin{pf}
As above we define a map
\begin{align*}
F: S^1\times S^1\times \SU(2) &\longrightarrow \SU(2)\times \SU(2)\\
(\lambda,\mu,A) &\longmapsto \left( A\diag(\lambda,\lambda^{-1})A^{-1},A\diag(\mu,\mu^{-1})A^{-1}\right)
\end{align*}
We consider the standard $\R$-basis of $\su(2)$:
$$
\left\{i=\left(\begin{array}{cc}i & 0\\ 0 & -i\end{array}\right),
j=\left(\begin{array}{cc}0 & 1\\ -1 & 0\end{array}\right), k=\left(\begin{array}{cc}0 & i\\ i & 0\end{array}\right)\right\}.
$$
Now $\R\cdot i$ exponentiates to
$Z(\rho(\pi_1T^2))=\{\diag(\eta,\eta^{-1})\vert\eta \in S^1\}$, the
stabilizer of $\rho$ under the conjugation action of $\SU(2)$. It is easily
verified that 
$$\left\{dF(1,0,0),dF(0,1,0),dF(0,0,j),dF(0,0,k)\right\}_{(\lambda,\mu,1)}$$
 are linearly
independent if $\lambda\neq\pm 1$ or $\mu\neq\pm 1$. The
result follows as above.
\end{pf}
\begin{cor} $\chi_i=[\iota_{T^2}^*\hol_i]$ is a smooth point of $X(\pi_1T^2,\SU(2))$. The local $\R$-dimension of $X(\pi_1T^2,\SU(2))$ around $\chi_i$ equals $2$. Furthermore $T_{\chi_i}X(\pi_1T^2,\SU(2))$ may be identified with $H^1(\pi_1T^2,\su(2))$.
\end{cor}
\begin{pf}
The restriction of $F$ to $S^1\times S^1\times\{1\}$ provides a local slice to
the action through $\rho_i$, upon which the stabilizer of $\rho$ acts trivially. The tangent space to the orbit through
$\rho_i$ may be
identified with $B^1(\pi_1T^2,\su(2))$. From the cohomology
computations we have $\dim_{\R}H^1(\pi_1T^2,\su(2))=2$.
\end{pf}\\
\\
For $\gamma \in \Gamma$ we define a function $t_{\gamma}:
R(\Gamma,\SU(2))\rightarrow \R$ by $t_{\gamma}(\rho)=\tr \rho(\gamma)$. If
$\rho$ is a smooth point of $R(\Gamma,\SU(2))$, then $t_{\gamma}$ is
smooth near $\rho$.
Since $\tr$ is invariant under conjugation, $t_{\gamma}$
descends to a map on the quotient $X(\Gamma,\SU(2))$, which we again refer to
as $t_{\gamma}$. If
$\chi=[\rho]$ is a smooth point of $X(\Gamma,\SU(2))$, then $t_{\gamma}$ is
smooth in a neighbourhood of $\chi$.

For a representation $\rho=(\rho_1,\rho_2):\Gamma \rightarrow \SU(2) \times \SU(2)$
and $\gamma\in\Gamma$ let $T^i_{\gamma}(\rho)=t_{\gamma}(\rho_i)$. This defines
an $\R^2$-valued function $T_{\gamma}=(T_{\gamma}^1,T_{\gamma}^2)$ on $R(\Gamma,\SU(2)\times\SU(2))$, which we view as a ``complex'' trace function.\\
\\
Let $\rho=\iota_{T^2}^*\hol$ and let $z=(z_1,z_2) \in Z^1(\pi_1T^2,\su(2)\oplus\su(2))$ be
given.
The infinitesimal change of the trace of $\rho(\gamma)$ is given as
$$
dT_{\gamma}(z)=\left( dt_{\gamma}(z_1), dt_{\gamma}(z_2) \right)
$$
We wish to apply this to $z_{ang}, z_{shr}, z_{tws}$ and $z_{len}$. Let
$\lambda \in \pi_1T^2$ be the meridian and $\mu \in \pi_1T^2$ the longitude.
We assume that
$$
\rho(\lambda)=\left( \left( \begin{array}{cc}
\eta_1 & 0\\
0 & \overline{\eta}_1 \end{array}\right), \left( \begin{array}{cc}
\eta_2 & 0\\
0 & \overline{\eta}_2 \end{array}\right) \right) \in \SU(2)\times\SU(2)
$$
and
$$
\rho(\mu)=\left( \left( \begin{array}{cc}
\xi_1 & 0\\
0 & \overline{\xi}_1 \end{array}\right), \left( \begin{array}{cc}
\xi_2 & 0\\
0 & \overline{\xi}_2 \end{array}\right) \right) \in \SU(2)\times\SU(2)
$$
with $\xi_1=\xi_2=:\xi$ and $\xi \neq \pm 1$, since $\rho(\mu)$ is a nontrivial rotation. Then $\rho$ preserves the pair of axes $\{\gamma, \gamma^{\perp}\}$, where $\gamma=\C\cap\Sph^3$ and $\gamma^{\perp}=\C j\cap\Sph^3$. If we use cylindrical coordinates
$(r,\theta,z)$ around $\gamma$,
then we have already observed that
\begin{align*}
\sigma_{\partial/\partial \theta} 
&=\left(  \frac{1}{2}\left( \begin{array}{cc}
i & 0\\
0 & -i \end{array}\right),  \frac{1}{2}\left( \begin{array}{cc}
i & 0\\
0 & -i \end{array}\right) \right)
 \in \su(2)\oplus \su(2)
\end{align*}
and
\begin{align*}
\sigma_{\partial/\partial z} 
&= \left( \frac{1}{2}\left( \begin{array}{cc}
i & 0\\
0 & -i \end{array}\right), \frac{1}{2}\left( \begin{array}{cc}
-i & 0\\
0 & i \end{array}\right) \right)
\in \su(2)\oplus \su(2)\,.
\end{align*}
In particular, this implies that $\sigma_{\partial/\partial \theta}+\sigma_{\partial/\partial z}\in\Gamma(U_{\varepsilon}(\Sigma),\E_1)$,
and on the other hand $\sigma_{\partial/\partial \theta}-\sigma_{\partial/\partial z}\in\Gamma(U_{\varepsilon}(\Sigma),\E_2)$. Therefore we have
$$
\omega_{tws}+\omega_{len} \in \Omega^1(U_{\varepsilon}(\Sigma),\E_1)
$$
and
$$
\omega_{tws}-\omega_{len} \in \Omega^1(U_{\varepsilon}(\Sigma),\E_2)\,.
$$
Again, we concentrate on the value of the cocycles $z_{ang}, z_{shr}, z_{tws}$ and $z_{len}$
on the meridian $\mu \in \pi_1 T^2$.
We obtain
\begin{align*}
z_{ang}(\mu) &= \left( \frac{\alpha}{2} \left( \begin{array}{cc}
i & 0\\
0 & -i \end{array}\right), \frac{\alpha}{2}\left( \begin{array}{cc}
i & 0\\
0 & -i \end{array}\right) \right) \in \su(2)\oplus\su(2)\\
z_{shr}(\mu) &= \left(\frac{\alpha}{2} \left( \begin{array}{cc}
i & 0\\
0 & -i \end{array}\right), \frac{\alpha}{2}\left( \begin{array}{cc}
-i & 0\\
0 & i \end{array}\right) \right) \in \su(2)\oplus\su(2) \,,
\end{align*}
while
$$
z_{tws}(\mu)=z_{len}(\mu)=0\,.
$$
As a consequence we obtain for the infinitesimal change of trace
\begin{align*}
dT_{\mu}(z_{ang}) &= \alpha(-\Im \xi, -\Im \xi)\in\R^2\\
dT_{\mu}(z_{shr}) &= \alpha(-\Im \xi, +\Im \xi)\in\R^2\,,
\end{align*}
while
$$
dT_{\mu}(z_{tws})=dT_{\mu}(z_{len})=0\,.
$$
Note that $\Im \xi=\frac{1}{2i}( \xi-\overline{\xi} ) \neq 0$ since $\xi \neq \pm 1$. Since the cohomology classes
of the cocycles $\{z_{ang}, z_{shr}, z_{tws},z_{len}\}$ provide a $\R$-basis of
$H^1(\pi_1T^2,\su(2)\oplus\su(2))$, we obtain as a
consequence of the above
 calculations:
\begin{lemma}
The function $t_{\mu}$ has $\R$-rank $1$ at $\chi_i=[\iota_{T^2}^*
\hol_i]$. In particular, the level-set $V_i=\{t_{\mu}\equiv t_{\mu}(\chi_i)\}$ is locally around
$\chi_i$ a smooth, half-dimensional submanifold of $X(\pi_1T^2, \SU(2))$. Furthermore
the cohomology class of the cocycle $z_{tws}+z_{len}$ provides a $\R$-basis of $T_{\chi_1}V_1$, the cohomology class of the cocycle $z_{tws}-z_{len}$ provides a $\R$-basis of $T_{\chi_2}V_2$.
\end{lemma}

\subsubsection{The higher genus case}
Let $\iota:F_g\rightarrow M$ be the inclusion of a boundary component of higher
genus $g\geq 2$. $\iota$ induces a group homomorphism $\iota_*:\pi_1F_g\rightarrow \pi_1M$ and a map $\iota^*:R(\pi_1M,G)\rightarrow R(\pi_1F_g,G)$ for $G=\SL_2(\C)$ or $\SU(2)$ respectively.
\begin{lemma} Let $\rho: \pi_1F_g \rightarrow \SL_2(\C)$ be an irreducible representation. Then $\rho$ is a smooth point of $R(\pi_1F_g,\SL_2(\C))$. The local $\C$-dimension of $R(\pi_1F_g,\SL_2(\C))$ around $\rho$ equals $6g-3$. $T_{\rho}R(\pi_1F_g,\SL_2(\C))$ may be identified with $Z^1(\pi_1F_g,\sl_2(\C))$. 
\end{lemma}
\begin{pf}
We identify $R(\pi_1F_g,\SL_2(\C))$ with the (affine algebraic) set 
$$
\left\{(A_1,B_1,\ldots,A_g,B_g)\in \SL_2(\C)^{2g}\vert
f(A_1,B_1,\ldots,A_g,B_g)=1\right\},
$$ 
where $f(A_1,B_1,\ldots,A_g,B_g)=[A_1,B_1]\cdot\ldots\cdot[A_g,B_g]$.
$\ker d_{\rho}f$ may be identified with the space of $1$-cocycles
$Z^1(\pi_1F_g,\sl_2(\C))$. We know that $\dim_{\C}
H^1(\pi_1F_g,\sl_2(\C))=6g-6$ from the cohomology computations. Since 
$\rho$ is irreducible, we have $Z^0(\pi_1F_g,\sl_2(\C))=0$, which implies $\dim_{\C}
Z^1(\pi_1F_g,\sl_2(\C))=6g-3$. Hence
$\rank_{\C}d_{\rho}f=3$, i.e.~$d_{\rho}f$ is surjective. Now the implicit
function theorem implies that $R(\pi_1F_g,\SL_2(\C))$ is smooth at $\rho$ with $T_{\rho}R(\pi_1F_g,\SL_2(\C))=Z^1(\pi_1F_g,\sl_2(\C))$. 
\end{pf}
\begin{cor}
Let $\rho=\iota_{F_g}^*\hol: \pi_1F_g\rightarrow \SL_2(\C)$. Then $\rho$ is a smooth point of $R(\pi_1F_g,\SL_2(\C))$. The local $\C$-dimension of $R(\pi_1F_g,\SL_2(\C))$ around $\rho$ equals $6g-3$. Furthermore $T_{\rho}R(\pi_1F_g,\SL_2(\C))$ may be identified with $Z^1(\pi_1F_g,\sl_2(\C))$. 
\end{cor}
\begin{pf}
Clearly $\rho$ is irreducible: If $v\in \Sigma$ is a singular vertex and we restrict $\hol$ further to $U_{\varepsilon}(v)$, which deformation-retracts to a pair of pants $P\subset F_g$, then $\iota_P^*\hol$ preserves a point $p \in\Hyp^3$. Now if $\rho$ was reducible, then $\iota_P^*\hol$ would preserve a geodesic, which is a contradiction. 
\end{pf}\\
\\
The following is a well-known fact about the action of $\SL_2(\C)$ on the
irreducible part of $R(\Gamma, \SL_2(\C))$, for convenience of the reader we give a proof:
\begin{lemma}\label{Leeb}
The action of $\SL_2(\C)$ on $R_{irr}(\Gamma, \SL_2(\C))$ is proper.
\end{lemma}
\begin{pf}
Let $X$ be a $G$-space. If we have a continuous $G$-equivariant map from $X$ to a proper
$G$-space $Y$, then $X$ itself will be a proper $G$-space.
We construct a continuous, $\SL_2(\C)$-equivariant map
\begin{align*}
R_{irr}(\Gamma, \SL_2(\C)) & \longrightarrow \Hyp^3\\
\rho & \longmapsto \cen(\rho)\,,
\end{align*}
where the "center" of a representation will be the point in $\Hyp^3$, which is
displaced the least in average by the generators of the group. More precisely,
let us fix a presentation $\langle \gamma_1, \ldots, \gamma_n\vert (r_i)_{i\in
I}\rangle$ of $\Gamma$. Note that the (modified) displacement function of $A \in
\SL_2(\C)$
\begin{align*}
\delta_A: \Hyp^3  &\longrightarrow \R\\
x &\longmapsto \cosh d(x,Ax)-1
\end{align*}
is a convex function in general. It is strictly convex if $A$ is parabolic. If
$A$ is semisimple, it is strictly convex along any geodesic different from the axis
of $A$. We define
$$
f_{\rho}(x)=\frac{1}{n}\sum_{i=1}^n \delta_{\rho(\gamma_i)}\,.
$$
If we have a sequence $x_n \in \Hyp^3$ which converges to $x_{\infty} \in
\partial_{\infty}\Hyp^3$, then since $\rho$ is irreducible, there has to be
at least one $\rho(\gamma_i)$ that does not fix $x_{\infty}$. Then it follows
that $\delta_{\rho(\gamma_i)}(x_n)\rightarrow \infty$. Therefore
$f_{\rho}$ is proper.
If we take any geodesic $\gamma$, again since $\rho$ is
irreducible, there has to be at least one $\rho(\gamma_i)$ such that
$\delta_{\rho(\gamma_i)}$ is strictly convex along $\gamma$. Therefore $f_{\rho}$ is strictly
convex.

As a proper and strictly convex function, $f_{\rho}$ assumes its minimum at a unique
point in $\Hyp^3$, which we define to be the center of $\rho$.

If we have a sequence of representations $\rho_n$ converging to $\rho$ with respect to the
compact-open topology on $R_{irr}(\Gamma,\SL_2(\C))$, then $f_{\rho_n}$ converges
to $f_{\rho}$ uniformly on compact sets. Therefore the map $\cen$
is continuous.
Since $\delta_{BAB^{-1}}(x)=\delta_A(B^{-1}x)$ we obtain that
the map $\cen$ is
$\SL_2(\C)$-equivariant.

This together with the fact that the action of $\SL_2(\C)$ on $\Hyp^3$ is
proper proves the lemma.
\end{pf}
\begin{cor} $\chi=[\iota_{F_g}^*\hol]$ is a smooth point of
$X(\pi_1F_g,\SL_2(\C))$. The local $\C$-dimension of $X(\pi_1F_g,\SL_2(\C))$ around
$\chi$ equals $6g-6$. $T_{\chi}X(\pi_1F_g,\SL_2(\C))$ may be identified with
$H^1(\pi_1F_g,\sl_2(\C))$.
\end{cor}
\begin{pf}
Since the action of $\SL_2(\C)$ is proper, we have a local slice to the action. We recall that the
stabilizer of $\rho$, $Z(\rho(\pi_1F_g))$, equals $\{\pm 1\}$. Therefore $X(\pi_1F_g,\SL_2(\C))$ is locally around $\chi$ the quotient of a free $\PSL_2(\C)$
action and therefore smooth.
\end{pf}\\
\\
The meridian curves around the singularity give rise to a pair-of-pants
decomposition of $F_g$. Let $\{\mu_1,\ldots,\mu_N\}$ be the family of
meridians, where $N=3g-3$. This may be used to give an alternative construction of
$R(\pi_1(F_g),\SL_2(\C))$, which is better suited for certain purposes.

Let $P$ denote the thrice-punctured sphere, i.e.~a pair of pants. The fundamental
group of $P$ is the free group on $2$ generators. We will use the following
slightly redundant presentation:
$$
\pi_1P=\langle \mu_1,\mu_2,\mu_3 \vert \mu_1\mu_2\mu_3 = 1 \rangle\,.
$$ 
It follows that 
$$
R(\pi_1P,\SL_2(\C))=\{(A_1,A_2,A_3) \in \SL_2(\C)^3 \vert A_1A_2A_3=1\}.
$$
Clearly the map $f:\SL_2(\C)^3 \rightarrow \SL_2(\C), (A_1,A_2,A_3) \mapsto
A_1A_2A_3$ is a submersion, such that $R(\pi_1P,\SL_2(\C))=f^{-1}(1)$ is a
smooth submanifold of $\C$-dimension $6$. 

Let $\iota_i: S^1 \rightarrow P$ be the inclusion of the $i$-th boundary
circle. Then the induced map $\iota_i^*: R(\pi_1P,\SL_2(\C)) \rightarrow
R(\pi_1S^1,\SL_2(\C))$ corresponds to the projection $pr_i : R(\pi_1P,\SL_2(\C)) \rightarrow
\SL_2(\C), (A_1,A_2,A_3) \mapsto A_i$, which is also a submersion.

The verification of the following statement is elementary and left to the reader:
\begin{lemma}\label{pants}
Let $\rho=\iota_P^*\hol$ be the restriction of the holonomy of a hyperbolic
cone-manifold structure to a pair of pants $P$. Then the differentials $\{ dt_{\mu_1}, dt_{\mu_2}, dt_{\mu_3} \}$ are $\C$-linearly
independent in $T^*_{\rho}R(\pi_1P,\SL_2(\C))$.
\end{lemma}
Since $\rho=\iota_P^*\hol$ is irreducible, we can use Lemma \ref{Leeb} to conclude that
$\chi=[\rho]$ is a smooth point in $X(\pi_1P,\SL_2(\C))$. The local
$\C$-dimension of $X(\pi_1P,\SL_2(\C))$ around $\chi$ is $3$. The functions
$\{t_{\mu_1},t_{\mu_2},t_{\mu_3}\}$ are local holomorphic coordinates on
$X(\pi_1P,\SL_2(\C))$ near $\chi$.\\
\\
We build up $R(\pi_1F_g,\SL_2(\C))$ from $R(\pi_1P,\SL_2(\C))$ using two
basic operations:
\begin{enumerate}
\item glue a pair of pants $P$ to a connected surface with boundary $S$ along a boundary
circle, call the resulting connected surface $S'$
\item glue a connected surface $S$ along two different boundary circles, call the resulting connected surface $S'$
\end{enumerate} 
In the first case $\pi_1S'=\pi_1S\amalg_{\pi_1S^1} \pi_1P$ by van Kampen's
theorem and we have
\begin{align*}
R(\pi_1S',\SL_2(\C))
& =  R(\pi_1S,\SL_2(\C)) \times_{R(\pi_1S^1,\SL_2(\C))} R(\pi_1P,\SL_2(\C))
\end{align*}
via the maps $$\iota^*_{S^1 \hookrightarrow S}:R(\pi_1S,\SL_2(\C)) \rightarrow R(\pi_1S^1,\SL_2(\C))$$
and $$\iota^*_{S^1 \hookrightarrow P}:R(\pi_1P,\SL_2(\C)) \rightarrow
R(\pi_1S^1,\SL_2(\C))\,,$$ 
which will be
transversal since the latter one is a submersion. Therefore $\rho=\iota_{S'}^*\hol$
is a smooth point in $R(\pi_1S',\SL_2(\C))$ since $\rho_S=\iota^*_S\hol$ is a smooth point
in $R(\pi_1S,\SL_2(\C))$ and  $\rho_P=\iota^*_P\hol$ is a smooth point
in $R(\pi_1P,\SL_2(\C))$.\\
\\
In the second case $\pi_1S'$ splits as an HNN-extension of $\pi_1S$. More precisely, if
$\mu_1, \mu_2 \in \pi_1S$ are the loops around the boundary circles, which will be
identified, then $\pi_1S'=\langle \pi_1S, \lambda \vert \lambda \mu_1
\lambda^{-1}=\mu_2\rangle$.  In this case we have
\begin{align*}
R(\pi_1S',\SL_2(\C)) &= \{ (\rho_S, B)
\vert B\rho_S(\mu_1)B^{-1}=\rho_S(\mu_2) \}\\ 
&\subset R(\pi_1S,\SL_2(\C)) \times \SL_2(\C)
\end{align*}
as a consequence. We show that the map 
\begin{align*}
f:R(\pi_1S,\SL_2(\C))\times \SL_2(\C) & \longrightarrow \SL_2(\C)\\
(\rho_S,B) & \longmapsto B\rho_S(\mu_1)B^{-1}\rho_S(\mu_2)^{-1}
\end{align*} 
is a submersion near $\rho=\iota_{S'}^*\hol$. This implies that $\rho=(\rho_S,B)$ is
a smooth point in $R(\pi_1S',\SL_2(\C))$.

Surjectivity of $df$ at $\rho$ can be established as follows:
Let $A_1=\rho_S(\mu_1)$ and $A_2=\rho_S(\mu_2)$. Clearly the map $B\mapsto BA_1B^{-1}A_2^{-1}$ has $\C$-rank $2$.
Since $\{dt_{\mu_1}, dt_{\mu_2}\}$ are linearly independent, we can construct a deformation $t \mapsto (\rho_S)_t$ with $(\rho_S)_t(\mu_2)=A_2$ and $dt_{\mu_1}(\dot{\rho}_S) \neq 0$. This deformation will be transverse to $\im(B \mapsto BA_1B^{-1}A_2^{-1})$.\\
\\
From the construction given above the following is immediate:
\begin{lemma}
The differentials $\{dt_{\mu_1},\ldots,dt_{\mu_N}\}$ with $N=3g-3$ are linearly
independent over $\C$ in $T^*_{\rho}R(\pi_1F_g,\SL_2(\C))$ for $\rho=\iota_{F_g}^*\hol$.
\end{lemma}
Clearly
$$
z_{tws}^i(\mu_j)=\int_{\mu_j}\omega_{tws}^i=0
$$
and 
$$
z_{len}^i(\mu_j)=\int_{\mu_j}\omega_{len}^i=0\,.
$$
Therefore 
$$
dt_{\mu_j}(z_{tws}^i)=0
$$
and
$$
dt_{\mu_j}(z_{len}^i)=0\,.
$$
As a consequence of this we obtain the following lemma.
\begin{lemma}
The level-set $V=\{t_{\mu_1}\equiv t_{\mu_1}(\chi), \ldots,
t_{\mu_N}\equiv t_{\mu_N}(\chi) \}$ is locally around $\chi=[\iota^*_{F_g}\hol]$ a smooth,
half-dimensional submanifold of $X(\pi_1F_g,\SL_2(\C))$. Furthermore, the cohomology classes
of the cocycles 
$\{z_{len}^1,\ldots,z_{len}^N\}$ provide a $\C$-basis of $T_{\chi}V$. The cohomology classes of the cocycles 
$\{z_{tws}^1,z_{len}^1,\ldots,z_{tws}^N,z_{len}^N\}$ provide a $\R$-basis for $T_{\chi}V$.
\end{lemma}
We now turn to the spherical case.
\begin{lemma} Let $\rho: \pi_1F_g \rightarrow \SU(2)$ be an irreducible representation. Then $\rho$ is a smooth point of $R(\pi_1F_g,\SU(2))$. The local $\R$-dimension of $R(\pi_1F_g,\SU(2))$ around $\rho$ equals $6g-3$. $T_{\rho}R(\pi_1F_g,\SU(2))$ may be identified with $Z^1(\pi_1F_g,\su(2))$. 
\end{lemma}
\begin{pf}
This follows as in the case of $\SL_2(\C)$ from the cohomology computations and
the implicit function theorem.
\end{pf}
\begin{cor}
Let $\rho_i=\iota_{F_g}^*\hol_i: \pi_1F_g\rightarrow
\SU(2)$. Then $\rho_i$ is a smooth point of $R(\pi_1F_g,\SU(2))$. The local $\R$-dimension of $R(\pi_1F_g,\SU(2))$ around $\rho_i$ equals $6g-3$. Furthermore $T_{\rho_i}R(\pi_1F_g,\SU(2))$ may be identified with $Z^1(\pi_1F_g,\su(2))$. 
\end{cor}
\begin{pf}
Clearly the $\rho_i$ are both irreducible: If $v\in \Sigma$ is a singular vertex and we restrict $\hol=(\hol_1,\hol_2)$ further to $U_{\varepsilon}(v)$, which deformation-retracts to a pair of pants $P\subset F_g$, then $\iota_P^*\hol$ preserves a point $p \in \Sph^3$. Without loss of generality we may assume that $p=1\in\Sph^3\subset\H$. Then since
$$
\Stab_{\SU(2)\times\SU(2)}(1)=\{(A,A) : A\in\SU(2) \}\,,
$$
we obtain that $\iota_P^*\hol_1=\iota_P^*\hol_2$. Now if $\rho_1$ or $\rho_2$ were reducible, then $\iota_P^*\hol$ would preserve a geodesic, which is a contradiction. 
\end{pf}
\begin{cor} $\chi_i=[\iota_{F_g}^*\hol_i]$ is a smooth point of $X(\pi_1F_g,\SU(2))$. The local $\R$-dimension of $X(\pi_1F_g,\SU(2))$ around $\chi_i$ equals $6g-6$. $T_{\chi_i}X(\pi_1F_g,\SU(2))$ may be identified with $H^1(\pi_1F_g,\su(2))$.
\end{cor}
\begin{pf}
Since the group $\SU(2)$ is compact, the properness of the action is granted. We recall that the stabilizer of $\rho_i$, $Z(\rho_i(\pi_1F_g))$,
equals $\{\pm 1\}$. Therefore $X(\pi_1F_g,\SU(2))$ is near $\chi_i$ a
quotient of a free $\PSU(2)$ action and therefore smooth.
\end{pf}
\begin{lemma}
The differentials $\{dt_{\mu_1},\ldots,dt_{\mu_N}\}$ with $N=3g-3$ are linearly
independent in $T^*_{\rho_i}R(\pi_1F_g,\SU(2))$ for $\rho_i=\iota_{F_g}^*\hol_i$. 
\end{lemma}
\begin{pf}
The arguments in the hyperbolic case apply without essential change.
\end{pf}\\
\\
We obtain finally:
\begin{lemma}
The level-set $V_i=\{t_{\mu_1}\equiv t_{\mu_1}(\chi_i), \ldots,
t_{\mu_N}\equiv t_{\mu_N}(\chi_i) \}$ is locally around $\chi_i=[\iota^*_{F_g}\hol_i]$ a smooth,
half-dimensional submanifold of $X(\pi_1F_g,\SU(2))$.
Furthermore the cohomology classes of the cocycles $\{ z_{tws}^1+z_{len}^1, \ldots,
z_{tws}^N+z_{len}^N \}$ provide a basis for $T_{\chi_1}V_1$, while
the cohomology classes of the cocycles $\{ z_{tws}^1-z_{len}^1, \ldots, z_{tws}^N-z_{len}^N \}$ provide a basis for $T_{\chi_2}V_2$.
\end{lemma}

\subsection{Local rigidity}

\begin{lemma}\label{cohomology}
Let $C$ be a hyperbolic or a spherical cone-3-manifold with cone-angles
$\leq\pi$. Then:
\begin{enumerate}
\item The natural map $H^1(M,\E) \rightarrow H^1(\partial M_{\varepsilon},\E)$ is
injective.
\item $\dim H^1(M,\E)=\frac{1}{2}\dim H^1(\partial M_{\varepsilon},\E)$.
\end{enumerate}
In the spherical case, the assertions hold for the parallel subbundles $\E_i\subset\E$ individually.
\end{lemma}
\begin{pf}
Let us look at a part of the long exact cohomology sequence of the pair
$(M_{\varepsilon},\partial M_{\varepsilon})$ with coefficients in $\E$. The natural map
$q: H^1(M_{\varepsilon},\partial M_{\varepsilon},\E)\rightarrow H^1(M_{\varepsilon},\E)$ factors through $L^2$-cohomology,
since $H^1(M_{\varepsilon},\partial M_{\varepsilon},\E)=H^1_{cp}(M,\E)$:
$$\xymatrix{
 \ar[r] & H^1(M_{\varepsilon},\partial M_{\varepsilon},\E) \ar[r]^q & H^1(M_{\varepsilon},\E) \ar[r]^r & H^1(\partial
M_{\varepsilon},\E) \ar[r] & \\
 & H^1_{cp}(M,\E) \ar@{=}[u]\ar[r] & H^1_{L^2}(M,\E) \ar[u] & &}
$$
Since by our vanishing theorem $H^1_{L^2}(M,\E)=0$, we have that $q$ is the zero
map and $r:
H^1(M_{\varepsilon},\E) \rightarrow H^1(\partial M_{\varepsilon},\E)$ is injective.

Since the Killing form $B$ on $\E$ (resp. the parallel metric $h^{\E}$ in the spherical
case) provides a non-degenerate coefficient
pairing, we can apply Poincar\'{e} duality to conclude that $H^2(M_{\varepsilon},\partial
M_{\varepsilon},\E)\cong H^1(M_{\varepsilon},\E)^*$ and $H^2(M_{\varepsilon},\E)\cong H^1(M_{\varepsilon},\partial M_{\varepsilon},\E)^*$. The
Poincar\'{e} duality isomorphisms are natural, such that we obtain the following
commutative diagram:
$$\xymatrix{
 & H^1(M_{\varepsilon},\E)^* \ar[r]^{q^*} & H^1(M_{\varepsilon},\partial M_{\varepsilon},\E)^* & \\
 \ar[r] & H^2(M_{\varepsilon},\partial M_{\varepsilon},\E)
\ar[u]^{\cong}_{\text{P.D.}}\ar[r] & H^2(M_{\varepsilon},\E)
\ar[u]^{\cong}_{\text{P.D.}} \ar[r] & \\}
$$
Since $q^*=0$, we obtain the following short exact sequence:
$$\xymatrix{
 & & & H^1(M_{\varepsilon},\E)^* &\\
0 \ar[r] & H^1(M_{\varepsilon},\E) \ar[r] & H^1(\partial M_{\varepsilon}, \E) \ar[r] & H^2(M_{\varepsilon},\partial M_{\varepsilon},\E)
\ar[r]\ar[u]^{\cong}_{\text{P.D.}} & 0}
$$
This implies that $\dim H^1(M_{\varepsilon},\E)=\frac{1}{2}\dim H^1(\partial M_{\varepsilon},\E)$. In the spherical case these arguments apply to the parallel subbundles $\E_i\subset\E$.
\end{pf}
\subsubsection{The hyperbolic case}
The following is a well-known fact about the holonomy representation of a
hyperbolic cone-manifold structure, for convenience of the reader we give
a proof:
\begin{lemma}
The holonomy of a hyperbolic cone-manifold structure $\hol:\pi_1M\rightarrow \SL_2(\C)$ is irreducible.
\end{lemma}
\begin{pf}
Let us assume that the holonomy representation is reducible. Then there is a
point $x_{\infty} \in \partial_{\infty} \Hyp^3$ fixed by the holonomy. The
volume decreasing flow, which moves each point $x$ with unit speed towards
$x_{\infty}$ along the unique geodesic connecting $x$ and $x_{\infty}$, may
then be pulled back via the developing map to a volume decreasing flow on
$M$. This is a contradiction since $M$ has finite volume.
\end{pf}
\begin{lemma}
Let $\hol:\pi_1M\rightarrow \SL_2(\C)$ be the 
holonomy of a hyperbolic cone-manifold structure with cone-angles $\leq\pi$. Then $\hol$ is
a smooth point of $R(\pi_1M,\SL_2(\C))$. The $\C$-dimension of
$R(\pi_1M,\SL_2(\C))$ around $\hol$ equals $\tau+3-\frac{3}{2}\chi(\partial
M_{\varepsilon})$, where $\tau$ is the number of torus components contained in $\partial M_{\varepsilon}$. $T_{\hol}R(\pi_1M,\SL_2(\C))$ may be identified with $Z^1(\pi_1M,\sl_2(\C))$.
\end{lemma}
\begin{pf}
We follow the discussion in M. Kapovich's book (cf.~\cite{Kap}), which
essentially amounts to a
transversality argument. The key to the proof is the following splitting of $M_{\varepsilon}$:
\begin{lemma}{\rm\cite[Lm. 8.46]{Kap}}
There is a system of disjoint $1$-handles $\{ H_1,\ldots,H_t\}$ in $M_{\varepsilon}$
attached to $\partial M_{\varepsilon}$ such that $M_1:=M_{\varepsilon}\setminus
\interior(\cup_iH_i)$ is a handlebody.
\end{lemma}
As a consequence $M_{\varepsilon}$ may be written as a union
$$M_{\varepsilon}=M_1\cup_SM_2\,,$$ 
where $S$ is a surface of genus $g=1+t-\chi(\partial M_{\varepsilon})/2$. $M_2$ is homotopy
equivalent to the wedge product of the components of $\partial M_{\varepsilon}$ and $t-b+1$
circles, where $b$ is the number of components of $\partial M_{\varepsilon}$.
Therefore we obtain by van Kampen's theorem
$$
\pi_1M_{\varepsilon} = \pi_1M_1 \amalg_{\pi_1S} \pi_1M_2\,,
$$
where $\pi_1M_1$ is the free group on $g$ generators, and $\pi_1 M_2$ splits as a
free product of the fundamental groups of the components of $\partial M_{\varepsilon}$ and
$t-b+1$ $\Z$-factors. 
Consequently we obtain for the representation varieties
$$
R(\pi_1M_{\varepsilon},\SL_2(\C)) = R(\pi_1M_1, \SL_2(\C)) \times_{R(\pi_1S,\SL_2(\C))}
R(\pi_1M_2, \SL_2(\C))
$$
via the maps 
$$
\res_1 : R(\pi_1M_1,\SL_2(\C)) \rightarrow R(\pi_1S,\SL_2(\C))
$$
and
$$
\res_2 : R(\pi_1M_2,\SL_2(\C)) \rightarrow R(\pi_1S,\SL_2(\C))\,.
$$
$R(\pi_1M_1, \SL_2(\C))$ and $R(\pi_1M_2, \SL_2(\C))$ are smooth near the
restriction of the holonomy of a hyperbolic cone-manifold structure. Note that $\pi_1S$ surjects onto $\pi_1 M_{\varepsilon}$. Since $\hol$ is irreducible, this will also be the case for $\iota_S^*\hol$, which is therefore seen to be a smooth point of $R(\pi_1S,\SL_2(\C))$.

Therefore it is sufficient to show
that $\res_1$ and $\res_2$ meet transversally at $\iota^*_S \hol$. This will
follow from the equation
\begin{align*}
& \dim_{\C}Z^1(\pi_1M_1,\sl_2(\C)) +
\dim_{\C}Z^1(\pi_1M_2,\sl_2(\C)) \\
=& \dim_{\C}Z^1(\pi_1S,\sl_2(\C)+ \dim_{\C} Z^1(\pi_1M_{\varepsilon},\sl_2(\C))\,,
\end{align*}
if we use the identification
\begin{align*}
Z^1(\pi_1M_{\varepsilon},\sl_2(\C))&=\{(z_1,z_2)
\,\vert\, d\res_1(z_1)=d\res_2(z_2) \}\\
& \subset
Z^1(\pi_1M_1,\sl_2(\C))\oplus Z^1(\pi_1M_2,\sl_2(\C))\,.
\end{align*}
To obtain the desired equation, we have to calculate the dimensions of the
cocycle spaces. Note that $Z^1(\Gamma\amalg \Gamma',\g)=Z^1(\Gamma,\g)\oplus Z^1(\Gamma',\g)$. 
\begin{enumerate}
\item $\dim_{\C}Z^1(\pi_1M_1,\sl_2(\C))=3+3t-\frac{3}{2}\chi(\partial M_{\varepsilon})$, since
$\pi_1M_1$ is the free group on $g=1+t-\chi(\partial M_{\varepsilon})/2$
generators.  

\item $\dim_{\C}Z^1(\pi_1M_2,\sl_2(\C))=\tau-3\chi(\partial M_{\varepsilon})+3t+3$,
since we have $\dim_{\C}Z^1(\pi_1T^2,\sl_2(\C))=4$ at $\iota_{T^2}^*\hol$,
$\dim_{\C}Z^1(\pi_1F_g,\sl_2(\C))=-3\chi(F_g)+3$ at $\iota_{F_g}^*\hol$ and
the fundamental group of a wedge of $t-b+1$ circles is the free group on that
number of generators.

\item $\dim_{\C}Z^1(\pi_1S,\sl_2(\C))=6t-3\chi(\partial M_{\varepsilon})+3$, since 
$\iota^*_S\hol$ is irreducible.
\item  $\dim_{\C}Z^1(\pi_1M_{\varepsilon},\sl_2(\C))=\tau-\frac{3}{2}\chi(\partial M_{\varepsilon})+3$, since by Lemma \ref{cohomology}  
$\dim_{\C}H^1(M_{\varepsilon},\E)=\frac{1}{2}\dim_{\C}H^1(\partial
M_{\varepsilon},\E)$, furthermore $\hol$
is irreducible, therefore $Z^0(\pi_1M_{\varepsilon},\sl_2(\C))=0$.
\end{enumerate}
This finishes the proof.
\end{pf}
\begin{cor}
$\chi=[\hol]$ is a smooth point of $X(\pi_1M,\SL_2(\C)$. The $\C$-dimension of
$X(\pi_1M,\SL_2(\C)$ around $\chi$ equals  $\tau-\frac{3}{2}\chi(\partial
M_{\varepsilon})$, where $\tau$ is the number of torus components in $\partial
M_{\varepsilon}$. $T_{\chi}X(\pi_1M,\SL_2(\C))$ may be identified with $H^1(\pi_1M,\sl_2(\C))$.
\end{cor}
\begin{pf}
 $Z(\hol(\pi_1M))=\{\pm 1\}$ since $\hol$ is irreducible. Using 
Lemma \ref{Leeb} we proceed in the same way as in the surface case.
\end{pf}\\
\\
We are now ready to state and prove the main result in the hyperbolic case:
\begin{thm}
Let $C$ be a hyperbolic cone-3-manifold with cone-angles $\leq\pi$. Let $\{\mu_1, \ldots, \mu_N\}$ be the family of meridians, where $N$ is the number of edges contained in $\Sigma$. Then the map $$X(\pi_1M,\SL_2(\C))\rightarrow \C^N, \chi \mapsto
(t_{\mu_1}(\chi),\ldots,t_{\mu_N}(\chi))$$ is locally biholomorphic near $\chi=[\hol]$.
\end{thm}
\begin{pf}
Without loss of generality we may assume that $\Sigma$ is connected. Then we
have to consider two cases:
\begin{enumerate}
\item $\Sigma$ is a circle, i.e.~$\partial M_{\varepsilon} = T^2$
\item $\Sigma$ is a connected, trivalent graph, i.e.~$\partial M_{\varepsilon} = F_g$ 
\end{enumerate}
Let us recall what we already know. The level-set of the trace functions
$$V=\{t_{\mu_1}\equiv t_{\mu_1}(\chi),\ldots,t_{\mu_N}\equiv t_{\mu_N}(\chi)\}$$
is a smooth, half-dimensional submanifold of $X(\pi_1\partial M_{\varepsilon},\SL_2(\C))$
in each case, since the differentials
$
\{ dt_{\mu_1}, \ldots, dt_{\mu_N} \}
$
are $\C$-linearly independent in $H^1(\pi_1\partial
M_{\varepsilon},\SL_2(\C))^*$ at $\chi$.
If we work in the de-Rham realization of $H^1(\pi_1\partial
M_{\varepsilon},\SL_2(\C))$, the classes of the differential forms
$$
\{\omega_{len}^1, \ldots, \omega_{len}^N \}
$$
provide a $\C$-basis of $T_{\chi}V$. Furthermore, these forms are $L^2$-bounded on $U_{\varepsilon}(\Sigma)$.

On the other hand, the restriction map $H^1(M,\E) \rightarrow H^1(\partial
M_{\varepsilon},\E)$ is injective with half-dimensional image. This means that
$X(\pi_1M,SL_2(\C))$ is immersed into $X(\pi_1\partial M_{\varepsilon},SL_2(\C))$ as a
half-dimensional submanifold. 

We claim that the submanifolds $V$ and $X(\pi_1M,SL_2(\C))$ are transversal in $X(\pi_1\partial
M_{\varepsilon},SL_2(\C))$ at $\chi$.
It is sufficient to show that $T_{\chi}V$ and $\im(H^1(M,\E) \rightarrow
H^1(\partial M_{\varepsilon},\E))$ intersect trivially in
$H^1(\partial M_{\varepsilon},\E)$.

 Let $\omega \in \Omega^1(M,\E)$ be a closed form such that
$[\omega]\vert_{\partial M_{\varepsilon}} \in T_{\chi}V$. In particular, since
 the forms $\omega_{len}^i$ are $L^2$-bounded on $U_{\varepsilon}(\Sigma)$,
$\omega + d\sigma$ will be  $L^2$-bounded on $U_{\varepsilon}(\Sigma)$ for some $\sigma
\in \Gamma(U_{\varepsilon}(\Sigma),\E)$.  
We choose a cut-off function $\varphi$, which is $1$ in a neighbourhood of
$\Sigma$ and which is supported in $U_{\varepsilon}(\Sigma)$. Then $\varphi
\sigma$ extends to a section on $M$, such that $\omega+d(\varphi\sigma)$ is
$L^2$-bounded on $M$. Since $H^1_{L^2}(M,\E)=0$, this implies that $[\omega]=0$
in $H^1(M,\E)$ and therefore $[\omega]\vert_{\partial M_{\varepsilon}}=0$.   

It follows that the differentials
$
\{ dt_{\mu_1}, \ldots, dt_{\mu_N} \}
$
are $\C$-linearly independent already in $H^1(\pi_1M,\SL_2(\C))^*$.
\end{pf}\\
\\
The complex length of the $i$-th meridian is related to its
trace via 
$$
t_{\mu_i}(\rho)=\pm 2 \cosh(\mathcal{L}(\rho(\mu_i))/2)\,.
$$
Locally the set of representations $\rho: \pi_1M \rightarrow \SL_2(\C))$ such
that $\mathcal{L}(\rho(\mu_i))$ is purely imaginary for all $i\in\{1,\ldots,N\}$
corresponds to hyperbolic cone-manifold structures on $M$. The cone-angle
$\alpha_i$ is just the imaginary part of $\mathcal{L}(\rho(\mu_i))$, therefore
we obtain using Lemma \ref{complexlengthhyp}:
\begin{cor}[local rigidity] 
Let $C$ be a hyperbolic cone-3-mani-fold with cone-angles $\leq\pi$. 
Then the set of cone-angles $\{\alpha_1,\ldots,\alpha_N\}$,  where $N$ is the number of edges contained in $\Sigma$, provides a local parame-trization of the space of hyperbolic cone-manifold structures near the given structure on $M$. In particular, there are no deformations leaving the cone-angles fixed. 
\end{cor}
\subsubsection{The spherical case}
\begin{lemma}
\label{Porti}
Let $\hol:\pi_1M\rightarrow \SU(2)\times \SU(2)$ be the holonomy of a spherical cone-manifold structure. Then $\hol_1$ and $\hol_2$ are both non-abelian, unless $\Sigma$ is a link and $M$ is Seifert fibered.
\end{lemma}
\begin{pf}
Let us assume that $\hol_1$ is abelian. Then we may assume that the holonomy is contained in $S^1 \times \SU(2)$. This means that the Hopf-fibration on $\Sph^3 \subset \H$ obtained by left-multiplication with $S^1 \subset \H$ is preserved by the holonomy and may be pulled back via the developing map to a Seifert fibration on $M$. If $\hol_2$ is abelian, then the Hopf-fibration obtained by right-multiplication with $S^1 \subset \H$ will be invariant under the holonomy, and the same argument applies. In both cases the singular locus $\Sigma$ has to be a link, since in the presence of vertices $\hol_1$ and $\hol_2$ are clearly irreducible. 
\end{pf}
\begin{lemma}
Let $\hol_i:\pi_1M\rightarrow \SU(2)$ be a component of the holonomy of a
spherical cone-manifold structure with cone-angles $\leq\pi$. If $M$ is not Seifert fibered, then $\hol_i$ is a smooth point of $R(\pi_1M,\SU(2))$. The $\R$-dimension of
$R(\pi_1M,\SU(2))$ around $\hol_i$ equals $\tau+3-\frac{3}{2}\chi(\partial
M_{\varepsilon})$, where $\tau$ is the number of torus components in $\partial M_{\varepsilon}$. $T_{\hol_i}R(\pi_1M,\SU(2))$ may be identified with $Z^1(\pi_1M,\su(2))$.
\end{lemma}
\begin{pf}
The arguments in the hyperbolic case apply directly, 
the $\R$-dimensions of the $\su(2)$-cocycle spaces are equal to the $\C$-dimensions of the corresponding ${\sl}_2(\C)$-cocycle spaces.  
\end{pf}
\begin{cor}
$\chi_i=[\hol_i]$ is a smooth point of $X(\pi_1M,\SU(2))$. The $\R$-dimension of
$X(\pi_1M,\SU(2))$ around $\chi_i$ equals  $\tau-\frac{3}{2}\chi(\partial
M_{\varepsilon})$, where $\tau$ is the number of torus components in $\partial
M_{\varepsilon}$. $T_{\chi_i}X(\pi_1M,\SU(2))$ may be identified with $H^1(\pi_1M,\su(2))$.
\end{cor}
\begin{pf}
The action of $\SU(2)$ on $R(\pi_1M,\SU(2))$ is proper since $SU(2)$ is a compact group. Since $\hol_i$ is non-abelian by Lemma \ref{Porti}, we have that $Z(\hol_i(\pi_1M))=\{\pm 1\}$. Now the result follows as in the surface case.
\end{pf}\\ 
\\
The main result in the spherical case is the following theorem: 
\begin{thm}
Let $C$ be a spherical cone-3-manifold with cone-angles $\leq\pi$, which is not Seifert fibered. Let $\{\mu_i,\ldots,\mu_N\}$ be the family of meridians, where
$N$ is the number of edges contained in $\Sigma$.
Then the map $$X(\pi_1M,\SU(2))\rightarrow \R^N, \chi_i \mapsto
(t_{\mu_1}(\chi_i),\ldots,t_{\mu_N}(\chi_i))$$ is a local diffeomorphism near
$\chi_i=[\hol_i]$ for $i \in \{1,2\}$.
\end{thm}
\begin{pf}
The proof proceeds exactly along the same lines as in the hyperbolic case. The level-sets of the trace-functions
$$V_i=\{t_{\mu_1}\equiv t_{\mu_1}(\chi_i), \ldots,
t_{\mu_N}\equiv t_{\mu_N}(\chi_i) \}$$ are smooth, half-dimensional submanifolds of $X(\pi_1M_{\varepsilon},\SU(2))$ near $\chi_i$ for $i\in\{1,2\}$. 
The classes of the differential forms
$$
\{ \omega^1_{tws}+\omega^1_{len}, \ldots, \omega^N_{tws}+\omega^N_{len}\}
$$
provide a basis for $T_{\chi_1}V_1$, while the classes of the forms
$$
\{ \omega^1_{tws}-\omega^1_{len}, \ldots, \omega^N_{tws}-\omega^N_{len}\}
$$
provide a basis for $T_{\chi_1}V_1$. These forms are $L^2$-bounded on $U_{\varepsilon}(\Sigma)$. The same argument as in the hyperbolic case shows, that $T_{\chi_i}V_i$ and $\im(H^1(M,\E_i)\rightarrow H^1(\partial M_{\varepsilon},\E_i))$ are transversal for $i\in\{1,2\}$. It follows that the differentials
$
\{ dt_{\mu_1}, \ldots, dt_{\mu_N} \}
$
are $\R$-linearly independent already in $H^1(\pi_1M,\SU(2))^*$ at $\chi_i$ for $i\in\{1,2\}$. 
\end{pf}\\
\\
Locally around $\hol$ the set of representations $\rho=(\rho_1,\rho_2)$ such that
$t_{\mu_i}(\rho_1)=t_{\mu_i}(\rho_2)$, equivalently $\mathcal{L}_1(\rho(\mu_i))=0$, for all $i\in\{1,\ldots,N\}$ corresponds to
spherical cone-manifold structures on $M$. The cone-angle $\alpha_i=\mathcal{L}_2(\rho(\mu_i))$ is related to the
trace of the meridian via 
$$
t_{\mu_i}(\rho_1)=t_{\mu_i}(\rho_2)=\pm 2 \cos(\alpha_i/2)\,.
$$
Therefore we obtain using Lemma \ref{complexlengthsph}:
\begin{cor}[local rigidity]
Let $C$ be a spherical cone-3-manifold with cone-angles $\leq\pi$,
which is not Seifert fibered. 
Then the set of cone-angles $\{\alpha_1,\ldots,\alpha_N\}$,  where $N$ is the number of edges contained in $\Sigma$, provides a local parametrization of the space of spherical cone-manifold structures near the given structure on $M$. In particular, there are no deformations leaving the cone-angles fixed. 
\end{cor}


\end{document}